\documentclass[a4paper,10pt]{amsart}
\usepackage{latexsym}
\usepackage{amssymb}
\usepackage{mathrsfs}
\usepackage{amscd}
\usepackage{color}
\usepackage{epsfig}
\usepackage{graphicx}
\usepackage[all]{xy}

\theoremstyle{plain}
\newtheorem{thm}{Theorem}[section]
\newtheorem{cor}[thm]{Corollary}
\newtheorem{lem}[thm]{Lemma}
\newtheorem{prop}[thm]{Proposition}
\newtheorem{Thm}{Theorem}
\newtheorem{Cor}[Thm]{Corollary}

\theoremstyle{definition}
\newtheorem{defn}[thm]{Definition}
\newtheorem{rem}[thm]{Remark}
\newtheorem{ex}[thm]{Example}

\newcommand{\C}{{\mathbb C}}
\newcommand{\R}{{\mathbb R}}
\newcommand{\Q}{{\mathbb Q}}
\newcommand{\Z}{{\mathbb Z}}
\newcommand{\tabtopsp}[1]{\vbox{\vbox to#1{}\vbox to1em{}}}

\def\t{\mathfrak t}

\def\CA{\mathcal{A}}
\def\CB{\mathcal{B}}

\def\CI{\mathcal{I}}

\def\CL{\mathcal{L}}

\def\CN{\mathcal{N}}

\def\CS{\mathcal{S}}
\def\CT{\mathcal{T}}
\def\CU{\mathcal{U}}

\def\SH{\mathscr{H}}
\def\SK{\mathscr{K}}
\def\SS{\mathscr{S}}

\def\om{\omega}

\def\<{\left\langle}
\def\>{\right\rangle}

\def\Aut{\operatorname{Aut}}

\def\id{\operatorname{id}}
\def\SL{\operatorname{SL}}
\def\GL{\operatorname{GL}}
\def\pr{\operatorname{pr}}

\def\Vect{\operatorname{Vect}}

\newcommand{\abs}[1]{\lvert#1\rvert}
\newcommand{\norm}[1]{\lVert#1\rVert}

\numberwithin{equation}{section}

\begin{document}
\title[Local torus actions]{Local torus actions modeled on the standard representation}
\author[T. Yoshida]{Takahiko Yoshida} 


\address{Graduate School of Mathematical Sciences, The University of Tokyo, 8-1 Komaba 3-chome, Meguro-ku, Tokyo, 153-8914, Japan}
\email{takahiko@ms.u-tokyo.ac.jp}
\curraddr{Meiji Institute for Advanced Study of Mathematical Sciences, Meiji University, 1-1-1 Higashimita, Tama-ku, Kawasaki, 214-8571, Japan}
\email{takahiko@math.meiji.ac.jp}


\thanks{The author was partially supported by Research Fellowship of the Japan Society for the Promotion of Science
for Young Scientists.}

\begin{abstract}
We introduce the notion of a {\it local torus action modeled on the standard representation} (for simplicity, we call it a local torus action). It is a generalization of a locally standard torus action and also an underlying structure of a locally toric Lagrangian fibration. For a local torus action, we define two invariants called a {\it characteristic pair} and an {\it Euler class of the orbit map}, and prove that local torus actions are classified topologically by them. As a corollary, we obtain a topological classification of locally standard torus actions, which is a generalization of the topological classification of quasi-toric manifolds by Davis and Januszkiewicz~\cite[Proposition 1.8]{DJ} and of effective $T^2$-actions on four-dimensional manifolds without nontrivial finite stabilizers by Orlik and Raymond~\cite{OR}. We investigate locally toric Lagrangian fibrations from the viewpoint of local torus actions. We give a necessary and sufficient condition in order that a local torus action becomes a locally toric Lagrangian fibration. Locally toric Lagrangian fibrations are classified by Boucetta-Molino~\cite{BM} up to fiber-preserving symplectomorphisms. We shall reprove the classification theorem of locally toric Lagrangian fibrations by refining the proof of the classification theorem of local torus actions. We also investigate the topology of a manifold equipped with a local torus action when the Euler class of the orbit map vanishes. 
\end{abstract}
\maketitle

\section{Introduction}\label{introduction}
Let $S^1$ be the unit circle in $\C$ and $T^n:=(S^1)^n$ the $n$-dimensional compact torus. $T^n$ acts on the $n$-dimensional complex vector space $\C^n$ by coordinatewise complex multiplication. This action is called the {\it standard representation of $T^n$}. In this paper we focus on manifolds which are locally modeled on the standard representation of $T^n$. A typical example is a nonsingular toric variety. $T^n$ acts on an $n$-dimensional toric variety as a subgroup of the $n$-dimensional complex torus $(\C^*)^n$. If the toric variety is nonsingular, then for each point of it we can take a coordinate neighborhood $(U, \rho, \varphi)$, where $U$ is a $T^n$-invariant connected open neighborhood of the point, $\rho$ is an automorphism of $T^n$, and $\varphi$ is a $\rho$-equivariant diffeomorphism from $U$ to an open set of $\C^n$ invariant under the standard representation of $T^n$. The latter means that $\varphi(u\cdot x)=\rho(u)\cdot \varphi(x)$ for $u\in T^n$ and $x\in U$. In general, a torus action which has an atlas consisting of such coordinate neighborhoods is said to be {\it locally standard} and such an atlas is called a {\it standard atlas}. This is one of the starting point of the pioneer work~\cite{DJ} of Davis and Januszkiewicz, in which they focused on locally standard torus actions whose orbit spaces are simple convex polytopes and showed that they still have fascinating combinatorial properties of toric varieties. (A manifold with this torus action is now called a {\it quasi-toric manifold}.) After their work, topological generalizations of the original toric theory have been actively studied~\cite{BP, BR, H2, HM, HM2, Ma, MP}. 

In this paper, for a topological space which is not necessarily equipped with a global torus action we consider a generalization of a standard atlas. 
\begin{defn}\label{localaction}
Let $X$ be a compact Hausdorff space. A {\it weakly standard $C^r$ $(0\le r\le \infty )$ atlas} of $X$ is an atlas $\{ (U_{\alpha}^X, \varphi^X_{\alpha})\}_{\alpha \in \CA}$ which satisfies the following properties: 
\begin{enumerate}
\item for each $\alpha$, $\varphi^X_{\alpha}$ is a homeomorphism from $U_{\alpha}^X$ to an open set of $\C^n$ which is invariant with respect to the standard representation of $T^n$,  
\item on each nonempty overlap $U_{\alpha \beta}^X:=U_{\alpha}^X\cap U_{\beta}^X$, 
\begin{enumerate}
\item $\varphi^X_{\alpha}(U_{\alpha \beta}^X)$ and $\varphi^X_{\beta}(U_{\alpha \beta}^X)$ are also invariant with respect to the standard representation of $T^n$ and 
\item there exists an automorphism $\rho_{\alpha \beta}$ of $T^n$ such that the overlap map $\varphi^X_{\alpha \beta}:=\varphi^X_{\alpha}\circ (\varphi^X_{\beta})^{-1}$ is $\rho_{\alpha \beta}$-equivariant $C^r$ diffeomorphic with respect to the standard representation of $T^n$.\footnote{In this paper, for simplicity, we assume that every overlap of two open sets is connected.} 
\end{enumerate}
\end{enumerate}
Two weakly standard $C^r$ atlases $\{ (U_{\alpha}^X, \varphi^X_{\alpha})\}_{\alpha \in \CA}$ and $\{ (V_{\beta}^X, \psi^X_{\beta})\}_{\beta \in \CB}$ of $X$ are defined to be {\it equivalent} if on each nonempty overlap $U_{\alpha}^X\cap V_{\beta}^X$, there exists an automorphism $\rho$ of $T^n$ such that $\varphi^X_{\alpha}\circ (\psi^X_{\beta})^{-1}$ is $\rho$-equivariant $C^r$ diffeomorphic. We call an equivalence class of weakly standard $C^r$ atlases a {\it $C^r$ local $T^n$-action on $X$ modeled on the standard representation}, or a local $T^n$-action on $X$ if there are no confusions and denote it by $\CT$. 
\end{defn}

For each local torus action $\CT$ there uniquely exists a weakly standard atlas in $\CT$ that contains all weakly standard atlases in $\CT$. We call it a {\it maximal} weakly standard atlas in $\CT$.
\begin{defn}\label{locally-equiv-iso}
Let $(X_i,\CT_i)$ ($i=1,2$) be a $2n$-dimensional manifold equipped with a $C^r$ local $T^n$-action $\CT_i$ and $\{ (U_{\alpha}^{X_1}, \varphi^{X_1}_{\alpha})\}_{\alpha \in \CA} \in \CT_1$ and $\{ (U_{\beta}^{X_2}, \varphi^{X_2}_{\beta})\}_{\beta \in \CB}\in \CT_2$ the maximal weakly standard atlases of $X_1$ and $X_2$. A {\it $C^r$ isomorphism $f_X\colon (X_1,\CT_1)\to (X_2,\CT_2)$} between $(X_1,\CT_1)$ and $(X_2,\CT_2)$ is a $C^r$ diffeomorphism $f_X\colon X_1\to X_2$ which satisfies the following conditions, namely, for each nonempty overlap $U_{\alpha}^{X_1}\cap f_X^{-1}(U_{\beta}^{X_2})\neq \emptyset$
\begin{enumerate}
\item $\varphi^{X_1}_\alpha(U_{\alpha}^{X_1}\cap f_X^{-1}(U_{\beta}^{X_2}))$ and $\varphi^{X_2}_\beta(f_X(U_{\alpha}^{X_1})\cap U_{\beta}^{X_2})$ are invariant with respect to the standard representation of $T^n$ and 
\item there exists an automorphism $\rho$ of $T^n$ such that $\varphi^{X_2}_{\beta}\circ f_X\circ (\varphi^{X_1}_{\alpha})^{-1}$ is $\rho$-equivariant. 
\end{enumerate}
$(X_1,\CT_1)$ and $(X_2,\CT_2)$ are said to be {\it $C^r$ isomorphic} if there exists a $C^r$ isomorphism between $(X_1,\CT_1)$ and $(X_2,\CT_2)$. 
\end{defn}
The purpose of this paper is to classify local torus actions topologically in terms of certain invariants. This is an improvement on the previous work~\cite{Y1}. The main motivation of this work is to develop the equivariant theory of local torus actions and its application to the geometry and topology of Lagrangian fibrations. In \cite{Y3} we also discuss the lifting problem of local torus actions to fiber bundles. 

It is obvious that a standard atlas of a locally standard torus action satisfies the conditions in Definition~\ref{localaction}. In this sense, the notion of a local torus action is a generalization of that of a locally standard torus action. As in the usual group action case, we can define the orbit space and the orbit map for a local $T^n$-action $\CT$ on a manifold $X$. In fact, in Section~\ref{orbit-obst} we construct an $n$-dimensional $C^0$ manifold $B_X$ with corners and a $C^0$ open map $\mu_X \colon X\to B_X$ which is locally identified with the orbit map of the standard representation of $T^n$. 

Locally standard torus actions provide several examples of local torus actions but not all local torus actions come from locally standard torus actions. Let $\Aut (T^n)$ be the group of automorphisms of $T^n$ as a group. Then a nontrivial $T^n$-bundle on an $n$-dimensional closed manifold whose structure group is $\Aut (T^n)$ is equipped with a local $T^n$-action which is not induced by any locally standard $T^n$-action. In general, for a local $T^n$-action $\CT$ on a manifold $X$, a weakly standard atlas $\{ (U_{\alpha}^X, \varphi^X_{\alpha})\}_{\alpha \in \CA}\in \CT$ induces an atlas $\{ (U_{\alpha}^B,\varphi^B_{\alpha} )\}$ of $B_X$, and automorphisms $\rho_{\alpha \beta}$ in (ii) of Definition~\ref{localaction} form a \v{C}ech one-cocycle on $\{ U_{\alpha}^B\}_{\alpha \in \CA}$ with values in $\Aut (T^n)$, hence it defines a cohomology class in the first \v{C}ech cohomology of $B_X$ with coefficients in $\Aut (T^n)$. We show that it is the obstruction for a local $T^n$-action to come from a locally standard $T^n$-action (see Proposition~\ref{ob1}). 

To give a topological classification of local torus actions, let us first recall the topological classification of quasi-toric manifolds by Davis and Januszkiewicz. Let $X$ be a quasi-toric manifold with $T^n$-action. By definition, the quotient space $X/T^n$ is a simple convex polytope. Let $F_1,\ldots ,F_d$ be facets of $X/T^n$. For each facet $F_a$ there uniquely exists a circle subgroup $S^1_a$ of $T^n$ such that the preimage of $F_a$ by the natural projection is fixed by $S^1_a$. Let $\Lambda$ be the lattice of integral elements in the Lie algebra of $T^n$. We denote by $L_a$ the rank-one sublattice in $\Lambda$ spanned by the integral element which generates $S^1_a$. Then, by associating $L_a$ with $F_a$, we obtain the map from the set of all facets of the orbit space to the set of rank-one sublattices in $\Lambda$. This map is called a {\it characteristic function}. We remark that the characteristic function is identified with the unimodular rank-one sublattice bundle $\coprod_aF_a\times L_a$ in the trivial $\Lambda$-bundle $\coprod_aF_a\times \Lambda$. In this paper we use this identification. In \cite{DJ} for a quasi-toric manifold Davis and Januszkiewicz constructed a new quasi-toric manifold, called a canonical model, by using the characteristic function, and proved that a quasi-toric manifold is equivariant homeomorphic to the canonical model of it. Masuda and Panov generalized their technique to locally standard torus manifolds and proved that the same is true for a locally standard torus manifold provided that the second cohomology group of the orbit space vanishes. 

We shall classify local torus actions topologically by generalizing their method. For a local $T^n$-action $\CT$ on a manifold $X$ and a weakly standard atlas $\{ (U_{\alpha}^X, \varphi^X_{\alpha})\}_{\alpha \in \CA}\in \CT$, the \v{C}ech one-cocycle $\{ \rho_{\alpha \beta}\}$ corresponds to the principal Aut$(T^n)$-bundle which we denote by $\pi_{P_X}\colon P_X\to B_X$ on $B_X$. By Proposition~\ref{ob1}, when $(X, \CT )$ comes from a locally standard $T^n$-action, $P_X$ is the product bundle $P_X=B_X\times \Aut (T^n)$. Let $\pi_{\Lambda_X}\colon \Lambda_X\to B_X$ be the associated $\Lambda$-bundle of $P_X$ by the natural isomorphism between $\Aut (T^n)$ and the automorphism group $\GL (\Lambda)$ of $\Lambda$. Since $B_X$ is a manifold with corners, there is a natural stratification on $B_X$. We denote by ${\mathcal S}^{(n-1)}B_X$ the codimension one part of the natural stratification of $B_X$. For $(X, \CT )$, let $\pi_{T_X}\colon T_X\to B_X$ the associated $T^n$-bundle of $P_X$ by the natural action of $\Aut (T^n)$ on $T^n$. By construction, $T_X$ acts fiberwise on $X$ and this fiberwise action is locally identified with the standard representation of $T^n$. Hence a characteristic function is generalized to a rank one subbundle $\pi_{\CL_X}\colon \CL_X\to \CS^{(n-1)}B_X$ of the restriction of $\pi_{\Lambda_X}\colon \Lambda_X\to B_X$ to ${\mathcal S}^{(n-1)}B$. We call $\pi_{\CL_X}\colon \CL_X\to \CS^{(n-1)}B_X$ the {\it characteristic bundle} and also call the pair $(P_X, \CL_X)$ the {\it characteristic pair}. Similar to quasi-toric manifolds, we can construct a new $C^0$ manifold $X_{(P_X, \CL_X)}$ equipped with a local torus action by using $(P_X, \CL_X)$. We still call $X_{(P_X, \CL_X)}$ a canonical model. But unlike the above case, in general, $X_{(P_X, \CL_X)}$ is not $C^0$ isomorphic to $(X, \CT )$ even if $(X, \CT )$ comes from a locally standard torus action. So, in order to represent the difference between them, we define another topological invariant $e(X,\CT )$ called the {\it Euler class of $\mu_X$}. It is a natural generalization of an Euler class of a principal torus bundle. We also show that $e(X,\CT )$ is an obstruction in order that $\mu_X$ admits a $C^0$ section. With these preliminaries we can state a topological classification theorem of local torus actions (Theorem~\ref{classification}). 
\begin{Thm}[a topological classification of local torus actions]
Local torus actions are classified by the characteristic pairs and the Euler classes of the orbit maps up to $C^0$ isomorphisms. 
\end{Thm}
We focus on locally standard torus actions. Then we can also obtain the following corollary (Corollary~\ref{classification-locallystandard}). 
\begin{Cor}[a topological classification of locally standard torus actions]
Locally standard torus actions are classified by the characteristic bundles and the Euler classes of the orbit maps up to equivariant homeomorphisms. 
\end{Cor}
This is a generalization of the topological classification of quasi-toric manifolds by Davis and Januszkiewicz~\cite{DJ} and of effective $T^2$-actions on four-dimensional manifolds without nontrivial finite stabilizers by Orlik and Raymond~\cite{OR}. 

In case of $e(X,\CT )=0$, we shall investigate the topology of a manifold $(X,\CT )$ with a local torus action by using algebraic topological methods. We show that if the zero-dimensional part $\CS^{(0)}B_X$ of $B_X$ is nonempty, then the fundamental group of $X$ is isomorphic to the one of $B_X$ (Theorem~\ref{pi_1}). We give a way to compute cohomology groups and $K$-groups by using the Atiyah-Hirzebruch-Leray spectral sequence for $\mu_X$ (Theorem~\ref{H-spectral} and Theorem~\ref{K-spectral}). In the case where $X$ is four-dimensional and both $X$ and $B_X$ are oriented, we observe the signature of $X$. 

Another important class of local torus actions appears in Lagrangian fibrations. Let $\om_{\C^n}$ be the standard symplectic structure on $\C^n$, 
\begin{equation}\label{om_C}
\om_{\C^n}:=\frac{1}{2\pi \sqrt{-1}}\sum_{k=1}^ndz_k\wedge d\overline{z}_k.
\end{equation}
The standard representation of $T^n$ preserves $\om_{\C^n}$ and the map $\mu_{\C^n}\colon \C^n\to \R^n$ which is defined by 
\begin{equation}\label{std-moment}
\mu_{\C^n}(z):=(\abs{z_1}^2, \ldots ,\abs{z_n}^2)
\end{equation}
for $z=(z_1, \ldots ,z_n)\in \C^n$ is a moment map of the standard representation of $T^n$. Note that the image of $\mu_{\C^n}$ is the $n$-dimensional standard nonnegative cone 
\begin{equation}\label{positive-cone}
\R^n_{\ge 0}:=\{ \xi =(\xi_1, \ldots ,\xi_n)\in \R^n \colon \xi_i\ge 0\ i=1, \ldots ,n\} . 
\end{equation}
Let $(X,\om)$ be a $2n$-dimensional closed $C^{\infty}$ symplectic manifold and $B$ an $n$-dimensional $C^{\infty}$ manifold with corners. A $C^{\infty}$ map $\mu \colon (X, \om )\to B$ is called a {\it locally toric Lagrangian fibration} if it is locally identified with $\mu_{\C^n}\colon (\C^n,\om_{\C^n})\to \R^n_{\ge 0}$ (for the precise definition see Definition~\ref{def-ltlf} ). It is a natural generalization of a moment map of a symplectic toric manifold. In the case of $\partial B=\emptyset$, it is a nonsingular Lagrangian fibration. Conversely, by the Arnold-Liouville theorem~\cite{Ar}, a nonsingular Lagrangian fibration with closed connected fibers on a closed manifold is also such an example. We will see in Proposition~\ref{locally-toric-Lag} that for a locally toric Lagrangian fibration $\mu \colon (X, \om )\to B$, $X$ admits a $C^{\infty}$ local torus action. We also give a necessary and sufficient condition in order that a $C^{\infty}$ local torus action $(X,\CT )$ admits a symplectic structure $\om$ so that $\mu_X \colon (X, \om )\to B_X$ is a locally toric Lagrangian fibration (Theorem~\ref{om}). In particular, as is well known, if $\mu_X \colon (X, \om )\to B_X$ is a locally toric Lagrangian fibration, then $B_X$ is equipped with a rigid structure called an {\it integral affine structure} (see Definition~\ref{integral-affine-structure}, also consult \cite{Du} and \cite[Lemma 2.5]{S} for more details). We show in Section~\ref{symp} that there is a characteristic pair associated with an integral affine structure, and the canonical model constructed by this characteristic pair admits a symplectic structure so that the orbit map is a locally toric Lagrangian fibration. We shall refine the method used to prove Theorem~\ref{classification} to obtain the classification theorem of locally toric Lagrangian fibrations (Theorem~\ref{symplectic-classification}). This theorem has been obtained by Boucetta and Molino in \cite{BM}. 
%

This paper is organized as follows. In the next section we define the orbit space and the orbit map of a local torus action and investigate their properties. In Section~\ref{example} we see several examples. We give an obstruction in order that a local torus action is induced by a locally standard torus action. We also show that a locally toric Lagrangian fibration admits a $C^{\infty}$ local torus action. In Section~\ref{cp-cm}, we introduce a characteristic pair and construct a canonical model from a characteristic pair. We also explain how a characteristic pair is associated with a local torus action. In Section~\ref{orbit-map-section}, we define the Euler class of the orbit map and prove that the orbit map has a $C^0$ section if and only if it vanishes. Section~\ref{classification thm} is devoted to the topological classification. In Section~\ref{symp} we investigate locally toric Lagrangian fibrations from the viewpoint of local torus actions. We give a necessary and sufficient condition in order that a manifold equipped with a local torus action admits a symplectic structure that makes the orbit map a locally toric Lagrangian fibration. We also reprove the classification theorem of locally toric Lagrangian fibrations by Boucetta and Molino. Finally, in Section~\ref{topology}, we investigate the topology such as fundamental groups, cohomology groups, and $K$-groups. In Appendix A we describe some facts about $K$-theory of low dimensional CW complexes which are used in this section. We also observe the signature in the oriented, four-dimensional case. 

\subsection{Conventions}
Throughout this paper we employ the vector notation in order to represent elements of $\C^n$, namely, $z=(z_1, \ldots ,z_n) \in \C^n$. The similar notation is also used for $T^n=(S^1)^n$, $\R^n$, etc. In this paper, all manifolds, maps, and local $T^n$-actions are assumed to be of class $C^0$ unless otherwise stated. 
\\

{\bf Acknowledgment.} I am thankful to Professor Margaret Symington for informing me of her work about Lagrangian fibrations and torus actions and to Professor Hisaaki Endo for explaining me Meyer's signature cocycle. I am also grateful to Professor Mikio Furuta and Professor Mikiya Masuda for helpful suggestions and for constant encouragement. Professor Furuta is my mentor and this work is motivated by his suggestion. 

\section{Orbit structures}\label{orbit-obst}
\subsection{The standard representation}
In this subsection, we briefly summarize the fundamental facts about the standard representation of a torus, which is the local model in this paper. Recall that the standard representation of $T^n$ is the $T^n$-action on $\C^n$ defined by 
\begin{equation}\label{std-representation}
u\cdot z:=(u_1z_1, \ldots ,u_nz_n)
\end{equation}
for $u=(u_1, \ldots ,u_n)\in T^n$ and $z=(z_1,\ldots , z_n)\in \C^n$. The orbit space $\C^n/T^n$ of the action is equipped with the natural stratification whose $k$-dimensional strata consists of $k$-dimensional orbits. Let $\R^n_{\ge 0}$ be the standard $n$-dimensional nonnegative cone which is defined by \eqref{positive-cone}. It also has the natural stratification with respect to the number of components $\xi_i$ of the coordinates $(\xi_1\,\ldots ,\xi_n)$ which are equal to zero. 
We denote the $k$-dimensional part of $\R^n_{\ge 0}$ by $\CS^{(k)}\R^n_{\ge 0}$. Namely, $\xi \in \R^n_{\ge 0}$ is an element of $\CS^{(k)}\R^n_{\ge 0}$ if and only if the number of nonzero components $\xi_i$ is equal to $k$. Then the map $\mu_{\C^n}\colon \C^n\to \R^n_{\ge 0}$ defined by \eqref{std-moment} is invariant with respect to the standard representation of $T^n$ and induces a homeomorphism from $\C^n/T^n$ to $\R^n_{\ge 0}$ which preserves stratifications. It is easy to see the following proposition. 
\begin{prop}\label{observation}
For any $\xi \in \R^n_{\ge 0}$, all points of $\mu_{\C^n}^{-1}(\xi )$ have a common stabilizer. If $\xi $ lies in $\CS^{(k)}\R^n_{\ge 0}$, then it is an $(n-k)$-dimensional subtorus of $T^n$. In particular, for two points $z_1$, $z_2\in \C^n$ with $\mu_{\C^n}(z_1)=\mu_{\C^n}(z_2)$, there exists an element $u\in T^n$ such that $z_2=u\cdot z_1$ and $u$ is unique up to the common stabilizer of the $T^n$-action on $\mu_{\C^n}^{-1}(\mu_{\C^n}(z_1))$. 
\end{prop}

Let $\rho$ be an automorphism in $\Aut (T^n)$. Let $f\colon \C^n \to \C^n$ be a $\rho$-equivariant homeomorphism. Then, the following proposition is obvious. 
\begin{prop}\label{observation2} 
For any $z\in \C^n$, $\rho$ sends the stabilizer of $z$ isomorphically to the one of $f(z)$. 
\end{prop}

\subsection{Orbit spaces and orbit maps}
Let $(X,\CT )$ be a $2n$-dimensional manifold equipped with a $C^r$ local $T^n$-action. For $(X,\CT)$ we define the orbit space and the orbit map in the following way. Let $\{ (U_{\alpha}^X, \varphi^X_{\alpha})\}_{\alpha \in \CA}\in \CT$ be the maximal weakly standard atlas of $X$. We endow each quotient space $\varphi^X_{\alpha}(U_{\alpha}^X)/T^n$ with the quotient topology induced from the topology of $\varphi^X_{\alpha}(U_{\alpha}^X)$ by the natural projection $\pi_{\alpha} \colon \varphi^X_{\alpha}(U_{\alpha}^X)\to \varphi^X_{\alpha}(U_{\alpha}^X)/T^n$. By the property (ii) for each overlap $U_{\alpha \beta}^X$, $\varphi^X_{\alpha \beta}$ induces a homeomorphism from $\varphi^X_{\beta}(U_{\alpha \beta}^X)/T^n$ to $\varphi^X_{\alpha}(U_{\alpha \beta}^X)/T^n$. We define two elements $b_{\alpha}\in \varphi^X_{\alpha}(U_{\alpha}^X)/T^n$ and $b_{\beta}\in \varphi^X_{\beta}(U_{\beta}^X)/T^n$ to be equivalent, or $b_{\alpha}\sim_{orb}b_{\beta}$ if $b_{\alpha}\in \varphi^X_{\alpha}(U_{\alpha \beta}^X)/T^n$, $b_{\beta} \in \varphi^X_{\beta}(U_{\alpha \beta}^X)/T^n$ and the map induced by $\varphi^X_{\alpha \beta}$ sends $b_{\beta}$ to $b_{\alpha}$. It is an equivalence relation on the disjoint union $\coprod_{\alpha}\left( \varphi^X_{\alpha}(U_{\alpha}^X)/T^n\right)$. We define the {\it orbit space} $B_X$ of the local $T^n$-action to be the quotient space 
\[
B_X:=\coprod_{\alpha}\left( \varphi^X_{\alpha}(U_{\alpha}^X)/T^n\right)/\sim_{orb}
\]
together with the quotient topology. It is easy to see that $B_X$ is a Hausdorff space and $\{ \varphi^X_{\alpha}(U_{\alpha}^X)/T^n\}_{\alpha \in \CA}$ is an open covering of $B_X$. By the construction of $B_X$, the map $\coprod_{\alpha}\pi_{\alpha} \circ \varphi^X_{\alpha}\colon \coprod_{\alpha}U_{\alpha}^X\to \coprod_{\alpha}\left( \varphi^X_{\alpha}(U_{\alpha}^X)/T^n\right)$ induces a $C^0$ open map from $X$ to $B_X$. We call it the {\it orbit map} of the local $T^n$-action and denote it by $\mu_X \colon X\to B_X$. 
\begin{defn}[{\cite[Section 6]{Dav}}]\label{C^0mfdwithcorners}
Let $B$ be a Hausdorff space. A {\it structure of an $n$-dimensional $C^0$ manifold with corners} on $B$ is a system of coordinate neighborhoods modeled on open subsets of $\R^n_{\ge 0}$ so that overlap maps are homeomorphisms which preserve the natural stratifications induced from the one of $\R^n_{\ge 0}$. 
\end{defn}
\begin{prop}\label{B}
$B_X$ is endowed with a structure of an $n$-dimensional $C^0$ manifold with corners. 
\end{prop}
\begin{proof}
A structure of an $n$-dimensional $C^0$ manifold with corners on $B_X$ is constructed as follows. We put $U_{\alpha}^B:=\varphi^X_{\alpha}(U_{\alpha}^X)/T^n$. The restriction of $\mu_{\C^n}$ to $\varphi^X_{\alpha}(U_{\alpha}^X)$ induces a homeomorphism from $U_{\alpha}^B$ to the open subset $\mu_{\C^n}(\varphi^X_{\alpha}(U_{\alpha}^X))$ of $\R^n_{\ge 0}$, which is denoted by $\varphi^B_{\alpha}$. By construction, on each overlap $U_{\alpha \beta}^B:=U_{\alpha}^B\cap U_{\beta}^B$, the overlap map $\varphi^B_{\alpha \beta}:=\varphi^B_{\alpha}\circ (\varphi^B_{\beta})^{-1}\colon \mu_{\C^n}(\varphi^X_{\beta}(U_{\alpha \beta}^X))\to \mu_{\C^n}(\varphi^X_{\alpha}(U_{\alpha \beta}^X))$ preserves the natural stratifications of $\mu_{\C^n}(\varphi^X_{\alpha}(U_{\alpha \beta}^X))$ and $\mu_{\C^n}(\varphi^X_{\beta}(U_{\alpha \beta}^X))$. Thus, $\{ (U_{\alpha}^B,\varphi^B_{\alpha})\}_{\alpha \in \CA}$ is the desired atlas. 
\end{proof}
\begin{rem}\label{std}
The atlas $\{ (U_{\alpha}^B,\varphi^B_{\alpha})\}_{\alpha \in \CA}$ of $B_X$ constructed in the proof of Proposition~\ref{B} has the following property:  
for each $\alpha$, $U_{\alpha}^X=\mu_X^{-1}(U_{\alpha}^B)$, 
$\varphi^X_{\alpha}(U_{\alpha}^X)=\mu_{\C^n}^{-1}(\varphi^B_{\alpha}(U_{\alpha}^B))$ 
and the following diagram commutes
\[
\xymatrix{
X\ar[d]^{\mu_X}\ar@{}[r]|{\supset}& \mu_X^{-1}(U_{\alpha}^B)\ar[r]^{\varphi^X_{\alpha}}\ar[d]^{\mu_X} & 
\mu_{\C^n}^{-1}(\varphi^B_{\alpha}(U_{\alpha}^B))\ar[d]^{\mu_{\C^n}}\ar@{}[r]|{\subset}& 
\C^n\ar[d]^{\mu_{\C^n}} \\
B_X\ar@{}[r]|\supset & U_{\alpha}^B\ar[r]^{\varphi^B_{\alpha}} & \varphi^B_{\alpha}(U_{\alpha}^B)\ar@{}[r]|\subset& 
\R^n_{\ge 0}. 
}
\]
\end{rem}
\begin{rem}
Let $(X_1,\CT_1)$ and $(X_2,\CT_2)$ be $2n$-dimensional manifolds equipped with $C^r$ local $T^n$-actions. Suppose that there is a $C^r$ isomorphism $f_X\colon (X_1,\CT_1)\to (X_2,\CT_2)$. Then, $f_X$ induces a stratification preserving homeomorphism $f_B\colon B_{X_1}\to B_{X_2}$ such that $f_X$ and $f_B$ satisfy $\mu_{X_2}\circ f_X=f_B\circ \mu_{X_1}$. In Section~\ref{classification thm}, we classify local torus actions up to $C^0$ isomorphisms. 
\end{rem}

Let $(X, \CT )$ be a $2n$-dimensional manifold equipped with a $C^r$ local $T^n$-action. By Proposition~\ref{B} $B_X$ is equipped with a natural stratification. We denote by $\CS^{(k)}B_X$ the $k$-dimensional part of $B_X$, namely, $\CS^{(k)}B_X$ consists of those points which are sent to points in $\CS^{(k)}\R^n_{\ge 0}$ by a local coordinate system. In particular, the top-dimensional part $\CS^{(n)}B_X$ is equal to the interior of $B_X$ (see Figure~1). 
\begin{figure}[hbtp]
\begin{center}
\input{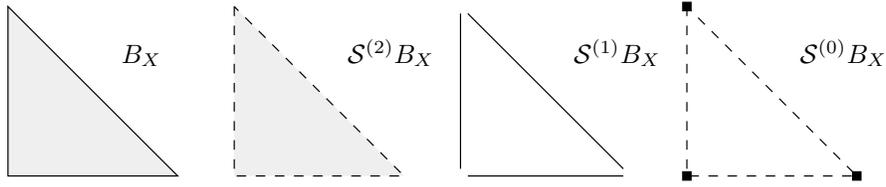}
\caption{The natural stratification when $B_X$ is a triangle}
\label{fig16}
\end{center}
\end{figure}

\begin{prop}\label{T^n-bundle} 
For each $k$, the restriction $\mu_X |_{\CS^{(k)}B_X} \colon \mu_X^{-1}(\CS^{(k)}B_X) \to \CS^{(k)}B_X$ of $\mu_X$ to $\CS^{(k)}B_X$ admits a structure of a $C^r$ fiber bundle whose fiber is $T^k$. In particular, $\mu_X |_{\CS^{(n)}B_X} \colon \mu_X^{-1}(\CS^{(n)}B_X) \to \CS^{(n)}B_X$ is a $C^r$ fiber bundle whose fiber is $T^n$ and whose structure group is contained in the semidirect product $T^n\rtimes \Aut (T^n)$ of $T^n$ and $\Aut (T^n)$. 
\end{prop} 
\begin{proof} 
It is sufficient to prove the proposition connected componentwise. Let $(\CS^{(k)}B_X)_a$ ($a=1$, $\ldots$, $m$) be the connected component of $\CS^{(k)}B_X$. We take a weakly standard atlas $\{ (U^X_{\alpha}, \varphi^X_{\alpha})\}_{\alpha \in \CA}\in \CT$. It induces an atlas $\{ (U^B_{\alpha}, \varphi^B_{\alpha})\}_{\alpha \in \CA}$ of $B_X$ which satisfies the property in Remark~\ref{std}. Let $(U_{\alpha}^B,\varphi^B_{\alpha})$ be a coordinate neighborhood with $(\CS^{(k)}B_X)_a\cap U_{\alpha}^B\neq \emptyset$. Then $(\CS^{(k)}B_X)_a\cap U_{\alpha}^B$ is represented as common zeros of some exactly $n-k$ components, say $\varphi^B_{\alpha ,i_1}$, $\cdots$, $\varphi^B_{\alpha ,i_{n-k}}$, of $\varphi^B_{\alpha}$, namely, 
\[
(\CS^{(k)}B_X)_a\cap U_{\alpha}^B=\{ b\in U_{\alpha}^B\colon \varphi^B_{\alpha ,i_1}(b)=\cdots =\varphi^B_{\alpha ,i_{n-k}}(b)=0\} .
\]
We put 
\[
T_{i_1, \ldots ,i_{n-k}}:=\{ u\in T^n\colon u_j=1\ \text{unless }j\neq i_1,\ldots ,i_{n-k}\} .
\]
Note that $T^n/T_{i_1, \ldots ,i_{n-k}}$ is identified with $T^k$. Let $\iota\colon \R^n_{\ge 0}\to \C^n$ be the section of $\mu_{\C^n}\colon\C^n \to \R^n_{\ge 0}$ which is defined by 
\begin{equation}\label{iota}
\iota (\xi):=(\xi_1^{1/2}, \ldots ,\xi_n^{1/2}). 
\end{equation}
By Proposition~\ref{observation} the equation 
\[
\varphi^X_{\alpha}(x)=u_{\alpha}(x)\cdot \iota \circ \varphi^B_{\alpha}\circ \mu_X (x) 
\]
uniquely defines the map $u_{\alpha}\colon \mu_X^{-1}((\CS^{(k)}B_X)_a\cap U_{\alpha}^B)\to T^n/T_{i_1, \ldots ,i_{n-k}}$. We define the local trivialization $\phi_{\alpha} \colon \mu_X^{-1}((\CS^{(k)}B_X)_a\cap U_{\alpha}^B)\to (\CS^{(k)}B_X)_a\cap U_{\alpha}^B\times T^k$ by 
\begin{equation}\label{local-triv}
\phi_{\alpha} (x):=\left( \mu_X(x), u_{\alpha}(x)\right) , 
\end{equation}
where we used the identification $T^n/T_{i_1,\ldots i_{n-k}}\cong T^k$. This proves the first half. Next we prove the other part. For $k=n$, a direct computation shows that on a nonempty overlap $U_{\alpha \beta}^B\setminus \partial B_X$, 
\[
\phi_{\alpha}\circ \phi_{\beta}^{-1}(b,u)
=(b, \rho_{\alpha \beta}(u)u_{\alpha \beta}(b)) 
\]
for $(b,u)\in U_{\alpha \beta}^B\setminus \partial B_X\times T^n$, where $u_{\alpha \beta}$ is the map $u_{\alpha \beta}\colon U_{\alpha \beta}^B\setminus \partial B_X\to T^n$ which is defined uniquely by  
\begin{equation}\label{T^n-part}
\varphi^X_{\alpha \beta}  \circ \iota \circ \varphi^B_{\beta}(b)=u_{\alpha \beta}(b) \cdot \iota \circ \varphi^B_{\alpha}(b). 
\end{equation}
This implies that the transition function with respect to $\phi_{\alpha}$ and $\phi_{\beta}$, which is denoted by $g_{\alpha \beta}$, is of the form $g_{\alpha \beta}(b)=(u_{\alpha \beta}(b), \rho_{\alpha \beta})\in T^n\rtimes \Aut (T^n)$. 
\end{proof}

\section{Examples}\label{example}
\subsection{Locally standard torus actions}
\begin{defn}[Locally standard torus actions]
Let $X$ be a $2n$-dimensional $C^{\infty}$ manifold equipped with a $C^{\infty}$ $T^n$-action. A {\it standard coordinate neighborhood} of $X$ is a triple $(U, \rho , \varphi)$ consisting of a $T^n$-invariant open set $U$ of $X$, an automorphism $\rho$ of $T^n$, and a $\rho$-equivariant diffeomorphism $\varphi$ from $U$ to some $T^n$-invariant open subset in $\C^n$. The action of $T^n$ on $X$ is said to be {\it locally standard} if every point in $X$ lies in some standard coordinate neighborhood, and an atlas which consists of standard coordinate neighborhoods is called a {\it standard atlas}. 
\end{defn}
See \cite{DJ, BP} for more details. A standard atlas is weakly standard. Therefore, a locally standard $T^n$-action on a closed $C^{\infty}$ manifold $X$ induces a local $T^n$-action on $X$. We give examples of locally standard torus actions. 
\begin{ex}
An effective $C^{\infty}$ $T^2$-action on a $4$-dimensional $C^{\infty}$ manifold $X$ without nontrivial finite stabilizers is locally standard because of the slice theorem. These actions have been studied by Orlik and Raymond in \cite{OR}. 
\end{ex}
\begin{ex}[Nonsingular toric varieties]
A complex $n$-dimensional toric variety is a normal complex algebraic variety $X$ of dimension $n$ with a $(\C^*)^n$-action having a dense orbit. $T^n$ acts on $X$ as a subgroup of $(\C^*)^n$. If $X$ is nonsingular, the $T^n$-action on $X$ is locally standard. In fact, it is well known that there is a one-to-one correspondence between toric varieties and fans, and top-dimensional cones in the fan associated with $X$ correspond to standard coordinate neighborhoods all of which cover $X$ since all cones are nonsingular. For toric varieties, see \cite{Dan, Ful, Oda}.  
\end{ex}
\begin{ex}[Quasi-toric manifolds]
A quasi-toric manifold is a $C^{\infty}$ manifold equipped with a locally standard torus action whose orbit space is combinatorially isomorphic to a simple convex polytope. A quasi-toric manifold was first introduced by Davis and Januszkiewicz in their pioneer work~\cite{DJ} as a topological generalization of a projective toric variety and now it plays a central role in toric topology. See \cite{DJ, BP} for more details. 
\end{ex}

Note that the existence of a local $T^n$-action does not necessarily imply an existence of a locally standard $T^n$-action. For example, a nontrivial fiber bundle on an $n$-dimensional closed manifold whose fiber is $T^n$ and whose structure group is $\Aut (T^n)$ is equipped with a local $T^n$-action which is not induced by any locally standard $T^n$-action. In general, for any $C^r$ local $T^n$-action $\CT$ on $X$, we take a weakly standard atlas $\{ (U_{\alpha}^X, \varphi^X_{\alpha})\}_{\alpha \in \CA}\in \CT$. We denote by $\{ (U_{\alpha}^B, \varphi^B_{\alpha})\}_{\alpha \in \CA}$ the atlas of $B_X$ which is induced by $\{ (U_{\alpha}^X, \varphi^X_{\alpha})\}_{\alpha \in \CA}$ and satisfies the property in Remark~\ref{std}. Assume that the index set $\CA$ is countably ordered. It is easy to see that the automorphisms $\rho_{\alpha \beta}$ of $T^n$ in the property ($\mathrm{ii}$) of Definition~\ref{localaction} form a \v{C}ech one-cocycle $\{ \rho_{\alpha \beta}\}$ on $\CU:=\{ U_{\alpha}^B\}_{\alpha \in \CA}$ with values in $\Aut (T^n)$. By definition, two one-cocycles $\{ \rho_{\alpha \beta}\}$ and $\{ \rho_{\alpha \beta}'\}$ on $\CU$ are equivalent if and only if there is a \v{C}ech zero-cochain $\{ \rho_{\alpha}\}$ such that 
\[
\rho_{\alpha \beta}'=\rho_{\alpha}\circ \rho_{\alpha \beta}\circ \rho_{\beta}^{-1}
\]
on each overlap $U_{\alpha \beta}^B$. We denote the set of equivalence classes by $H^1(\CU ;\Aut (T^n))$ and also denote its direct limit $\displaystyle \lim_{\to}H^1(\CU ;\Aut (T^n))$ with respect to refinements of $\CU$ by $H^1(B_X;\Aut (T^n))$. See \cite[Appendix A]{LM} for more details. $\{ \rho_{\alpha \beta}\}$ determines an equivalence class in $H^1(B_X;\Aut (T^n))$. It does not depend on the choice of equivalent weakly standard atlases and depends only on the local $T^n$-action. The equivalence class of $\{\rho_{\beta \alpha}\}$ in $H^1(B_X;\Aut (T^n))$ is an obstruction for the local $T^n$-action to be induced by a locally standard $T^n$-action. 
\begin{prop}\label{ob1}
A $C^r$ local $T^n$-action $\CT$ on $X$ is induced by some $C^r$ locally standard $T^n$-action if and only if $\{ \rho_{\alpha \beta}\}$ is equivalent to the trivial \v{C}ech one-cocycle in $H^1(B_X;\Aut (T^n))$, where the trivial \v{C}ech one-cocycle is the one whose values on all open set are equal to the identity map of $T^n$.  
\end{prop}
\begin{proof}
If a $C^r$ local $T^n$-action $\CT$ on $X$ is induced by some $C^r$ locally standard $T^n$-action, then there exists a standard atlas $\{ (U_{\alpha}^X,\rho_{\alpha},\varphi_{\alpha}^X)\}$ in $\CT$ such that $\rho_{\alpha \beta}$ is of the form $\rho_{\alpha \beta}=\rho_{\alpha}\circ \rho_{\beta}^{-1}$ on each overlap $U_{\alpha \beta}^B$. This implies that $\{ \rho_{\alpha \beta}\}$ is equivalent to the trivial \v{C}ech one-cocycle. 

Conversely, assume that $\{ \rho_{\alpha \beta}\}$ is equivalent to the trivial \v{C}ech one-cocycle in 
$H^1(B_X;\Aut (T^n))$. By replacing the open covering of $B_X$ by a refinement if necessary, we 
may assume that there exists a \v{C}ech zero-cochain $\{ \rho_{\alpha}\}$ on $\{U_{\alpha}^B\}$ such that $\rho_{\alpha}\circ \rho_{\beta}^{-1}=\rho_{\alpha \beta}$ for each $U_{\alpha \beta}^B$. Then we can recover the $T^n$-action globally on $X$ as 
\[
u\cdot x:=(\varphi^X_{\alpha})^{-1}(\rho_{\alpha}(u)\cdot \varphi^X_{\alpha}(x))
\]
for $u\in T^n$ and $x\in X$ provided that $x$ lies in $\mu_X^{-1}(U_{\alpha}^B)$. The action is well defined because if $x$ also lies in $\mu_X^{-1}(U_{\beta}^B)$,   
\[
\begin{split}
(\varphi^X_{\alpha})^{-1}(\rho_{\alpha}(u)\cdot \varphi^X_{\alpha}(x))&=
(\varphi^X_{\alpha})^{-1}(\rho_{\alpha}\circ \rho_{\beta}^{-1}\circ \rho_{\beta}(u)\cdot 
\varphi^X_{\alpha}\circ (\varphi^X_{\beta})^{-1}\circ \varphi^X_{\beta}(x)) \\
&=(\varphi^X_{\alpha})^{-1}(\rho_{\alpha \beta}(\rho_{\beta}(u))\cdot 
\varphi^X_{\alpha \beta}(\varphi^X_{\beta}(x))) \\
&=(\varphi^X_{\alpha})^{-1}(\varphi^X_{\alpha \beta}(\rho_{\beta}(u)\cdot \varphi^X_{\beta}(x))) \\
&=(\varphi^X_{\beta})^{-1}(\rho_{\beta}(u)\cdot \varphi^X_{\beta}(x)), 
\end{split}
\]
and for $u_1$, $u_2\in T^n$, 
\[
\begin{split}
u_1\cdot (u_2\cdot x)&=(\varphi^X_{\alpha})^{-1}
(\rho_{\alpha}(u_1)\cdot 
\varphi^X_{\alpha}((\varphi^X_{\alpha})^{-1}(\rho_{\alpha}(u_2)\cdot \varphi^X_{\alpha}(x)))) \\
&=(\varphi^X_{\alpha})^{-1}(\rho_{\alpha}(u_1u_2)\cdot \varphi^X_{\alpha}(x))\\
&=(u_1u_2)\cdot x. 
\end{split}
\]
\end{proof}

We can construct an example of a local torus action which is not induced by any global torus action. 
\begin{ex}\label{ex3.9}
Let $B$ be a two-dimensional $C^{\infty}$ torus with one boundary component and one corner points 
(see Figure~\ref{fig22}). 
\begin{figure}[hbtp]
\begin{center}
\input{fig22.pstex_t}
\caption{$B$}
\label{fig22}
\end{center}
\end{figure}
We construct a four-dimensional $C^{\infty}$ manifold $X$ equipped with a $C^{\infty}$ local $T^2$-action whose orbit space is $B$ as follows. First we focus on the interior $B\setminus \partial B$ of $B$ which is denoted by $B_1$. Let $\rho \colon\pi_1(B)(\cong \pi_1(B_1))\to SL_2(\Z )$ be the representation of the fundamental group which is defined by 
\[
\rho ([\alpha ])=
\begin{pmatrix}
1 & 0 \\
-1 & 1
\end{pmatrix},\ 
\rho ([\beta ])=
\begin{pmatrix}
1 & -1 \\
0 & 1
\end{pmatrix},\ 
\rho ([\gamma ])=
\begin{pmatrix}
3 & 1 \\
-1 & 0
\end{pmatrix},
\]
where $\alpha$, $\beta$, and $\gamma$ are representatives of generators of $\pi_1(B)$ which satisfy 
the relation $[\alpha][\beta][\alpha]^{-1}[\beta]^{-1}[\gamma]=1$ (see Figure~\ref{fig22}).  
We identify $T^2$ with $\R^2/\Z^2$ and denote by $\pi_T\colon T_{\rho}^2\to B_1$ the associated $T^2$-bundle 
of the universal covering of $B_1$ by $\rho$. 

Next we pay attention to a neighborhood of the boundary $\partial B$ of $B$. 
We define subsets $\overline{B}_2$, $U_1$, and $U_2$ of $\R^2_+$ by 
\[
\begin{split}
\overline{B}_2&=\{ \xi \in \R^2 \colon 0\le \xi_1 <4,\ 0\le \xi_2<1\} 
\cup \{ \xi \in \R^2 \colon 0\le \xi_1<1,\ 0\le \xi_2<4\} ,\\
U_1&=\{ \xi \in \R^2 \colon 3<\xi_1<4,\ 0\le \xi_2<1\} ,\\ 
U_2&=\{ \xi \in \R^2 \colon 0\le \xi_1<1,\ 3<\xi_2<4\} ,  
\end{split}
\]
and also define diffeomorphisms $\varphi^B \colon U_1\to U_2$ and 
$\varphi^X \colon\mu_{\C^2}^{-1}(U_1)\to \mu_{\C^2}^{-1}(U_2)$ by 
\[
\begin{split}
& \varphi^B (\xi )=(\xi_2, 7-\xi_1),\\
& \varphi^X (z)=\left( 
\dfrac{z_1^3z_2}{\abs{z_1}^3}, \sqrt{7-\abs{z_1}^2}\left( \dfrac{z_1}{\abs{z_1}}\right)^{-1}\right) . 
\end{split}
\]
Note that $\varphi^B$ and $\varphi^X$ satisfy 
\begin{equation}\label{commute}
\mu_{\C^2}\circ \varphi^X =\varphi^B\circ \mu_{\C^2} .
\end{equation}
We denote by $X_2$ the manifold which is obtained from $\mu_{\C^2}^{-1}(\overline{B}_2)$ by gluing $\mu_{\C^2}^{-1}(U_1)$ and $\mu_{\C^2}^{-1}(U_2)$ with $\varphi^X$ and denote by $B_2$ the surface with one corner 
which is obtained from $\overline{B}_2$ by gluing $U_1$ and $U_2$ with $\varphi^B$ (see Figure~3). 
\begin{figure}[hbtp]
\begin{center}
\input{fig4.pstex_t}
\caption{$B_2$}
\label{fig4}
\end{center}
\end{figure}
$B_2$ can be identified with a neighborhood of $\partial B$. By \eqref{commute}, $\mu_{\C^2}$ descends to the map from $X_2$ to $B_2$. We denote it by $\mu_2\colon X_2\to B_2$. It is easy to see that the restriction $\mu_2 |_{\mu_2^{-1}(B_1\cap B_2)}\colon\mu_2^{-1}(B_1\cap B_2)\to B_1\cap B_2$ is a fiber bundle with fiber $T^2$ and structure group $SL_2(\Z )$ which is isomorphic to $\pi_T |_{B_1\cap B_2}\colon T_{\rho}^2 |_{B_1\cap B_2}\to B_1\cap B_2$ because they have the same monodromy. Thus we can glue $\pi_T\colon T_{\rho}^2\to B_1$ and $\mu_2 \colon X_2\to B_2$ together to get the map $\mu \colon X\to B$. By the construction, $X$ is a $C^{\infty}$ four-dimensional manifold equipped with a local $T^2$-action whose orbit space $B_X$ and orbit map $\mu_X$ are equal to $B$ and $\mu$, respectively. Note that $\mu$ is $C^{\infty}$ in this example. 
\end{ex}

\subsection{Locally toric Lagrangian fibrations}
Let $(X,\om)$ be a $2n$-dimensional closed $C^{\infty}$ symplectic manifold and $B$ an $n$-dimensional $C^{\infty}$ manifold with corners. 
\begin{defn}[Locally toric Lagrangian fibrations~{\cite[Definition 2.7]{Ham}}]\label{def-ltlf}
A $C^{\infty}$ map $\mu \colon (X, \om )\to B$ is called a {\it locally toric Lagrangian fibration} if there exists a system $\{ (U_{\alpha},\varphi^B_{\alpha})\}$ of coordinate neighborhoods of $B$ modeled on $\R^n_{\ge 0}$, and for each $\alpha$ there exists a symplectomorphism $\varphi^X_{\alpha}\colon (\mu^{-1}(U_{\alpha}), \om )\to (\mu_{\C^n}^{-1}(\varphi^B_{\alpha}(U_{\alpha})), \om_{\C^n})$ such that $\mu_{\C^n}\circ \varphi^X_{\alpha}=\varphi^B_{\alpha}\circ \mu$. 
\end{defn}
\begin{ex}
A moment map of a symplectic toric manifold is a locally toric Lagrangian fibration.  
\end{ex}
\begin{ex}
Let $\om_{\R^n\times T^n}$ be the standard symplectic structure on $\R^n\times T^n$ which is defined by 
\begin{equation}\label{om_cotangent}
\om_{\R^n\times T^n}=\sum_{k=1}^nd\theta_k\wedge d\xi_k ,  
\end{equation}
where $(\xi_1,\ldots ,\xi_n)$ is the standard coordinates of $\R^n$ and $(\theta_1,\ldots ,\theta_n)$ is the angle coordinates of $T^n$ with period $1$, which means $(e^{2\pi \theta_1}, \ldots ,e^{2\pi \theta_n})\in T^n$. Then the projection $\pr_1\colon (\R^n\times T^n, \om_{\R^n\times T^n})\to \R^n$ to the first factor is a nonsingular Lagrangian fibration. It is easy to see that $\pr_1\colon (\R^n\times T^n, \om_{\R^n\times T^n})\to \R^n$ is a noncompact locally toric Lagrangian fibration.  
\end{ex}

The Arnold-Liouville theorem~\cite{Ar} says that $\pr_1\colon (\R^n\times T^n, \om_{\R^n\times T^n})\to \R^n$ is the local model of a nonsingular Lagrangian fibration. 
\begin{thm}[Arnold-Liouville~\cite{Ar}]
Let $\mu \colon (X,\om )\to B$ be a nonsingular Lagrangian fibration with closed connected fibers on a closed manifold. Then for each $b\in B$ there exists a coordinate neighborhood $(U,\varphi^B)$ and there also exists a symplectomorphism $\varphi^X\colon (\mu^{-1}(U),\om) \to (\varphi^B(U)\times T^n,\om_{\R^n\times T^n})$ such that $\pr_1\circ \varphi^X=\varphi^B\circ \mu$. 
\end{thm}
Therefore, a nonsingular Lagrangian fibration with closed connected fibers on a closed manifold is also a locally toric Lagrangian fibration. Nonsingular Lagrangian fibrations on closed oriented surfaces have been investigated by Mishachev. In \cite{Mis} he showed that the complete list of nonsingular Lagrangian fibrations on closed oriented surfaces. According to it, if a closed oriented surface $B$ is a base of a nonsingular Lagrangian fibration, then $B$ is $T^2$. Moreover, concrete constructions of all these have been given in \cite{Ge}. 

In the following example we give a locally toric Lagrangian fibration with singular fibers which is not a symplectic toric manifold.
\begin{ex}
For a sufficiently small positive number $\varepsilon$ ($1\gg \varepsilon >0$), let $U$ be the open subset of $(\C^2 ,\om_{\C^2})$ which is defined by  
\[
U:=\{ z\in \C^2\colon z_1\neq 0, \abs{z_2}^2>\abs{z_1}^2-(1+\varepsilon ) \} . 
\]
$S^1$ acts on $(U, \om_{\C^2})$ by 
\[
u\cdot z:=(z_1, u^{-1}z_2)
\]
for $u\in S^1$ and $z\in U$. This action is Hamiltonian and the map $\mu_{S^1}\colon U\to \R$ defined by 
\[
\mu_{S^1}(z):=1-\abs{z_2}^2 
\]
is a moment map of this action. We denote by $\overline{X}$ the cut space by the symplectic cutting with respect to the $S^1$-action, namely, $\overline{X}$ is the quotient space  
\[
\overline{X}:=\{ (z, w)\in U\times \C \colon 1-\abs{z_2}^2-\abs{w}^2=0\} /S^1
\]
by the $S^1$-action 
\[
u\cdot (z,w):=(z_1, u^{-1}z_2, u^{-1}w) . 
\]
By \cite{L}, $\overline{X}$ is a symplectic manifold and the restriction of the standard $T^2$-action on $\C^2$ to 
$U$ induces a Hamiltonian $T^2$-action on $\overline{X}$ with a moment map 
$\overline{\mu}([z, w]):=(\abs{z_1}^2, \abs{z_2}^2)$. We denote by $\overline{B}$ the image of $\overline{\mu}$. 
$\overline{B}$ is written by 
\[
\overline{B}=\{ \xi \in \R^2\colon \xi_1>0, 0\le \xi_2\le 1, \xi_2>\xi_1-(1+\varepsilon )\} . 
\]
Let $\overline{X}_1$ and $\overline{X}_2$ be the open subspaces of $\overline{X}$ which are defined by 
\[
\begin{split}
\overline{X}_1&:=\{ [z, w]\in \overline{X}\colon \abs{z_1}^2<\varepsilon \} ,\\ 
\overline{X}_2&:=\{ [z, w]\in \overline{X}\colon \abs{z_1}^2-(1+\varepsilon )<\abs{z_2}^2<\abs{z_1}^2-1\} . 
\end{split}
\]
We define the diffeomorphism $\varphi^X \colon \overline{X}_1\to \overline{X}_2$ by 
\[
\varphi^X([z, w]):=\left[ \left( (z_1/\abs{z_1})\sqrt{\abs{z_1}^2+\abs{z_2}^2+1}, 
\left( z_1/\abs{z_1}\right)^{-1}z_2\right) , w\right] . 
\]
It is easy to see that $\varphi^X$ is well defined and preserves symplectic structures of $\overline{X}_1$ and $\overline{X}_2$ induced by the symplectic structure on $\overline{X}$. The images of 
$\overline{X}_1$ and $\overline{X}_2$ by $\overline{\mu}$ are open subsets of $\overline{B}$ which are written by 
\[
\overline{\mu}(\overline{X}_1)= \{ \xi \in \overline{B}\colon \xi_1< \varepsilon  \} ,\ \  
\overline{\mu}(\overline{X}_2)=\{ \xi \in \overline{B}\colon \xi_2< \xi_1-1 \} . 
\]
$\varphi^X$ covers the diffeomorphism $\varphi^B\colon \overline{\mu}(\overline{X}_1)\to\overline{\mu}(\overline{X}_2)$ by 
\[
\varphi^B(\xi ):=(\xi_1+\xi_2+1, \xi_2). 
\]
This means that they satisfy 
\begin{equation}\label{com}
\overline{\mu}\circ \varphi^X=\varphi^B\circ \overline{\mu}. 
\end{equation}
We define $X$ to be the quotient space 
\[
X:=\overline{X}/\sim_X, 
\]
where $[z_1, w_1]\sim_X [z_2, w_2]$ if and only if $[z_i, w_i]\in \overline{X}_i$ for $i=1,2$ and 
$[z_2, w_2]=\varphi^X([z_1, w_1])$ and also define $B$ to be the quotient space 
\[
B:=\overline{B}/\sim_B, 
\]
where $\xi\sim_B \eta$ if and only if $\xi \in \overline{\mu}(\overline{X}_1)$, $\eta \in \overline{\mu}(\overline{X}_2)$, and $\eta=\varphi^B(\xi )$. Since $\overline{\mu}$ satisfies \eqref{com}, it descends to the map $\mu \colon X\to B$. By the construction, the symplectic structure on $\overline{X}$ induces a symplectic structure $\om$ on $X$ so that $\mu \colon (X,\om )\to B$ is a locally toric Lagrangian fibration. 
\end{ex}

We show that a locally toric Lagrangian fibration is endowed with a $C^{\infty}$ local torus action. First we recall the following fact about automorphisms of the local models of a nonsingular Lagrangian fibration. In this lemma we identify an automorphism of $T^n$ with an element of $\GL_n(\Z )$. 
\begin{lem}[{\cite[Lemma 2.5]{S}}]\label{local-auto}
Let $\varphi$ be a fiber-preserving symplectomorphism of $\pr_1\colon (\R^n\times T^n, \om_{\R^n\times T^n})\to \R^n$. Then, $\varphi$ is of the form $\varphi (\xi ,u)=(\rho^{-T}(\xi)+c, \rho(u)f(\xi ))$ for some $\rho \in \Aut (T^n)$, $c\in \R^n$ and a Lagrangian section $f$ of $\pr_1\colon (\R^n\times T^n, \om_{\R^n\times T^n})\to \R^n$, where $\rho^{-T}$ is the transpose inverse of $\rho$. The latter means that $f^*\om_{\R^n\times T^n}$ vanishes. 
\end{lem}

We show the following proposition. 
\begin{prop}\label{locally-toric-Lag}
Let $\mu \colon (X,\om)\to B$ be a locally toric Lagrangian fibration on an $n$-dimensional base $B$ and $\{ (U_{\alpha},\varphi^B_{\alpha},\varphi^X_{\alpha})\}$ a system of local identifications of $\mu$ with $\mu_{\C^n}$. Then, on each nonempty overlap $U_{\alpha \beta}:=U_{\alpha}\cap U_{\beta}$, there exists an element $\rho_{\alpha \beta}\in \Aut (T^n)$ such that $\varphi^X_{\alpha}\circ (\varphi^X_{\beta})^{-1}$ is $\rho_{\alpha \beta}$-equivariant. 
\end{prop}
\begin{proof}
First we focus on the interior $B\setminus \partial B$ of $B$. Since the restriction of $\mu \colon (X,\om)\to B$ to $B\setminus \partial B$ is a nonsingular Lagrangian fibration, it is locally identified with the local model. In fact, we can construct local identifications explicitly. For each $\alpha$, $\varphi^X_{\alpha}$ sends $(\mu^{-1}(U_{\alpha}\setminus \partial B), \om)$ fiber-preserving symplectomorphically to $(\mu_{\C^n}^{-1}(\varphi^B_{\alpha}(U_{\alpha}\setminus \partial B)),\om_{\C^n})$ which covers $\varphi^B_{\alpha}$. We define a fiber preserving symplectomorphism $\nu_{\C^n}\colon (\varphi^B_{\alpha}(U_{\alpha}\setminus \partial B)\times T^n,\om_{\R^n\times T^n})\to (\mu_{\C^n}^{-1}(\varphi^B_{\alpha}(U_{\alpha}\setminus \partial B)),\om_{\C^n})$ by 
\begin{equation}\label{nu_{C^n}}
\nu_{\C^n}(\xi ,u):=u\cdot \iota (\xi ), 
\end{equation}
where $\iota$ is the map defined by \eqref{iota}. Note that $T^n$ acts on $\varphi^B_{\alpha}(U_{\alpha}\setminus \partial B)\times T^n$ by the multiplication to the second factor and $\nu_{\C^n}$ is equivariant with respect to this action and the standard representation of $T^n$. Thus $\nu_{\C^n}^{-1}\circ \varphi^X_{\alpha}$ is a local identification. 

On each nonempty overlap $U_{\alpha \beta}$, by applying Lemma~\ref{local-auto} to the fiber preserving symplectomorphism $\nu_{\C^n}^{-1}\circ \varphi^X_{\alpha}\circ (\varphi^X_{\beta})^{-1}\circ \nu_{\C^n}\colon (\varphi^B_{\beta}(U_{\alpha \beta}\setminus \partial B),\om_{\R^n\times T^n})\to (\varphi^B_{\alpha}(U_{\alpha \beta}\setminus \partial B),\om_{\R^n\times T^n})$, on $U_{\alpha \beta}$ there exists an automorphism $\rho_{\alpha \beta}\in \Aut (T^n)$ such that  $\nu_{\C^n}^{-1}\circ \varphi^X_{\alpha}\circ (\varphi^X_{\beta})^{-1}\circ \nu_{\C^n}$ is $\rho_{\alpha \beta}$-equivariant. Since $\nu_{\C^n}$ is equivariant and $\mu^{-1}(U_{\alpha \beta}\setminus \partial B)$ is open dense in $\mu^{-1}(U_{\alpha \beta})$, the overlap map $\varphi^X_{\alpha}\circ (\varphi^X_{\beta})^{-1}\colon \mu_{\C^n}^{-1}(\varphi^B_{\beta}(U_{\alpha \beta}))\to \mu_{\C^n}^{-1}(\varphi^B_{\alpha}(U_{\alpha \beta}))$ is also $\rho_{\alpha \beta}$-equivariant. 
\end{proof}
This implies that $\{ (\mu^{-1}(U_{\alpha}),\varphi^X_{\alpha})\}$ is a weakly standard atlas of $X$. Hence, $X$ is equipped with a $C^{\infty}$ local $T^n$-action. It is obvious that the orbit space and the orbit map of the local $T^n$-action on $X$ are naturally identified with $B$ and $\mu$. In Section~\ref{symp}, we will see a necessary and sufficient condition in order that a manifold with a local torus action becomes a locally toric Lagrangian fibration.

\subsection{New local torus actions from given ones}
When a local torus action is given, we can construct new local torus actions from the given one. 
\begin{ex}[fiber products]
Let $(X,\CT )$ be a $2n$-dimensional manifold equipped with a $C^r$ local $T^n$-action. Suppose that $f\colon B'\to B_X$ is a stratification preserving locally homeomorphism from an $n$-dimensional $C^0$ manifold $B'$ with corners to $B_X$. Then it is obvious that the fiber product $f^*X:=\{ (b', x)\in B'\times X\colon f(b')=\mu_X(x)\}$ of $f$ and $\mu_X$ is equipped with a $C^r$ local $T^n$-action whose orbit space is $B'$. 
\end{ex}
\begin{ex}[blowing-ups]\label{blowing-up}
Let $(X,\CT )$ be a $2n$-dimensional manifold equipped with a $C^{\infty}$ local $T^n$-action. Let $x\in X$. Suppose that there exists a coordinate neighborhood $(U^X,\varphi^X)$ in a weakly standard atlas of $\CT$ such that $\varphi^X$ sends $x$ to the origin of $\C^n$. For a sufficiently small positive real number $\varepsilon >0$, we denote by $\varphi^X(U^X)_{\ge \varepsilon}$ the quotient space of $\{ (z, w)\in \varphi^X(U^X)\times \C \colon \norm{z}^2-\lvert w\rvert^2=\varepsilon \}$ by the circle action defined by 
\begin{equation*}
t\cdot (z,w):=(tz_1, \ldots ,tz_n, t^{-1}w). 
\end{equation*}
It is easy to see that $\varphi^X(U^X)_{\ge \varepsilon}$ is smooth and $T^n$ acts on $\varphi^X(U^X)_{\ge \varepsilon}$ by 
\[
u\cdot [z,w]:=[u\cdot z, w]. 
\]
Let $\overline{D}_{\varepsilon}(0)$ be the closed disc in $\C^n$ centered at the origin with radius $\varepsilon^{1/2}$. We define the diffeomorphism $f_{\varepsilon}$ from the open set 
$\varphi^X(U^X)\setminus \overline{D}_{\varepsilon}(0)$ to 
$\{ [z, w]\in \varphi^X(U^X)_{\ge \varepsilon}\colon w\neq 0\}$ by 
\begin{equation*}
f_{\varepsilon}(z)=\left[ z, (\norm{z}^2-\varepsilon)^{1/2}\right]
\end{equation*}
for $z\in \varphi^X(U^X)\setminus \overline{D}_{\varepsilon}(0)$. It is easy to see that $f_{\varepsilon}$ is an equivariant diffeomorphism. By removing $(\varphi^X_{\alpha})^{-1}(\overline{D}_{\varepsilon}(0))$ for a sufficiently small $\varepsilon >0$ from $X$ and gluing $X\setminus (\varphi^X_{\alpha})^{-1}(\overline{D}_{\varepsilon}(0))$ and $\varphi^X(U^X)_{\ge \varepsilon}$ by $f_{\varepsilon}\circ \varphi^X_{\alpha}$, we can obtain a new manifold equipped with a $C^{\infty}$ local $T^n$-action. We call it the {\it blowing-up} of $X$ at $x$. This is a $C^{\infty}$ aspect of the symplectic blowing-up by Guillemin and Sternberg in \cite{GuS1}. 
\end{ex}
\begin{ex}[connected sums]
Let $(X_1,\CT_1)$ and $(X_2,\CT_2)$ be $2n$-dimensional manifolds equipped with $C^r$ local $T^n$-actions. Suppose that $(U^{X_1},\varphi^{X_1})$ and $(U^{X_2}$, $\varphi^{X_2})$ are coordinate neighborhoods in weakly standard atlases of $\CT_1$ and $\CT_2$ such that both $\varphi^{X_1}$ and $\varphi^{X_2}$ send $U^{X_1}$ and $U^{X_2}$ $C^r$ diffeomorphically to a same $T^n$-invariant open set of $\C^n$. Then we can perform the connected sum $X_1\# X_2$ and the obtained manifold $X_1\# X_2$ has a $C^r$ local $T^n$-action whose orbit space is the connected sum $B_{X_1}\# B_{X_2}$. 
\end{ex}

\section{Characteristic pairs and canonical models}\label{cp-cm}
In this section, we introduce a characteristic pair and construct a canonical model from a characteristic pair. 
Both of them will play important roles in the rest of this paper. We also show that a characteristic pair is associated with a local torus action. 

\subsection{Characteristic pairs}\label{char-pair}
Let $B$ be an $n$-dimensional $C^0$ manifold with corners. We assume that $\partial B\neq \emptyset$. We denote by ${\mathcal S}^{(k)}B$ the $k$-dimensional part of the natural stratification of $B$. Let $\Lambda$ be the lattice of integral elements in the Lie algebra $\t$ of $T^n$, namely, 
\[
\Lambda :=\{ t\in \t \colon \exp t =1\} .
\]
Since the differential of any automorphism of $T^n$ at the unit element preserves $\Lambda$, by assigning to an automorphism of $T^n$ its differential at the unit element, there is a natural homomorphism from $\Aut (T^n)$ to $\GL (\Lambda)$. It is an isomorphism. In fact, it follows from the surjectivity of the exponential map of $T^n$ and the equation $\varphi \circ \exp =\exp \circ d\varphi$ for any automorphism $\varphi \in \Aut (T^n)$. In the rest of this paper, we identify $\Aut (T^n)$ with $\GL (\Lambda)$ by this isomorphism. Let $\pi_P\colon P\to B$ be a principal $\Aut (T^n)$-bundle on $B$ and $\pi_{\Lambda}\colon \Lambda_P\to B$ the associated $\Lambda$-bundle of $P$ by the above isomorphism $\Aut (T^n)\cong \GL (\Lambda)$. Suppose that $\pi_{\mathcal L}\colon {\mathcal L}\to {\mathcal S}^{(n-1)}B$ is a rank one subbundle of the restriction $\pi_{\Lambda} |_{\CS^{(n-1)}B}\colon \Lambda_P|_{\CS^{(n-1)}B}\to \CS^{(n-1)}B$ of $\pi_{\Lambda}\colon \Lambda_P \to B$ to ${\mathcal S}^{(n-1)}B$. For each $k$ and any point $b\in \CS^{(k)}B$, let $U$ be an open neighborhood of $b$ in $B$ on which there exists a local trivialization $\varphi^{\Lambda}\colon \pi_{\Lambda}^{-1}(U)\to U\times \Lambda$ of $\Lambda_P$. By shrinking $U$ if necessary, we can assume that the intersection $U\cap \CS^{(n-1)}B$ has exactly $n-k$ connected components, say, $(U\cap {\mathcal S}^{(n-1)}B)_1$, $\ldots$, $(U\cap {\mathcal S}^{(n-1)}B)_{n-k}$ (see Figure~4). 
\begin{figure}[hbtp]
\begin{center}
\input{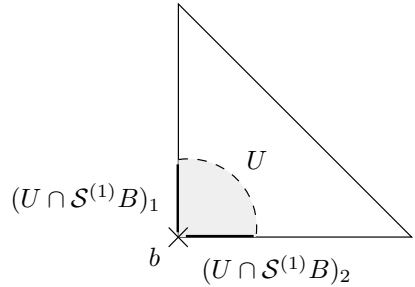}
\caption{An open neighborhood $U$ when $B$ is a triangle and $b\in \CS^{(0)}B$}
\label{fig23}
\end{center}
\end{figure}
Since $\Lambda$ is discrete, for each $(U\cap {\mathcal S}^{(n-1)}B)_a$ there exists a rank one sublattice $L_a\subset \Lambda$ such that $\varphi^{\Lambda}$ sends $\pi_{\CL}^{-1}\left( (U\cap {\mathcal S}^{(n-1)}B)_a\right)$ fiberwise to $(U\cap{\mathcal S}^{(n-1)}B)_a\times L_a$ 
\[
\begin{matrix}
\pi_{\Lambda}^{-1}(U) & \overset{\varphi^{\Lambda}}{\cong} & U\times \Lambda \\
\cup & & \cup \\
\pi_{\Lambda}^{-1}\left( (U\cap {\mathcal S}^{(n-1)}B)_a\right) & \cong & 
(U\cap {\mathcal S}^{(n-1)}B)_a\times \Lambda \\
\cup & & \cup \\
\pi_{\CL}^{-1}\left( (U\cap {\mathcal S}^{(n-1)}B)_a\right) & \cong & 
(U\cap {\mathcal S}^{(n-1)}B)_a\times L_a.  
\end{matrix}
\]
\begin{defn}\label{unimodular}
$\pi_{\CL}\colon \CL \to \CS^{(n-1)}B$ is said to be {\it unimodular} if for each $k$ and any point $b\in \CS^{(k)}B$, the sublattice $L_1+\cdots +L_{n-k}$ generated by the above rank one sublattices $L_1$, $\ldots$, $L_{n-k}$ is a rank $n-k$ direct summand of $\Lambda$. (In \cite{DJ} such a sublattice is called an {\it $(n-k)$-dimensional unimodular subspace} of $\Lambda$.) 
\end{defn}
Of course, $L_1$, $\ldots$, $L_{n-k}$ themselves depend on the choice of a neighborhood $U$ and a local trivialization $\varphi^{\Lambda}$, but Definition~\ref{unimodular} does not depend on the choice of them because unimodularity is invariant by an automorphism of $\Lambda$. 
\begin{defn}
Let $\pi_{\mathcal L}\colon {\mathcal L}\to {\mathcal S}^{(n-1)}B$ be a unimodular rank one subbundle of $\pi_{\Lambda} |_{\CS^{(n-1)}B}\colon \Lambda_P|_{\CS^{(n-1)}B}\to \CS^{(n-1)}B$. Then the pair $(P, \CL)$ of $\pi_P\colon P\to B$ and $\pi_{\mathcal L}\colon {\mathcal L}\to {\mathcal S}^{(n-1)}B$ is called a {\it characteristic pair} and $\pi_{\mathcal L}\colon {\mathcal L}\to {\mathcal S}^{(n-1)}B$ is called a {\it characteristic bundle}. 
\end{defn}

Let $(X,\CT )$ be a $2n$-dimensional manifold equipped with a local $T^n$-action. We show that there is a characteristic pair associated with $(X,\CT )$. Let $\{ (U_{\alpha}^X,$ $\varphi^X_{\alpha})\}_{\alpha \in \CA}\in \CT$ be the maximal weakly standard atlas. It induces an atlas $\{ (U_{\alpha}^B, \varphi^B_{\alpha})\}_{\alpha \in \CA}$ of $B_X$ which satisfies the property in Remark~\ref{std} and also determines a \v{C}ech one-cocycle $\{ \rho_{\alpha \beta}\}$ on $\{ U_{\alpha}^B\}_{\alpha \in \CA}$ with coefficients in $\Aut (T^n)$. It defines the principal $\Aut (T^n)$-bundle $\pi_{P_X}:P_X\to B_X$ on $B_X$ by setting 
\begin{equation}\label{P_X}
P_X:=\left( \coprod_{\alpha}U_{\alpha}^B\times \Aut (T^n)\right) /\sim_P,  
\end{equation}
where $(b_{\alpha}, h_{\alpha})\in U_{\alpha}^B\times \Aut (T^n)\sim_P(b_{\beta}, h_{\beta})\in U_{\beta}^B\times \Aut (T^n)$ if and only if $b_{\alpha}=b_{\beta}\in U_{\alpha \beta}^B$ and $h_{\alpha}=\rho_{\alpha \beta}\circ h_{\beta}$. The bundle projection $\pi_{P_X}$ is defined in the obvious way. For each $\alpha$, every point in $\pi_{P_X}^{-1}(U_{\alpha}^B)$ has a unique representative which lies in $U_{\alpha}^B\times \Aut (T^n)$. By associating a point in $\pi_{P_X}^{-1}(U_{\alpha}^B)$ with the unique representative, we define the local trivialization of $P_X$ on $U_{\alpha}^B$ which is denoted by $\varphi^P_{\alpha}:\pi_{P_X}^{-1}(U_{\alpha}^B)\to U_{\alpha}^B\times \Aut (T^n)$. Note that the following equality holds for $(b, h)\in U_{\alpha \beta}^B\times \Aut (T^n)$ 
\begin{equation}\label{overlap-P}
\varphi^P_{\alpha \beta}(b,h):=\varphi^P_{\alpha}\circ ( \varphi^P_{\beta})^{-1}(b,h)=(b, \rho_{\alpha \beta}\circ h) .
\end{equation}
We denote by $\pi_{\Lambda_X}:\Lambda_X\to B_X$ the $\Lambda$-bundle associated with $P_X$. We show that the property (ii) in Definition~\ref{localaction} determines a unique unimodular rank one subbundle of $\pi_{\Lambda_X}|_{\CS^{(n-1)}B_X}\colon \Lambda_X|_{\CS^{(n-1)}B_X}\to \CS^{(n-1)}B_X$. For each coordinate neighborhood $(U_{\alpha}^B, \varphi^B_{\alpha})$ of $B_X$ with $U_{\alpha}^B\cap \CS^{(n-1)}B_X\neq \emptyset$, $T^n$ acts on the preimage $\mu_{\C^n}^{-1}\left( \varphi^B_{\alpha}(U_{\alpha}^B\cap \CS^{(n-1)}B_X)\right)$ as the restriction of the standard representation of $T^n$. For simplicity, we assume that the intersection $U_{\alpha}^B\cap \CS^{(n-1)}B_X$ is connected. (Otherwise, we may consider componentwise.) Then, $\mu_{\C^n}^{-1}\left( \varphi^B_{\alpha}(U_{\alpha}^B\cap \CS^{(n-1)}B_X)\right)$ is fixed by a circle subgroup of $T^n$. We denote the circle subgroup by $S^1_{\alpha}$ and also denote the rank one sublattice of $\Lambda$ spanned by the integral element which generates $S^1_{\alpha}$ by $\CL_{\alpha}$. Suppose that $(U_{\alpha}^B, \varphi^B_{\alpha})$ and $(U_{\beta}^B, \varphi^B_{\beta})$ are coordinate neighborhoods satisfying the above conditions and the intersection $U_{\alpha \beta}^B\cap \CS^{(n-1)}B_X$ is nonempty. 
\begin{lem}\label{patch}
Under the identification $\Aut (T^n)\cong \GL (\Lambda)$, the automorphism $\rho_{\alpha \beta}$ of $T^n$ in the property $(\mathrm{ii})$ of Definition~\textup{\ref{localaction}} sends $\CL_{\beta}$ isomorphically to $\CL_{\alpha}$. 
\end{lem}
\begin{proof}
Since $\varphi^X_{\alpha \beta}$ is a $\rho_{\alpha \beta}$-equivariant homeomorphism, the automorphism 
$\rho_{\alpha \beta}$ of $T^n$ sends $S^1_{\beta}$ isomorphically to $S^1_{\alpha}$. See Proposition~\ref{observation2}. Therefore, $\rho_{\alpha \beta}$ also sends $\CL_{\beta}$ isomorphically to $\CL_{\alpha}$. 
\end{proof}

By the construction of $\pi_{\Lambda_X}\colon \Lambda_X\to B_X$, there exists a local trivialization $\varphi^{\Lambda}_{\alpha}\colon$ $\pi^{-1}_{\Lambda_X}(U_{\alpha}^B)\to U_{\alpha}^B\times \Lambda$ of $\pi_{\Lambda_X}\colon \Lambda_X\to B_X$ on each $U_{\alpha}^B$ such that on an overlap $U_{\alpha \beta}^B$ the transition function with respect to $\varphi^{\Lambda}_{\alpha}$ and $\varphi^{\Lambda}_{\beta}$ is $\rho_{\alpha \beta}$. We take a subsystem $\{ (U_{\alpha_i}^B, \varphi^B_{\alpha_i})\}_{i\in \CI}$ of $\{ (U_{\alpha}^B, \varphi^B_{\alpha})\}_{\alpha \in \CA}$ which covers $\CS^{(n-1)}B_X$. By Lemma~\ref{patch}, we can obtain the rank one subbundle $\pi_{\CL_X}\colon \CL_X \to \CS^{(n-1)}B_X$ of $\pi_{\Lambda_X}|_{\CS^{(n-1)}B_X}\colon \Lambda_X|_{\CS^{(n-1)}B_X}$ $\to \CS^{(n-1)}B_X$ by setting 
\begin{equation}\label{CL_X}
\CL_X:=\left( \coprod_iU_{\alpha_i}^B\cap \CS^{(n-1)}B_X\times \CL_{\alpha_i}\right) /\sim_L, 
\end{equation}
where $(b_i, l_i)\in U_{\alpha_i}^B\cap \CS^{(n-1)}B_X\times \CL_{\alpha_i}\sim_L (b_j, l_j)\in U_{\alpha_j}^B\cap \CS^{(n-1)}B_X\times \CL_{\alpha_j}$ if and only if $b_i=b_j$ and $l_i=\rho_{\alpha_i\alpha_j}(l_j)$. By the construction, it is easy to see that $\pi_{\CL_X}\colon \CL_X \to \CS^{(n-1)}B_X$ is unimodular. As a summary, we have the following proposition. 
\begin{prop}\label{char}
Associated with a local $T^n$-action $\CT$ on $X$, there exists a characteristic pair $(P_X, \CL_X)$, where $P_X$ and $\CL_X$ are defined by \textup{\eqref{P_X}} and \textup{\eqref{CL_X}}, respectively. 
\end{prop}

\begin{ex}
For a $2n$-dimensional manifold $X$ equipped with a locally standard $T^n$-action, $\pi_{P_X}\colon P_X\to B_X$ is the product principal $\Aut (T^n)$-bundle $P_X=B_X\times \Aut (T^n)$. Let $(\CS^{(n-1)}B_X)_a$ ($a=1, \ldots , k$) be the connected component of $\CS^{(n-1)}B_X$. Then each preimage $\mu_X^{-1}((\CS^{(n-1)}B_X)_a)$ is fixed by a circle subgroup $S^1_a$ of $T^n$. Let $L_a$ be the rank one sublattice in $\Lambda$ corresponding to $S^1_a$. Then, $\CL_X$ is the disjoint union $\coprod_a(\CS^{(n-1)}B_X)_a\times L_a$. In particular, if $X$ is a quasi-toric manifold, $\CL_X$ is essentially nothing but the characteristic function\footnote{Recall that the characteristic function is the map that associates $L_a$ with each $(\CS^{(n-1)}B_X)_a$.}. Since $\CL_X$ is unimodular for a locally standard $T^n$-action on $X$, it is easy to see that $B_X$ is {\it nice}, namely, each connected component of a $k$-dimensional part $\CS^{(k)}B$ of $B$ is contained in exactly $n-k$ facets.  
\end{ex}
\begin{ex}\label{CL-ex3.9}
In the case of Example~\ref{ex3.9}, we identify $\Lambda$ with $\Z \oplus \Z$ and also identify $\Aut (T^2)$ with $\GL_2(\Z )$. Then, $\pi_{P_X}\colon P_X\to B$ is the principal $\Aut (T^2)$-bundle associated with the universal covering of $B$ by $\rho$. The restriction of $\pi_{\Lambda_X}\colon \Lambda_X\to B$ to $\CS^{(1)}B$ (in this case, $n=2$) and the characteristic bundle $\pi_{\CL_X}\colon \CL_X\to \CS^{(1)}B$ are obtained by gluing two pairs 
\[
\left( \{ \xi \in \overline{B}_2\colon 0<\xi_1,\ \xi_2=0\} \times \Lambda ,\{ \xi \in \overline{B}_2\colon 0<\xi_1,\ \xi_2=0\} \times (\{ 0\}\oplus\Z ) \right) 
\]
and 
\[
\left( \{ \xi \in \overline{B}_2\colon \xi_1=0,\ 0<\xi_2\} \times \Lambda ,\{ \xi \in \overline{B}_2\colon \xi_1=0,\ 0<\xi_2\} \times ( \Z \oplus\{ 0\}) \right) 
\]
of trivial $\Lambda$-bundles and their rank one subbundles with the map $\varphi_{\Lambda}\colon \{ \xi \in U_1\colon 0<\xi_1,\ \xi_2=0\}\times \Lambda \to \{ \xi \in U_2\colon \xi_1=0,\ 0<\xi_2\} \times \Lambda$ defined by $\varphi_{\Lambda}(\xi ,\lambda):=(\varphi_B(\xi ), \rho (\gamma )(\lambda ))$. 
\end{ex}

For $i=1,2$, let $B_i$ be an $n$-dimensional $C^0$ manifold with corners and $(P_i, \CL_i)$ a characteristic pair on $B_i$. 
\begin{defn}\label{iso-char-pair}
An {\it isomorphism} $f_P\colon (P_1, \CL_1)\to (P_2, \CL_2)$ between characteristic pairs is a bundle isomorphism $f_P\colon P_1\to P_2$ which covers a stratification preserving homeomorphism $f_B\colon B_1\to B_2$ such that the lattice bundle isomorphism $f_{\Lambda}:\Lambda_{P_1}\to \Lambda_{P_2}$ induced by $f_P$ sends $\CL_1$ isomorphically to $\CL_2$. $(P_1, \CL_1)$ and $(P_2, \CL_2)$ are {\it isomorphic} if there exists an isomorphism between them. 
\end{defn}
The isomorphism class of the characteristic pair $(P_X, \CL_X)$ is an invariant of a local $T^n$-action on $X$. 
\begin{lem}\label{inv1}
For $i=1,2$, let $(X_i,\CT_i)$ be a $2n$-dimensional manifold with a local $T^n$-action. If there is a $C^0$ isomorphism $f_X\colon (X_1,\CT_1)\to (X_2,\CT_2)$, then $f_X$ induces an isomorphism $f_{P_X}\colon (P_{X_1}, \CL_{X_1})\to (P_{X_2}, \CL_{X_2})$ between characteristic pairs associated with $X_1$ and $X_2$. 
\end{lem}
\begin{proof}
Let $\{ (U_{\beta}^{X_1}, \varphi^{X_1}_{\beta})\}_{\beta \in \CB}\in \CT_1$ and $\{ (U_{\alpha}^{X_2}, \varphi^{X_2}_{\alpha})\}_{\alpha \in \CA}\in \CT_2$ be maximal weakly standard atlases of $X_1$ and $X_2$, and $\{ (U_{\beta}^{B_1}, \varphi^{B_1}_{\beta})\}_{\beta \in \CB}$ and $\{ (U_{\alpha}^{B_2}, \varphi^{B_2}_{\alpha})\}_{\alpha \in \CA}$ atlases of $B_{X_1}$ and $B_{X_2}$ induced by $\{ (U_{\beta}^{X_1}, \varphi^{X_1}_{\beta})\}_{\beta \in \CB}$ and $\{ (U_{\alpha}^{X_2}, \varphi^{X_2}_{\alpha})\}_{\alpha \in \CA}$, respectively. Suppose that $f_X\colon (X_1,\CT_1)\to (X_2,\CT_2)$ is a $C^0$ isomorphism and $f_B\colon B_{X_1}\to B_{X_2}$ is the homeomorphism induced by $f_X$. By definition, on each nonempty overlap $U^{B_1}_{\beta}\cap (f_B)^{-1}(U_{\alpha}^{B_2})$, there exists an automorphism $\rho_{\alpha \beta}^f$ of $T^n$ such that $\varphi^{X_2}_{\alpha}\circ f_X\circ (\varphi^{X_1}_{\beta})^{-1}$ is $\rho_{\alpha \beta}^f$-equivariant. It is easy to see that the following equality holds 
\begin{equation}\label{rho-iso}
\rho^f_{\alpha_0, \beta_0}\circ \rho^{X_1}_{\beta_0\beta_1}
=\rho^{X_2}_{\alpha_0\alpha_1}\circ \rho^f_{\alpha_1\beta_1}
\end{equation}
on a nonempty intersection $U^{B_1}_{\beta_0 \beta_1}\cap (f_B)^{-1}(U^{B_2}_{\alpha_0\alpha_1})$, where $\rho^{X_1}_{\beta_0\beta_1}$ and $\rho^{X_2}_{\alpha_0\alpha_1}$ are automorphisms of $T^n$ in (ii) of Definition~\ref{localaction} with respect to $X_1$ and $X_2$, respectively. We define the bundle isomorphism $(f_P)_{\alpha \beta}\colon U^{B_1}_{\beta}\cap f_B^{-1}(U^{B_2}_{\alpha})\times \Aut (T^n)\to f_B(U^{B_1}_{\beta})\cap U^{B_2}_{\alpha}\times \Aut (T^n)$ by 
\[
(f_P)_{\alpha \beta}(b,h):=(f_B(b),\rho_{\alpha \beta}^f\circ h) .
\]
Then, \eqref{rho-iso} implies that the equation 
\[
(f_P)_{\alpha_0 \beta_0}\circ \varphi^{P_{X_1}}_{\beta_0 \beta_1}=\varphi^{P_{X_2}}_{\alpha_0 \alpha_1}\circ (f_P)_{\alpha_1 \beta_1}
\]
holds on $U^{B_1}_{\beta_0 \beta_1}\cap f_B^{-1}(U^{B_2}_{\alpha_0\alpha_1})$, where $\varphi^{P_{X_1}}_{\beta_0 \beta_1}$ (resp. $\varphi^{P_{X_2}}_{\alpha_0 \alpha_1}$) is the overlap map defined by \eqref{overlap-P} for $\pi_{P_{X_1}}\colon P_{X_1}\to B_{X_1}$ (resp. $\pi_{P_{X_2}}\colon P_{X_2}\to B_{X_2}$) on $U^{B_1}_{\beta_0 \beta_1}$ (resp. $U^{B_2}_{\alpha_0\alpha_1}$). Therefore, we can patch them together to obtain the bundle isomorphism $f_P\colon P_{X_1}\to P_{X_2}$ which covers $f_B$. This proves the lemma. 
\end{proof}

\subsection{Canonical models}\label{sub-canonical}
In \cite[Section 1.5]{DJ}, Davis and Januszkiewicz constructed the canonical model of a quasi-toric manifold from the based polytope and the characteristic function. A similar construction can be done by using a characteristic pair in the following way. Let $B$ be an $n$-dimensional $C^0$ manifold with corners and $(P, \CL )$ a characteristic pair on $B$. We denote by $\pi_T:T_P\to B$ the $T^n$-bundle associated with $P$ by the natural action of $\Aut (T^n)$ on $T^n$. First we shall explain that for any $k$-dimensional part $\CS^{(k)}B$, $(P, \CL)$ determines a rank $n-k$ subtorus bundle of the restriction of $\pi_T:T_P\to B$ to $\CS^{(k)}B$. Let $\{ U_{\alpha}\}$ be an open covering of $B$ such that on each $U_{\alpha}$ there exists a local trivialization $\varphi^P_{\alpha}\colon \pi_P^{-1}(U_{\alpha})\to U_{\alpha}\times \Aut (T^n)$. On each nonempty overlap $U_{\alpha \beta}$ we denote by $\rho_{\alpha \beta}$ the transition function with respect to $\varphi^P_{\alpha}$ and $\varphi^P_{\beta}$, namely,  
\[
\varphi^P_{\alpha}\circ (\varphi^P_{\beta})^{-1}(b,f)=(b, \rho_{\alpha \beta}f)
\]
for $(b,f)\in U_{\beta}\times \Aut (T^n)$. Note that $\rho_{\alpha \beta}$ is locally constant since $\Aut (T^n)$ is discrete. For simplicity we assume that each $U_{\alpha \beta}$ is connected so that $\rho_{\alpha \beta}$ can be thought of as an element of $\Aut (T^n)$. $\varphi^P_{\alpha}$ induces a local trivializations of the associated bundles $T_P$ and $\Lambda_P$ which are denoted by $\varphi^T_{\alpha}:\pi_T^{-1}(U_{\alpha})\to U_{\alpha}\times T^n$ and $\varphi^{\Lambda}_{\alpha}:\pi_{\Lambda}^{-1}(U_{\alpha})\to U_{\alpha}\times \Lambda$, respectively. For $\CS^{(k)}B$ we take $U_{\alpha}$ with $U_{\alpha}\cap \CS^{(k)}B\neq \emptyset$. By replacing $U_{\alpha}$ by a sufficiently small one if necessary, we may assume that the intersection $U_{\alpha}\cap \CS^{(n-1)}B$ has exactly $n-k$ connected components, say $(U_{\alpha}\cap \CS^{(n-1)}B)_1$, $\cdots$, $(U_{\alpha}\cap \CS^{(n-1)}B)_{n-k}$. For $k=n$, this means that $U_{\alpha}$ is contained in $\CS^{(n)}B$. For $k<n$, there are $n-k$ rank one sublattices $L_1$, $\ldots$, $L_{n-k}$ of $\Lambda$ such that for $a=1$, $\ldots\ $, $n-k$ $\varphi^{\Lambda}_{\alpha}$ sends the restriction of $\pi_{\CL}\colon \CL\to \CS^{(n-1)}B$ to $(U_{\alpha}\cap \CS^{(n-1)}B)_a$ isomorphically to the trivial rank one subbundle $(U_{\alpha}\cap \CS^{(n-1)}B)_a\times L_a$ of $(U_{\alpha}\cap \CS^{(n-1)}B)_a\times \Lambda$. Since $\CL$ is unimodular, $L_1$, $\ldots$, $L_{n-k}$ generate the $(n-k)$-dimensional subtorus of $T^n$ which is denoted by $Z_{U_{\alpha}\cap \CS^{(k)}B}$. For $k=n$, we define $Z_{U_{\alpha}\cap \CS^{(n)}B}$ to be the trivial subgroup. Note that when $(P,\CL )$, $\{ U_{\alpha}\}$, and $\varphi^P_{\alpha}$ are induced by some local $T^n$-action $\CT$ on $X$, $Z_{U_{\alpha}\cap \CS^{(k)}B_X}$ is the common ($n-k$)-dimensional stabilizer of $T^n$-action on $\mu_{\C^n}^{-1}(U_{\alpha}\cap \CS^{(k)}B_X)$. 

Suppose that another $U_{\beta}$ satisfies the above condition and $U_{\alpha \beta}\cap \CS^{(k)}B\neq \emptyset$. By the definition of $(P,\CL )$, $\rho_{\alpha \beta}$ sends $Z_{U_{\beta}\cap \CS^{(k)}B_X}$ isomorphically to $Z_{U_{\alpha}\cap \CS^{(k)}B_X}$. Hence, by patching them together by using $\rho_{\alpha \beta}$'s, we obtain a rank $n-k$ subtorus bundle, which is denoted by $\pi_{Z_{\CS^{(k)}B}}\colon Z_{\CS^{(k)}B}\to \CS^{(k)}B$, of the restriction of $\pi_T\colon T_P\to B$ to $\CS^{(k)}B$. 
\begin{defn}\label{equiv}
We define two elements $t$, $t'\in T_P$ to be {\it equivalent} or $t\sim_{can} t'$ if $\pi_T(t)=\pi_T(t')$ and $t't^{-1}\in \pi_{Z_{\CS^{(k)}B}}^{-1}(\pi_T(t))$ provided that $\pi_T(t)$ lies in $\CS^{(k)}B$. Note that a fiber of $\pi_T\colon T_P\to B$ is equipped with the structure of a group since its structure group is $\Aut (T^n)$. 
\end{defn}
We denote by $X_{(P, \CL )}$ the quotient space of $T_P$ by the equivalence relation. The bundle projection $\pi_T\colon T_P\to B$ descends to the map $\mu_{X_{(P, \CL )}}\colon X_{(P,\CL )}\to B$. Since $\CL$ is unimodular and $B$ is a manifold with corners, by the same way as in Davis and Januszkiewicz~\cite[Section 1.5]{DJ}, or Masuda and Panov~\cite[Section 3.2]{MP}, we can show that $X_{(P, \CL )}$ is equipped with a $C^0$ local $T^n$-action whose orbit space is $B$ and whose orbit map is $\mu_{X_{(P, \CL)}}$. 
\begin{defn}
We call $X_{(P,\CL )}$ the {\it canonical model} associated with $(P,\CL )$. In particular, when $(P,\CL )$ is the characteristic pair $(P_X,\CL_X)$ of a local $T^n$-action $\CT$ on a $2n$-dimensional manifold $X$, we also call $X_{(P_X,\CL_X)}$ the canonical model associated with $(X,\CT )$. 
\end{defn}

Note that by the construction, a fiber of $\mu_{X_{(P, \CL )}}\colon X_{(P,\CL )}\to B$ admits a group structure. 

We give properties of a canonical model. 
\begin{prop}\label{section-can}
For a characteristic pair $(P,\CL )$, $\mu_{X_{(P, \CL )}}\colon X_{(P,\CL )}\to B$ admits a $C^0$ section. 
\end{prop}
\begin{proof}
Since a fiber of $\pi_T\colon T_P\to B$ admits a structure of a group, it admits the section which assigns to an element $b\in B$ the unit element of $\pi_T^{-1}(b)$. The composition of this section and the natural projection from $T_P$ to $X_{(P,\CL )}$ defines the $C^0$ section $s$ of $\mu_{X_{(P, \CL )}}\colon X_{(P,\CL )}\to B$. 
\end{proof}

\begin{prop}\label{group-action}
For a $2n$-dimensional manifold $(X,\CT )$ equipped with a local $T^n$-action, we denote the associated $T^n$-bundle $T_{P_X}$ of $P_X$ by $\pi_{T_X}\colon T_X\to B_X$ for simplicity. Then $T_X$ acts fiberwise on $X$. Similarly $X_{(P_X,\CL_X )}$ also acts fiberwise on $X$. For any $b\in B_X$ the action of $\mu_{X_{(P, \CL )}}^{-1}(b)$ on $\mu_X^{-1}(b)$ is simply transitive. 
\end{prop}
\begin{proof}
Let $\{ (U_{\alpha}^X, \varphi^X_{\alpha})\}_{\alpha \in \CA}\in \CT$ be a weakly standard atlas of $X$ and $\{ (U_{\alpha}^B, \varphi^B_{\alpha})\}_{\alpha \in \CA}$ the atlas of $B_X$ induced by $\{ (U_{\alpha}^X, \varphi^X_{\alpha})\}_{\alpha \in \CA}$. For $b\in B_X$ which lies in $U_{\alpha}^B$, the action of $\pi_{T_X}^{-1}(b)$ on $\mu_X^{-1}(b)$ can be written by
\[
t\cdot x:=(\varphi^X_{\alpha})^{-1}(t_{\alpha}\cdot \varphi^X_{\alpha}(x))
\]
for $t\in \pi_{T_X}^{-1}(b)$ and $x\in \mu_X^{-1}(b)$, where $\varphi^T_{\alpha}(t)=(b,t_{\alpha})$. Recall that $\varphi^T_{\alpha}\colon \pi_{T_X}^{-1}(U_{\alpha}^B)\to U_{\alpha}^B\times T^n$ denotes the local trivialization of $T_X$ which is induced by the local trivialization $\varphi^P_{\alpha}$ of $P_X$. It does not depend on the choice of $U_{\alpha}^B$. In fact, suppose that $b$ also lies in $U_{\beta}^B$ and $\varphi^T_{\beta}(t)=(b,t_{\beta})$. Since the transition function of $\varphi^T_{\alpha}$ and $\varphi^T_{\beta}$ is $\rho_{\alpha \beta}$, we have $t_{\alpha}=\rho_{\alpha \beta}(t_{\beta})$. By using the $\rho_{\alpha \beta}$-equivariantness of $\varphi^X_{\alpha \beta}$, 
\[
\begin{split}
(\varphi^X_{\alpha})^{-1}(t_{\alpha}\cdot \varphi^X_{\alpha}(x))&=(\varphi^X_{\beta})^{-1}\circ \varphi^X_{\beta \alpha}(\rho_{\alpha\beta}(t_{\beta})\cdot \varphi^X_{\alpha \beta}\circ \varphi^X_{\beta}(x))\\
&=(\varphi^X_{\beta})^{-1}(t_{\beta}\cdot \varphi^X_{\beta}(x)). 
\end{split}
\]
By Proposition~\ref{observation} and the construction of $X_{(P_X,\CL_X)}$, for each $b\in B_X$ the action of $\pi_{T_X}^{-1}(b)$ on $\mu_X^{-1}(b)$ descends to the action of $\mu_{X_{(P_X, \CL_X )}}^{-1}(b)$ on $\mu_X^{-1}(b)$ which is simply transitive. 
\end{proof}

\begin{lem}\label{inv2}
For $i=1$, $2$, let $B_i$ be an $n$-dimensional $C^0$ manifold with corners and $(P_i, \CL_i)$ a characteristic pair on $B_i$. Then, any isomorphism $f_P\colon (P_1, \CL_1)\to (P_2, \CL_2)$ induces the $C^0$ isomorphism $f_{X_{(P,\CL )}}\colon X_{(P_1, \CL_1)}\to X_{(P_2, \CL_2)}$ between canonical models of $(P_1, \CL_1)$ and $(P_2, \CL_2)$. 
\end{lem}
\begin{proof}
It follows directly from the construction of a canonical model. 
\end{proof}
\begin{rem}\label{group-str}
If there is an isomorphism $f_P\colon (P_1, \CL_1)\to (P_2, \CL_2)$ between characteristic pairs, then the induced $C^0$ isomorphism $f_{X_{(P,\CL )}}\colon X_{(P_1, \CL_1)}\to X_{(P_2, \CL_2)}$ is a fiberwise group isomorphism. 
\end{rem}

\section{On sections of orbit maps}\label{orbit-map-section}
Let $(X,\CT)$ be a $2n$-dimensional manifold equipped with a local $T^n$-action. 
\begin{prop}\label{section}
If $\mu_X \colon X\to B_X$ has a $C^0$ section $s\colon B_X\to X$, then there exists a weakly standard $C^0$ atlas $\{ (U_{\alpha}^X, \varphi^X_{\alpha})\}_{\alpha \in \CA}\in \CT$ such that for each $\alpha$ the following diagram commutes: 
\[
\xymatrix{
\mu_X^{-1}(U_{\alpha}^B)\ar[r]^{\varphi^X_{\alpha}} & \mu_{\C^n}^{-1}(\varphi^B_{\alpha}(U_{\alpha}^B)) \\
U_{\alpha}^B\ar[u]^s\ar[r]^{\varphi^B_{\alpha}} & \varphi^B_{\alpha}(U_{\alpha}^B)\ar[u]^{\iota}, 
}
\]
where $\{ (U_{\alpha}^B,\varphi_{\alpha}^B)\}_{\alpha \in \CA}\in \CT$ is the atlas of $B_X$ induced by $\{ (U_{\alpha}^X, \varphi^X_{\alpha})\}_{\alpha \in \CA}\in \CT$ which satisfies the property in Remark~\textup{\ref{std}} and $\iota$ is the section defined by \eqref{iota}. 
\end{prop}
\begin{proof}
Let $\{ (U_{\alpha}^X,\psi^X_{\alpha})\} \in \CT$ be a weakly standard $C^0$ atlas of $X$ and $\{ (U_{\alpha}^B,\psi^B_{\alpha})\} \in \CT$ the atlas of $B_X$ induced by $\{ (U_{\alpha}^X, \psi^X_{\alpha})\}$ which satisfies the property in Remark~\textup{\ref{std}}. For each $\alpha$ the equality 
\[
\theta_{\alpha}(b)\cdot s(b)=(\psi^X_{\alpha})^{-1}\circ \iota \circ \psi^B_{\alpha}(b)
\]
for $b\in U_{\alpha}^B$ determines a unique local section $\theta_{\alpha}$ of $\mu_{X_{(P_X, \CL_X)}}\colon X_{(P_X,\CL_X)}\to B_X$ on $U_{\alpha}^B$. Now we define a new local coordinate system $\varphi^X_{\alpha}$ on $\mu_X^{-1}(U_{\alpha}^B)$ by 
\[
\varphi^X_{\alpha}(x):=\psi^X_{\alpha}(\theta_{\alpha}(\mu_X(x))\cdot x)
\]
for $x\in \mu_X^{-1}(U_{\alpha}^B)$. Then $\{ (U_{\alpha}^X,\varphi^X_{\alpha})\}$ is the required weakly standard atlas. 
\end{proof}
By Proposition~\ref{T^n-bundle}, the restriction $\mu_X |_{\CS^{(n)}B_X} \colon \mu_X^{-1}(\CS^{(n)}B_X) \to \CS^{(n)}B_X$ of $\mu_X$ to $\CS^{(n)}B_X$ of $B_X$ is a $T^n$-bundle whose structure group is contained in $T^n\rtimes \Aut (T^n)$. This proposition implies that the structure group can be reduced to a subgroup of $\Aut (T^n)$ if $\mu_X$ admits a section. 

The following proposition shows that a section of the orbit map is unique up to $C^0$ isomorphisms. 
\begin{prop}\label{section-preserving}
Suppose that $\mu_X \colon X\to B_X$ has two sections $s_1$ and $s_2$. Then there exists a $C^0$ isomorphism $f_X$ of $X$ such that $f_X$ covers the identity on $B_X$ and preserves sections, namely, $f_X \circ s_1=s_2$.  
\end{prop}
\begin{proof}
Since the fiberwise action of $X_{(P_X,\CL_X)}$ on $X$ is simply transitive, there exists a unique $C^0$ section $\theta$ of $\mu_{X_{(P_X, \CL_X )}}$ such that 
\[
\theta (b)\cdot s_1(b)=s_2(b). 
\]
Then the required $C^0$ isomorphism $f_X$ of $X$ can be obtained by 
\[
f_X(x):=\theta (\mu_X(x))\cdot x.
\]
\end{proof}

The following lemma is a straightforward generalization of the result \cite[Proposition 1.8]{DJ} by Davis and Januszkiewicz and the result \cite[Lemma 3.6]{MP} by Masuda and Panov. 
\begin{lem}\label{canonical-model}
$\mu_X\colon X\to B_X$ is equipped with a $C^0$ section if and only if there exists a $C^0$ isomorphism between $(X,\CT )$ and $X_{(P_X,\CL_X)}$ which covers identity on $B_X$. 
\end{lem}
\begin{proof}
The if part follows from Proposition~\ref{section-can}. Conversely suppose that $\mu_X \colon X\to B_X$ has a section which is denoted by $s$. Define the surjective map $\nu :T_X\to X$ by 
\begin{equation}\label{nu}
\nu (t):=t\cdot s(\pi_{T_X}(t)) . 
\end{equation}
It is easy to see from the construction of $X_{(P_X, \CL_X)}$ that $\nu\colon T_X\to X$ descends to the required $C^0$ isomorphism from $X_{(P_X, \CL_X)}$ to $X$. 
\end{proof}

Next, we investigate when $\mu_X \colon X\to B_X$ has a section. We assume that the index set $\CA$ of the weakly standard atlas $\{ (U_{\alpha}^X, \varphi^X_{\alpha})\}_{\alpha \in \CA}$ is countable ordered. By the construction of $X_{(P_X,\CL_X)}$, there exists a $C^0$ isomorphism $h_{\alpha}\colon \mu_X^{-1}(U_{\alpha}^B)\to \mu_{X_{(P_X,\CL_X)}}^{-1}(U_{\alpha}^B)$ covering the identity on each $U_{\alpha}^B$ such that $h_{\alpha}$ is equivariant with respect to the fiberwise action of $T_X$ and $X_{(P_X,\CL_X)}$. By Proposition~\ref{group-action}, on each nonempty overlap $U_{\alpha \beta}^B$ the equation 
\begin{equation}\label{theta}
h_{\alpha}\circ h_{\beta}^{-1}(x)=\theta_{\alpha \beta}^X(b) x
\end{equation}
for $b\in U_{\alpha \beta}^B$ and $x\in \mu_{X_{(P_X,\CL_X)}}^{-1}(b)$ determines a unique local section $\theta_{\alpha \beta}^X$ of $\mu_{X_{(P_X, \CL_X)}}\colon$ $X_{(P_X,\CL_X)}\to B_X$ on $U_{\alpha \beta}^B$. Let $\SS_{(P_X, \CL_X)}$ denote the sheaf of $C^0$ sections of $\mu_{X_{(P_X, \CL_X)}}\colon$ $X_{(P_X,\CL_X)}\to B_X$. Then local sections $\theta_{\alpha \beta}^X$ form a \v{C}ech one-chain $\{ \theta_{\alpha \beta}^X\}$ on $\{ U^B_{\alpha}\}$ with values in $\SS_{(P_X, \CL_X)}$. 
\begin{lem}
$\{ \theta_{\alpha \beta}^X\}$ is a cocycle. 
\end{lem}
\begin{proof}
It is clear by \eqref{theta} and Proposition~\ref{group-action}. 
\end{proof}

Let $H^1(B_X;\SS_{(P_X, \CL_X)})$ denote the first \v{C}ech cohomology group of $B_X$ with values in $\SS_{(P_X, \CL_X)}$. By the above lemma, $\{ \theta_{\alpha \beta}^X\}$ defines a cohomology class in $H^1(B_X;\SS_{(P_X, \CL_X)})$. We denote it by $e(X,\CT )$. It is easy to see that $e(X,\CT )$ does not depend on the choice of $h_{\alpha}$'s and depends only on the local $T^n$-action on $X$. 
\begin{defn}\label{Euler}
We call $e(X,\CT )$ the {\it Euler class of $\mu_X\colon X\to B_X$}. 
\end{defn}

\begin{thm}\label{obst-reduction}
$\mu_X \colon X\to B_X$ has a section if and only if $e(X,\CT )$ vanishes. 
\end{thm}
\begin{proof}
By the definition of $e(X,\CT )$, $(X,\CT)$ is $C^0$ isomorphic to $X_{(P_X, \CL_X)}$ which covers the identity on $B_X$ if and only if $e(X,\CT )$ vanishes. Then the theorem follows directly from this fact together with Lemma~\ref{canonical-model}. 
\end{proof}

\begin{ex}
For an effective $C^{\infty}$ $T^2$-action on a $4$-dimensional $C^{\infty}$ manifold $X$ without nontrivial finite stabilizers, Orlik and Raymond showed in \cite{OR} that the orbit map is equipped with a section. Thus, $e(X,\CT )$ vanishes. 
\end{ex}
\begin{ex}
In the case of a nonsingular toric variety, $e(X,\CT )$ vanishes. In fact, the system of standard coordinate neighborhoods induced by top-dimensional cones has the property in Proposition~\ref{section}. Similarly, $e(X,\CT )$ vanishes for a quasi-toric manifold. It follows from \cite[Lemma 1.4]{DJ}. 
\end{ex}
\begin{ex}
If the local $T^n$-action is induced by a locally standard $T^n$-action and $\partial B_X=\emptyset$, then $\mu_X\colon X\to B_X$ is a principal $T^n$-bundle. In this case, $e(X,\CT )$ is nothing but the Euler class of the principal $T^n$-bundle. 
\end{ex}

For $i=1,2$, let $B_i$ be an $n$-dimensional $C^0$ manifold with corners and $(P_i, \CL_i)$ a characteristic pair on $B_i$. Suppose that there exists an isomorphism $f_P\colon (P_1, \CL_1)\to (P_2, \CL_2)$. By Lemma~\ref{inv2}, it induces the isomorphism $f_P^*\colon H^1(B_2;\SS_{(P_2, \CL_2)})\to H^1(B_1;\SS_{(P_1, \CL_1)})$ between cohomology groups. In particular, by Lemma~\textup{\ref{inv1}} and Lemma~\textup{\ref{inv2}}, a $C^0$ isomorphism $f_X\colon (X_1,\CT_1)\to (X_2,\CT_2)$ induces an isomorphism $f_{P_X}^*\colon H^1(B_{X_2};\SS_{(P_{X_2}, \CL_{X_2})})\to H^1(B_{X_1};\SS_{(P_{X_1}, \CL_{X_1})})$. 
\begin{lem}\label{inv3}
For $i=1,2$, let $(X_i,\CT_i)$ be a $2n$-dimensional manifold equipped with a local $T^n$-action. If there is a $C^0$ isomorphism $f_X\colon X_1\to X_2$, then $f_{P_X}^*e(X_2,\CT_2 )=e(X_1,\CT_1 )$. 
\end{lem}

\section{The topological classification}\label{classification thm}
The following is the main theorem of this paper. 
\begin{thm}\label{classification}
For $i=1,2$, let $(X_i,\CT_i)$ be a $2n$-dimensional manifold with a local $T^n$-action. $(X_1,\CT_1)$ and $(X_2,\CT_2)$ are $C^0$ isomorphic if and only if there exists an isomorphism $f_P\colon (P_{X_1}, \CL_{X_1})\to (P_{X_2}, \CL_{X_2})$ between characteristic pairs associated with $(X_1,\CT_1)$ and $(X_2,\CT_2)$ such that $f_P^*e(X_2,\CT_2 )=e(X_1,\CT_1)$. Moreover, for any characteristic pair $(P, \CL)$ on an $n$-dimensional $C^0$ manifold $B$ with corners and for any cohomology class $e\in H^1(B;\SS_{(P,\CL )})$, there exists a $2n$-dimensional $C^0$ manifold $(X,\CT )$ equipped with a $C^0$ local $T^n$-action whose characteristic pair and the Euler class of the orbit map are equal to $(P,\CL)$ and $e$, respectively. 
\end{thm}
\begin{proof}
If there exists a $C^0$ isomorphism $f_X\colon (X_1,\CT_1)\to (X_2,\CT_2)$, then, by Lemma~\ref{inv1} and \ref{inv3}, $f_X$ induces an isomorphism $f_{P_X}\colon (P_{X_1},\CL_{X_1})\to (P_{X_2},\CL_{X_2})$ such that $f_{P_X}^*e(X_2,\CT_2 )=e(X_1,\CT_1)$. 

Conversely, suppose that there is an isomorphism $f_P\colon (P_{X_1},\CL_{X_1})\to (P_{X_2},\CL_{X_2})$ such that $f_P^*e(X_2,\CT_2 )=e(X_1,\CT_1)$. We denote by $f_B\colon B_{X_1}\to B_{X_2}$ the stratification preserving homeomorphism which is induced by $f_P$. Let $\{ (U_{\beta}^{X_1},\varphi^{X_1}_{\beta})\} \in \CT_1$, $\{ (U_{\alpha}^{X_2},\varphi^{X_2}_{\alpha})\} \in \CT_2$ be weakly standard atlases of $X_1$, $X_2$ and $\{ (U_{\beta}^{B_1},\varphi^{B_1}_{\beta})\}$, $\{ (U_{\alpha}^{B_2},\varphi^{B_2}_{\alpha})\}$ the atlases of $B_{X_1}$, $B_{X_2}$ induced by $\{ (U_{\beta}^{X_1},\varphi^{X_1}_{\beta})\}$, $\{ (U_{\alpha}^{X_2},\varphi^{X_2}_{\alpha})\}$ such that for each $\alpha$, $\beta$, we can take $C^0$ isomorphisms $h_{\beta}^1\colon \mu_{X_1}^{-1}(U_{\beta}^{B_1})\to \mu_{X_{(P_{X_1},\CL_{X_1})}}^{-1}(U_{\beta}^{B_1})$, $h_{\alpha}^2\colon \mu_{X_2}^{-1}(U_{\alpha}^{B_2})\to \mu_{X_{(P_{X_2},\CL_{X_2})}}^{-1}(U_{\alpha}^{B_2})$ which satisfy \eqref{theta} for $X_1$, $X_2$, respectively. By the assumption $f_P^*e(X_2,\CT_2 )=e(X_1,\CT_1)$, by replacing $\{ (U_{\beta}^{X_1},\varphi^{X_1}_{\beta})\}$ by a refinement if necessary, we may assume that there exists a \v{C}ech zero-cochain $\{ \theta_{\beta}\}$ on $\{ (U_{\beta}^{B_1},\varphi^{B_1}_{\beta})\}$ with values in $\SS_{X_{(P_{X_1},\CL_{X_1})}}$ such that 
\begin{equation}\label{condi2-T}
\theta^{X_1}_{\beta_0\beta_1}(b)=f_{X_{(P,\CL )}}^{-1}\left(\theta^{X_2}_{\alpha_0\alpha_1}\left(f_B(b)\right)\right)\theta_{\beta_1}(b)\theta_{\beta_0}(b)^{-1}
\end{equation}
for $b\in U_{\beta_0\beta_1}^{B_1}\cap f_B^{-1}(U_{\alpha_0\alpha_1}^{B_2})$. For each nonempty overlap $U_{\beta}^{B_1}\cap f_B^{-1}(U_{\alpha}^{B_2})$ we define the $C^0$ isomorphism $f_{\alpha \beta}\colon \mu_{X_1}^{-1}(U_{\beta}^{B_1}\cap f_B^{-1}(U_{\alpha}^{B_2}))\to \mu_{X_2}^{-1}(f_B(U_{\beta}^{B_1})\cap U_{\alpha}^{B_2})$ by 
\begin{equation*}
f_{\alpha \beta}(x):=(h^2_{\alpha})^{-1}\circ f_{X_{(P,\CL )}}\left( \theta_{\beta}(\mu_{X_1}(x)) h^1_{\beta}(x)\right) . 
\end{equation*}
By using \eqref{condi2-T} and \eqref{theta} for $X_1$ and $X_2$, for any $x\in \mu_{X_1}^{-1}(U_{\beta_0\beta_1}^{B_1}\cap f_B^{-1}(U_{\alpha_0\alpha_1}^{B_2}))$ with $b=\mu_{X_1}(x)$   
\[
\begin{split}
f_{\alpha_0\beta_0}(x)=&(h^2_{\alpha_0})^{-1}\circ f_{X_{(P,\CL )}}\left( \theta_{\beta_0}(b) h^1_{\beta_0}(x)\right)\\
&=(h^2_{\alpha_1})^{-1}\left( \theta^{X_2}_{\alpha_1\alpha_0}(f_B(b)) f_{X_{(P,\CL )}}\left( \theta_{\beta_0}(b)\theta^{X_1}_{\beta_0\beta_1}(b) h^1_{\beta_1}(x)\right) \right) \\
&=(h^2_{\alpha_1})^{-1}\left( \theta^{X_2}_{\alpha_1\alpha_0}(f_B(b)) f_{X_{(P,\CL )}}\left( \theta_{\beta_0}(b)\theta^{X_1}_{\beta_0\beta_1}(b)\right) f_{X_{(P,\CL )}}\left( h^1_{\beta_1}(x)\right) \right) \\
&=(h^2_{\alpha_1})^{-1}\left( f_{X_{(P,\CL )}}\left( \theta_{\beta_1}(b)\right)f_{X_{(P,\CL )}}\left( h^1_{\beta_1}(x)\right) \right) \\
&=(h^2_{\alpha_1})^{-1}\circ f_{X_{(P,\CL )}}\left( \theta_{\beta_1}(b) h^1_{\beta_1}(x)\right) \\
&=f_{\alpha_1\beta_1}(x), 
\end{split}
\]
where we used the equality $\theta_{\alpha_1\alpha_0}^{X_2}(f_B(b))=\theta_{\alpha_0\alpha_1}^{X_2}(f_B(b))^{-1}$ and the fact that $f_{X_{(P,\CL )}}\colon$ $X_{(P_1, \CL_1)}\to X_{(P_2, \CL_2)}$ is fiberwise group isomorphism (see Remark~\ref{group-str}). Then we can patch the local $C^0$ isomorphisms $f_{\alpha \beta}$ together to obtain the required $C^0$ isomorphism $f_X$. 

Suppose that $(P, \CL)$ is a characteristic pair on an $n$-dimensional $C^0$ manifold $B$ with corners and $e$ is an element of $H^1(B;\SS_{(P,\CL )})$. We take a representative $\{ \theta_{\alpha \beta}\}$ of $e$ on a sufficiently small open cover $\{ U_{\alpha}\}$ of $B$. Then we can construct a new $2n$-dimensional $C^0$ manifold $(X,\CT )$ equipped with a $C^0$ local $T^n$-action by setting 
\[
X:=\left( \coprod_{\alpha}\mu_{X_{(P,\CL )}}^{-1}(U_{\alpha}^B)\right) /\sim , 
\]
where $x_{\alpha}\in \mu_{X_{(P,\CL )}}^{-1}(U_{\alpha}^B)\sim x_{\beta}\in \mu_{X_{(P,\CL )}}^{-1}(U_{\beta}^B)$ if and only if $\mu_{X_{(P,\CL )}}(x_{\alpha})=\mu_{X_{(P,\CL )}}(x_{\beta})$ and $x_{\alpha}=\theta_{\alpha \beta}(\mu_{X_{(P,\CL )}}(x_{\alpha}))x_{\beta}$. It is easy to see that its characteristic pair and Euler class of the orbit map are equal to $(P,\CL)$ and $e$, respectively.  
\end{proof}

Next we focus on the case of locally standard torus actions. We take notice that if a manifold $X$ is equipped with a locally standard torus action, then, $P_X$ is the product bundle $P_X=B_X\times \Aut (T^n)$. In this case, we can obtain the following theorem. 
\begin{thm}
For $i=1$, $2$, let $X_i$ be a $2n$-dimensional manifold equipped with a locally standard $T^n$-action and $\rho$ an automorphism of $T^n$. $X_1$ and $X_2$ are $\rho$-equivariantly homeomorphic if and only if there exists a stratification preserving homeomorphism $f_B\colon B_{X_1}\to B_{X_2}$ such that the bundle isomorphism $f_P \colon P_{X_1}\to P_{X_2}$ defined by $f_P:=f_B\times \rho$ sends $\CL_{X_1}$ isomorphically to $\CL_{X_2}$ and $f_P^*e(X_2,\CT_2 )=e(X_1,\CT_1)$. Moreover, for any characteristic pair $(P, \CL)$, where $P$ is trivial, on an $n$-dimensional $C^0$ manifold $B$ with corners and for any cohomology class $e\in H^1(B;\SS_{(P,\CL )})$, there exists a $2n$-dimensional $C^0$ manifold equipped with a locally standard $T^n$-action whose characteristic pair and the Euler class of the orbit map are equal to $(P,\CL)$ and $e$, respectively. 
\end{thm}
The proof is almost same as that of Theorem~\ref{classification}. In particular, by putting $\rho =id_{T^n}$, we obtain the following corollary. 
\begin{cor}\label{classification-locallystandard}
Locally standard torus actions are classified by the characteristic bundle and the Euler class of the orbit map up to equivariant homeomorphisms. 
\end{cor}
\begin{rem}
Corollary~\ref{classification-locallystandard} is a generalization of the topological classification of quasi-toric manifolds by Davis and Januszkiewicz~\cite{DJ} and of effective $T^2$-actions on four-dimensional manifolds without nontrivial finite stabilizers by Orlik and Raymond~\cite{OR}. 
\end{rem}

\section{The symplectic case}\label{symp}
In this section we identify $\Lambda$ with $\Z^n$ by using the fixed decomposition $T^n= (S^1)^n$. By using the natural isomorphism $\Aut (T^n)\to \GL (\Lambda)$ together with this fact, we also identify $\Aut (T^n)$ with $\GL_n(\Z )$. 
\subsection{Integral affine structures}
Let $B$ be an $n$-dimensional $C^{\infty}$ manifold with corners and $\{ (U_{\alpha}^B, \varphi^B_{\alpha})\}$ an atlas of $B$ consisting of coordinate neighborhoods modeled on open subsets of $\R^n_{\ge 0}$. 
\begin{defn}\label{integral-affine-structure}
$\{ (U_{\alpha}^B, \varphi^B_{\alpha})\}$ is called an {\it integral affine structure} on $B$ if on each nonempty overlap $U_{\alpha \beta}^B\neq \emptyset$, there exists an element $A_{\alpha \beta}$ of $\GL_n(\Z )$ and there also exists a constant $c_{\alpha \beta}\in \R^n$ such that the overlap map $\varphi^B_{\alpha \beta}\colon \varphi^B_{\beta}(U_{\alpha \beta}^B)\to \varphi^B_{\alpha}(U_{\alpha \beta}^B)$ is of the form 
\begin{equation}\label{affinetrans}
\varphi^B_{\alpha \beta}(\xi )=A_{\alpha \beta}(\xi )+c_{\alpha \beta}.  
\end{equation}
\end{defn}

First we show the following lemma. 
\begin{lem}
For an $n$-dimensional $C^{\infty}$ manifold $B$ equipped with an integral affine structure $\{ (U_{\alpha}^B, \varphi^B_{\alpha})\}$, there is a characteristic pair on $B$ associated with $\{ (U_{\alpha}^B, \varphi^B_{\alpha})\}$. 
\end{lem}
\begin{proof}
By definition, the structure group of the cotangent bundle $T^*B$ of $B$ is reduced to $\GL_n(\Z )$. We denote by $\pi_P\colon P\to B$ the frame bundle of $T^*B$ with structure group $\GL_n(\Z )$ and also denote by $\pi_{\Lambda}\colon \Lambda_P\to B$ the associated $\Lambda$-bundle of $\pi_P\colon P\to B$. Let $(U_{\alpha}^B,\varphi^B_{\alpha})$ be a coordinate neighborhood in $\{ (U_{\alpha}^B, \varphi^B_{\alpha})\}$ with $U_{\alpha}^B\cap \CS^{(n-1)}B_X\neq \emptyset$. We may assume that the intersection $U_{\alpha}^B\cap \CS^{(n-1)}B_X$ is connected. (Otherwise, we may consider componentwise.) Then there exists a unique primitive vector $u_{\alpha}\in \Lambda$ such that $\varphi^B_{\alpha}(U_{\alpha}^B)$ and $\varphi^B_{\alpha}(U_{\alpha}^B\cap \CS^{(n-1)}B_X)$ are contained in the upper half space $\{ \xi \in \R^n \colon \< \xi , u_{\alpha}\> \ge 0\}$ and the hyperplane $\{ \xi \in \R^n \colon \< \xi , u_{\alpha}\> =0\}$ determined by $u_{\alpha}$, respectively. Suppose that $(U_{\beta}^B, \varphi^B_{\beta})$ is another coordinate neighborhoods satisfying the above conditions and the intersection $U_{\alpha \beta}^B\cap \CS^{(n-1)}B_X$ is nonempty. Let $u_{\beta}\in \Lambda$ be the corresponding primitive vector. Since the overlap map $\varphi^B_{\alpha \beta}$ is of the form \eqref{affinetrans} and since $\varphi^B_{\alpha \beta}$ sends $\{ \xi \in \R^n \colon \< \xi , u_{\beta}\> \ge 0\}$ and $\{ \xi \in \R^n \colon \< \xi , u_{\beta}\> =0\}$ diffeomorphically to $\{ \xi \in \R^n \colon \< \xi , u_{\alpha}\> \ge 0\}$ and $\{ \xi \in \R^n \colon \< \xi , u_{\alpha}\> =0\}$, respectively, we can show that $u_{\alpha}=A_{\alpha\beta}^{-T}(u_{\beta})$. In particular, $u_{\alpha}$'s form a section $s$ of $\pi_{\Lambda}|_{\CS^{(n-1)}B}\colon \Lambda_P|_{\CS^{(n-1)}B}\to \CS^{(n-1)}B$. Let $\pi_{\CL}\colon \CL\to \CS^{(n-1)}B$ be the rank one sublattice bundle of $\pi_{\Lambda}|_{\CS^{(n-1)}B}\colon \Lambda_P|_{\CS^{(n-1)}B}\to \CS^{(n-1)}B$ which is fiberwise generated by $s$. By construction, $\pi_{\CL}\colon \CL\to \CS^{(n-1)}B$ is unimodular, hence the pair $(P,\CL )$ is a characteristic pair on $B$. 
\end{proof}

Note that, by the proof of this lemma, the characteristic bundle $\pi_{\CL}\colon \CL\to \CS^{(n-1)}B$ admits a section which generates $\pi_{\CL}\colon \CL\to \CS^{(n-1)}B$ fiberwise. 

In Section~\ref{sub-canonical} we constructed the canonical model from a characteristic pair. Next we show that the canonical model constructed from the characteristic pair associated with an integral affine structure is a $C^{\infty}$ manifold and admits a symplectic structure such that the orbit map is a locally toric Lagrangian fibration. Let $B$ be an $n$-dimensional $C^{\infty}$ manifold equipped with an integral affine structure $\{ (U_{\alpha}^B, \varphi^B_{\alpha})\}$ and $(P,\CL )$ the characteristic pair associated with $\{ (U_{\alpha}^B, \varphi^B_{\alpha})\}$. Let $\pi_T\colon T_P\to B$ be the associated $T^n$-bundle of $\pi_P\colon P\to B$. Then we have the following exact sequence of associated fiber bundles of $P$
\[
\xymatrix{
0\ar[r] & \Lambda_P\ar[r] & T^*B\ar[r] & T_P\ar[r] & 0. 
}
\]
As is well known, $T^*B$ is equipped with the standard symplectic structure $\om_{T^*B}$. Since the natural fiberwise action of $\Lambda_P$ on $T^*B$ preserves $\om_{T^*B}$, $\om_{T^*B}$ descends to a symplectic structure on $T_P$, which is denoted by $\om_{T_P}$, so that $\pi_T\colon (T_P,\om_{T_P})\to B$ is a nonsingular Lagrangian fibration. Moreover, we can show the following lemma. 
\begin{lem}\label{can-lag}
The canonical model $X_{(P,\CL )}$ is a $C^{\infty}$ manifold and admits a symplectic structure $\om_{X_{(P,\CL )}}$ so that $\mu_{X_{(P,\CL )}}\colon (X_{(P,\CL )},\om_{X_{(P,\CL )}})\to B$ a locally toric Lagrangian fibration. 
\end{lem}
\begin{proof}
In this case $X_{(P,\CL )}$ can be reconstructed from $(T_P,\om_{T_P})$ by using the symplectic cutting technique in the following way. For each $(U_{\alpha}^B, \varphi^B_{\alpha})$, $\overline{\varphi}^T_{\alpha}\colon \pi_T^{-1}(U_{\alpha}^B)\to \varphi^B_{\alpha}(U_{\alpha}^B)\times T^n$ denotes the composition of $\varphi^B_{\alpha}\times \id_{T^n}\colon U_{\alpha}^B\times T^n\to \varphi^B_{\alpha}(U_{\alpha}^B)\times T^n$ and the local trivialization $\varphi^T_{\alpha}$ of $T_P$ on $U_{\alpha}^B$ induced by the local trivialization $\varphi^P_{\alpha}$ of $P$ on $U_{\alpha}^B$. By the construction of $\om_{T_P}$, $\overline{\varphi}^T_{\alpha}$ preserves the symplectic structures, namely, $(\overline{\varphi}^T_{\alpha})^*\om_{\R^n\times T^n}=\om_{T_P}$. We put 
\[
I_{\alpha}:=\{ i\in \{ 1, \ldots ,n\} \colon \varphi^B_{\alpha}(U_{\alpha}^B)\cap \{ \xi \in \R^n_{\ge 0} \colon \xi_i=0\} \neq \emptyset \}
\] 
and 
\[
T_{I_{\alpha}}:=\{ u\in T^n\colon u_i=1\ \text{unless}\ i\in I_{\alpha}\} . 
\]
Suppose that $I_{\alpha}=\{ i_1, \ldots ,i_k \}$ and $i_1<\ldots < i_k$. We take an open set $V_{\alpha}$ of $\R^n$ which satisfies $\varphi^B_{\alpha}(U_{\alpha}^B)=V_{\alpha}\cap \R^n_{\ge 0}$. $T_{I_{\alpha}}$ acts on $V_{\alpha}\times T^n\times \C^k$ by 
\[
u\cdot (\xi ,v, z):=\left( \xi ,uv, (u_{i_1}^{-1}z_1, \ldots ,u_{i_k}^{-1}z_k)\right) . 
\]
This action is Hamiltonian with respect to the symplectic structure $\om_{\R^n\times T^n}\oplus \om_{\C^k}$ and a moment map $\Phi_{\alpha}\colon (V_{\alpha}\times T^n\times \C^k, \om_{\R^n\times T^n}\oplus \om_{\C^k})\to \R^k$ is written by
\[
\Phi_{\alpha}(\xi ,v, z)=(\xi_{i_1}-\abs{z_1}^2,\ldots ,\xi_{i_k}-\abs{z_k}^2). 
\]
We take the symplectic quotient $\Phi_{\alpha}^{-1}(0)/T_{I_{\alpha}}$ of the $T_{I_{\alpha}}$-action at $0\in \R^k$. We denote by $\om_{\alpha}$ the reduced symplectic structure on $\Phi_{\alpha}^{-1}(0)/T_{I_{\alpha}}$. The natural action of $T^n$ on $V_{\alpha}\times T^n\times \C^k$ by the multiplication of the second factor induces a Hamiltonian $T^n$-action on $(\Phi_{\alpha}^{-1}(0)/T_{I_{\alpha}}, \om_{\alpha})$ with moment map $\mu_{\alpha}([\xi ,u,z])=\xi$. Note that the image of $\mu_{\alpha}$ is equal to $\varphi^B_{\alpha}(U_{\alpha}^B)$. It is obvious that this $T^n$-action on $(\Phi_{\alpha}^{-1}(0)/T_{I_{\alpha}}, \om_{\alpha})$ with moment map $\mu_{\alpha}$ is naturally identified with the standard $T^n$-action on $(\mu_{\C^n}^{-1}(\varphi^B_{\alpha}(U_{\alpha}^B)),\om_{\C^n})$ with moment map $\mu_{\C^n}$. By the construction of $X_{(P,\CL )}$ we can also show that $\overline{\varphi}^T_{\alpha}$ induces a natural homeomorphism $\varphi^{X_{(P,\CL )}}_{\alpha}\colon \mu_{X_{(P,\CL )}}^{-1}(U_{\alpha}^B)\to \Phi_{\alpha}^{-1}(0)/T_{I_{\alpha}}$ which satisfies $\mu_{\alpha}\circ \varphi^{X_{(P,\CL )}}_{\alpha}=\varphi^B_{\alpha}\circ \mu_{X_{(P,\CL )}}$. 

Suppose that $(U_{\alpha}^B, \varphi^B_{\alpha})$ and $(U_{\beta}^B, \varphi^B_{\beta})$ are coordinate neighborhoods with nonempty overlap $U_{\alpha \beta}^B$. To prove the lemma, it is sufficient to show that the overlap map $\varphi^{X_{(P,\CL )}}_{\alpha \beta}:=\varphi^{X_{(P,\CL )}}_{\alpha}\circ (\varphi^{X_{(P,\CL )}}_{\beta})^{-1}$ is a symplectomorphism from $(\mu_{\beta}^{-1}(\varphi^B_{\beta}(U_{\alpha \beta}^B)),\om_{\beta})$ to $(\mu_{\alpha}^{-1}(\varphi^B_{\alpha}(U_{\alpha \beta}^B)),\om_{\alpha})$. Since the overlap map $\overline{\varphi}^T_{\alpha \beta}:=\overline{\varphi}^T_{\alpha}\circ (\overline{\varphi}^T_{\beta})^{-1}$ is of the form $\overline{\varphi}^T_{\alpha \beta}(\xi ,u)=(A_{\alpha \beta}(\xi )+c_{\alpha \beta}, A_{\alpha \beta}^{-T}(u))$, it can be naturally extended to the symplectomorphism from $(V_{\beta}\cap (\varphi^B_{\alpha \beta})^{-1}(V_{\alpha})\times T^n, \om_{\R^n\times T^n})$ to $(\varphi^B_{\alpha \beta}(V_{\beta})\cap V_{\alpha}\times T^n, \om_{\R^n\times T^n})$, which is also denoted by the same notation $\overline{\varphi}^T_{\alpha \beta}$. Since $\overline{\varphi}^T_{\alpha \beta}$ is $A_{\alpha \beta}^{-T}$-equivariant with respect to the actions of $T_{I_{\beta}}$ and $T_{I_{\alpha}}$, it descends to a symplectomorphism from $(\mu_{\beta}^{-1}(\varphi^B_{\beta}(U_{\alpha \beta}^B)),\om_{\beta})$ to $(\mu_{\alpha}^{-1}(\varphi^B_{\alpha}(U_{\alpha \beta}^B)),\om_{\alpha})$ which is nothing but $\varphi^{X_{(P,\CL )}}_{\alpha \beta}$. This proves the lemma. 
\end{proof}

\begin{rem}
Let $B_1$ and $B_2$ be $n$-dimensional $C^{\infty}$ manifolds equipped with integral affine structures. Suppose that there exists a diffeomorphism $f_B\colon B_1\to B_2$ which preserves the integral affine structures. Then $f_B$ induces an isomorphism $f_P\colon (P_1,\CL_1)\to (P_2,\CL_2)$ between the characteristic pairs associated with the integral affine structures. Moreover, it is easy to see that the $C^0$ isomorphism $f_{X_{(P,\CL )}}\colon X_{(P_1, \CL_1)}\to X_{(P_2, \CL_2)}$ induced by $f_P$ is a fiber-preserving symplectomorphism.   
\end{rem}

\subsection{The necessary and sufficient condition}
Let $(X,\CT )$ be a $2n$-dimensional manifold equipped with a $C^{\infty}$ local $T^n$-action $\CT$. In this subsection, we investigate the condition in order that $X$ admits a symplectic structure $\om$ so that $\mu_X\colon (X,\om )\to B_X$ is a locally toric Lagrangian fibration. 
\begin{defn}
Let $\{ (U_{\alpha}^X, \varphi^X_{\alpha})\}_{\alpha \in \CA}\in \CT$ be a weakly standard $C^{\infty}$ atlas of $X$ and $\{ (U_{\alpha}^B, \varphi^B_{\alpha})\}_{\alpha \in \CA}$ the atlas of ${B_X}$ induced by $\{ (U_{\alpha}^X,$ $\varphi^X_{\alpha})\}_{\alpha \in \CA}$. For each nonempty overlap $U_{\alpha \beta}^X\neq \emptyset$, let $\rho_{\alpha \beta}\in \Aut (T^n)$ be the automorphism in \textup{(ii)} of Definition~\textup{\ref{localaction}} with respect to $\{ (U_{\alpha}^X, \varphi^X_{\alpha})\}_{\alpha \in \CA}$. We call $\{ (U_{\alpha}^B, \varphi^B_{\alpha})\}_{\alpha \in \CA}$ an {\it integral affine structure compatible with $\{ (U_{\alpha}^X,$ $\varphi^X_{\alpha})\}_{\alpha \in \CA}$} if $\{ (U_{\alpha}^B, \varphi^B_{\alpha})\}_{\alpha \in \CA}$ is an integral affine structure on $B_X$ and on each nonempty overlap $U_{\alpha \beta}^B\neq \emptyset$, the element $A_{\alpha \beta} \in \GL_n (\Z )$ in \eqref{affinetrans} is equal to $\rho_{\alpha \beta}^{-T}$. 
\end{defn}

\begin{lem}\label{integralaffine}
Suppose that there exists a symplectic structure $\om$ on $X$ and there also exists a weakly standard atlas $\{ (U_{\alpha}^X, \varphi^X_{\alpha})\}_{\alpha \in \CA}\in \CT$ of $X$ such that on each $U_{\alpha}^X$, $\varphi^X_{\alpha}$ preserves symplectic structures, namely, $\om ={\varphi^X_{\alpha}}^*\om_{\C^n}$. Then the atlas $\{ (U_{\alpha}^B, \varphi^B_{\alpha})\}_{\alpha \in \CA}$ of ${B_X}$ induced by $\{ (U_{\alpha}^X,$ $\varphi^X_{\alpha})\}_{\alpha \in \CA}$ is an integral affine structure compatible with $\{ (U_{\alpha}^X,$ $\varphi^X_{\alpha})\}_{\alpha \in \CA}$. In particular, $B_X$ becomes a $C^{\infty}$ manifold with corners. 
\end{lem}
\begin{proof}
By assumption, the restriction of $\mu_X \colon (X,\om)\to B_X$ to $B_X\setminus \partial B_X$ is a nonsingular Lagrangian fibration, and as in the proof of Proposition~\ref{locally-toric-Lag}, we can construct a symplectomorphism $\phi_{\alpha}\colon (\mu_X^{-1}(U_{\alpha}^B\setminus \partial B_X), \om )\to (\varphi^B_{\alpha}(U_{\alpha}^B\setminus \partial B_X)\times T^n, \om_{\R^n\times T^n})$ such that $\pr_1\circ \phi_{\alpha}=\varphi^B_{\alpha}\circ \mu_X$ for each $\alpha$. By applying Lemma~\ref{local-auto} to the overlap map $\phi_{\alpha \beta}:=\phi_{\alpha}\circ \phi_{\beta}^{-1}$ on each nonempty overlap $U_{\alpha \beta}^B\setminus \partial B_X$, we can show that the overlap map $\varphi^B_{\alpha \beta}$ of the base is of the form \eqref{affinetrans}. Since $U_{\alpha \beta}^B\setminus \partial B_X$ is open dense in $U_{\alpha \beta}^B$, $\varphi^B_{\alpha \beta}$ should be of the form \eqref{affinetrans} on the whole $U_{\alpha \beta}^B$. 
\end{proof}

Let $\{ (U_{\alpha}^X, \varphi^X_{\alpha})\}_{\alpha \in \CA}\in \CT$ be a weakly standard atlas of $X$. Suppose that the induced atlas $\{ (U_{\alpha}^B, \varphi^B_{\alpha})\}_{\alpha \in \CA}$ of $B_X$ is an integral affine structure compatible with $\{ (U_{\alpha}^X, \varphi^X_{\alpha})\}_{\alpha \in \CA}\in \CT$. Then the principal $\Aut (T^n)$-bundle $P_X$ is nothing but the frame bundle of $T^*B_X$ and the characteristic pair $(P_X,\CL_X)$ is nothing but the characteristic pair associated with $\{ (U_{\alpha}^B, \varphi^B_{\alpha})\}_{\alpha \in \CA}$. By Lemma~\ref{can-lag}, the canonical model $X_{(P_X,\CL_X)}$ is a $C^{\infty}$ symplectic manifold with symplectic structure $\om_{X_{(P_X,\CL_X)}}$ so that $\mu_{X_{(P_X,\CL_X)}}\colon (X_{(P_X,\CL_X)},\om_{X_{(P_X,\CL_X)}})\to B_X$ a locally toric Lagrangian fibration. In particular, the local isomorphism $h_{\alpha}\colon \mu_X^{-1}(U_{\alpha}^B)\to \mu_{X_{(P_X,\CL_X)}}^{-1}(U_{\alpha}^B)$ in Section~\ref{orbit-map-section} can be taken to be a $C^{\infty}$ isomorphism which covers the identity on each $U_{\alpha}^B$ and by \cite[Lemma 3.2]{HS} $\theta_{\alpha\beta}^X$ defined by \eqref{theta} can be also taken to be a $C^{\infty}$ local section of $T_X$ on $U_{\alpha \beta}^B$. Then a necessary and sufficient condition is given as follows. 
\begin{thm}\label{om}
Let $(X,\CT )$ be a $2n$-dimensional manifold equipped with a $C^{\infty}$ local $T^n$-action $\CT$. There exists a symplectic structure $\om$ on $X$ and there also exists a weakly standard atlas $\{ (U_{\alpha}^X, \varphi^X_{\alpha})\}_{\alpha \in \CA}\in \CT$ of $X$ such that on each $U_{\alpha}^X$, $\om ={\varphi^X_{\alpha}}^*\om_{\C^n}$ if and only if the atlas $\{ (U_{\alpha}^B, \varphi^B_{\alpha})\}_{\alpha \in \CA}$ of $B_X$ induced by $\{ (U_{\alpha}^X,$ $\varphi^X_{\alpha})\}_{\alpha \in \CA}$ is an integral affine structure compatible with $\{ (U_{\alpha}^X,$ $\varphi^X_{\alpha})\}_{\alpha \in \CA}$ and on each nonempty overlap $U_{\alpha \beta}^B$, $\theta_{\alpha \beta}^X$ is a Lagrangian section, namely, $(\theta_{\alpha \beta}^X)^*\om_{T_X}$ vanishes. 
\end{thm}

For nonsingular Lagrangian fibrations, this result is obtained by Duistermaat~\cite{Du}. See also \cite{S}, \cite{Mis}. Recently, in \cite{GaS} Gay and Symington showed the similar result for near-symplectic four-manifolds. 

\subsection{The symplectic classification}
Suppose that $(X,\CT )$ is a $2n$-dimensional $C^{\infty}$ manifold equipped with a $C^{\infty}$ local $T^n$-action $\CT$ satisfying the condition in Theorem~\ref{om}. Then $X$ is equipped with a symplectic structure $\om$ so that $\mu_X\colon (X,\om )\to B_X$ is a locally toric Lagrangian fibration, and the local sections $\theta_{\alpha \beta}^X$ define a \v{C}ech cohomology class $\lambda (X)\in H^1(B_X;\SS_{T_X}^{Lag})$ of $B_X$ with values in the sheaf $\SS_{T_X}^{Lag}$ of Lagrangian sections of $\pi_{T_X}\colon (T_X,\om_{T_X})\to B_X$. $\lambda (X)$ is called a {\it Lagrangian class} for $\mu_X\colon (X,\om )\to B_X$. Note that for two $2n$-dimensional $C^{\infty}$ manifolds $(X_1,\CT_1)$ and $(X_2,\CT_2)$ with $C^{\infty}$ local $T^n$-actions which satisfy the condition in Theorem~\ref{om}, if there exists a diffeomorphism $f_B\colon B_{X_1}\to B_{X_2}$ which preserves the integral affine structures, then $f_B$ induces an isomorphism $f_P^*\colon H^1(B_{X_2};\SS_{T_{X_2}}^{Lag})\to H^1(B_{X_1};\SS_{T_{X_1}}^{Lag})$. 

\begin{thm}[\cite{BM}]\label{symplectic-classification}
For $i=1,2$, let $(X_i,\CT_i)$ be a $2n$-dimensional $C^{\infty}$ manifold with a $C^{\infty}$ local $T^n$-action which is equipped with a symplectic structure $\om_i$ so that $\mu_{X_i}\colon (X_i,\om_i )\to B_{X_i}$ is a locally toric Lagrangian fibration. $\mu_{X_1}\colon (X_1,\om_1 )\to B_{X_1}$ and $\mu_{X_2}\colon (X_2,\om_2 )\to B_{X_2}$ are fiber-preserving symplectomorphic if and only if there exists a diffeomorphism $f_B\colon B_{X_1}\to B_{X_2}$ which preserves the integral affine structures such that $f_P^*\lambda (X_2)=\lambda (X_1)$. Moreover, for any $n$-dimensional $C^{\infty}$ manifold $B$ with corners equipped with an integral affine structure and for any cohomology class $\lambda \in H^1(B;\SS_{T_P}^{Lag})$, there exists a locally toric Lagrangian fibration whose integral affine structure and the Lagrangian class are the given ones. Here $\SS_{T_P}^{Lag}$ is the sheaf of Lagrangian sections of $\pi_T\colon (T_P,\om_{T_P})\to B$. 
\end{thm}
The proof is similar to Theorem~\ref{classification}.

\section{Topology}\label{topology}
Let $(X,\CT )$ be a $2n$-dimensional manifold with a local $T^n$-action. We assume that $X$ is connected (hence $B$ is also connected) and $e(X,\CT )$ vanishes. In this section, we investigate some topological invariants for $X$. 
\subsection{Fundamental groups}
In this subsection, we investigate the fundamental group of X. We fix a point $b_0$ in the interior $\CS^{(n)}B_X$ of $B_X$ and also fix points $x_0\in {\mu_X}^{-1}(b_0)$ and $t_0\in \pi^{-1}_{T_X}(b_0)$ which satisfy $\nu (t_0)=x_0$ as base points of $X$ and $T_X$, respectively, where $\nu\colon T_X\to X$ is the map defined by \eqref{nu}. Comparing the fundamental group of $T_X$ with that of $X$ by using the homomorphism induced from $\nu \colon T_X\to X$, we have the following result. 
\begin{thm}\label{pi_1}
Suppose that $\CS^{(0)}B_X$ is nonempty. Then $\mu_X$ induces an isomorphism ${\mu_X}_*\colon \pi_1(X, x_0) \cong \pi_1(B_X, b_0)$ of fundamental groups. 
\end{thm}
\begin{proof}
Since the structure group of $\pi_{T_X}\colon T_X\to B_X$ is $\Aut (T^n)$ it admits a section $s_T$. $s_T$ defines a section $s_X$ of $\mu_X$ by $s_X:=\nu \circ s_T$. Now we have the following commutative diagram of split short exact sequences for fundamental groups: 
\[
\xymatrix{
1 \ar[r] & \pi_1(\pi_{T_X}^{-1}(b_0),t_0)\ar[r]^{\iota_*} \ar[d]^{\kappa} & \pi_1(T_X, t_0) \ar[r]^{{\pi_{T_X}}_*} \ar[d]^{\nu_*} & \pi_1(B_X ,b_0) \ar[r] \ar[d]^{\id_{\pi_1(B_X ,b_0)}}& 1 \ \ \text{(exact)}\\
1 \ar[r] & \ker {\mu_X}_*\ar[r] & \pi_1(X, x_0) \ar[r]^{{\mu_X}_*} & \pi_1(B_X ,b_0) \ar[r] & 1 \ \ 
\text{(exact),} 
}
\]
where $\iota \colon \pi_{T_X}^{-1}(b_0)\hookrightarrow T_X$ is the natural inclusion and $\kappa \colon \pi_1(\pi_{T_X}^{-1}(b_0),t_0)\to \ker {\mu_X}_*$ is the homomorphism induced by the restriction $\nu |_{\pi_{T_X}^{-1}(b_0)} :\pi_{T_X}^{-1}(b_0)\to \mu^{-1}(b_0)$ of $\nu$ to $\pi_{T_X}^{-1}(b_0)$. First we claim that $\kappa$ is surjective. Note that it is equivalent to the surjectivity of $\nu_*$ since ${s_X}_*=\nu_*\circ {s_T}_*$. Since $b_0$ is in the interior of $B_X$, $\nu |_{\pi_{T_X}^{-1}(b_0)}$ is a homeomorphism which sends $t_0$ to $x_0$. Then it is sufficient to show that every element of $\ker {\mu_X}_*$ is represented by a loop in ${\mu_X}^{-1}(b_0)$. Let $\alpha \in \ker {\mu_X}_*$ and take its representative $a'\colon S^1 \to X$ with $a'(1)=x_0$. Then the map ${\mu_X} \circ a'\colon S^1\to B_X$ is homotopic to the constant map $b_0$. If necessary, by replacing a representative of $\alpha$, we can take a contractible open set $U$ located in the interior of $B_X$ such that the image of $a'$ is included in ${\mu_X}^{-1}(U)$. Since $U$ is in the interior of $B_X$, by Proposition~\ref{T^n-bundle}, ${\mu_X}^{-1}(b_0)$ is a deformation retract of ${\mu_X}^{-1}(U)$. We take a deformation retraction $h:I\times {\mu_X}^{-1}(U)\to {\mu_X}^{-1}(U)$ which satisfies $h(0, \cdot)=id_{{\mu_X}^{-1}(U)}$ and $h(1, \cdot): {\mu_X}^{-1}(U)\to {\mu_X}^{-1}(b_0)$, where $I$ implies the unit interval $[0,1]$. Then the map $f(s,u):=h(s, a'(u ))$ is the homotopy which connects $a'$ to the loop $a(u):=h(1, a'(u))$ in ${\mu_X}^{-1}(b_0)$ with $a(1)=x_0$. This proves the claim. 

To prove the theorem, it is sufficient to show that $\kappa$ is the constant map whose value is the unit element. By assumption, there exists a point $b_1$ in $\CS^{(0)}B_X$. The fiber ${\mu_X}^{-1}(b_1)$ consists of only one point which is denoted by $x_1$. We take a path $\gamma \colon I\to B_X$ with $\gamma (0)=b_0$ and $\gamma (1)=b_1$. 
(Since $B_X$ is a connected $C^0$ manifold, we can always take such a path.) Since $\pi_{T_X}\colon T_X\to B_X$ is a fiber bundle, there exists a continuous map $\widetilde{\gamma}\colon I\times \pi_{T_X}^{-1}(b_0)\to T_X$ which satisfies 
the following:
\[
\pi_{T_X}\circ \widetilde{\gamma}(s, t)=\gamma (s),\ \ \widetilde{\gamma}(0, t)=t
\]
for $s\in I$ and $t\in \pi_{T_X}^{-1}(b_0)$. $\widetilde{\gamma}$ is unique up to homotopy. For any $\alpha \in \pi_1(\pi_{T_X}^{-1}(b_0),t_0)$ and its representative $a\colon S^1\to \pi_{T_X}^{-1}(b_0)$ with $a(1)=t_0$, we define the map $\overline{a}\colon I\times S^1\to X$ by $\overline{a}(s, u):=\nu (\widetilde{\gamma}(s, a(u)))$. $\overline{a}$ satisfies 
\[
\overline{a}(0,u )=\nu \circ a(u),\ \ \overline{a}(1, u )\equiv x_1 
\]
for any $u \in S^1$ and $s\in I$. Then, $\overline{a}$ descends to the map from a two-dimensional closed disc $D$ which bounds $\nu \circ a$. This implies that $\kappa (\alpha )$ is the unit element.  
\end{proof}

\subsection{Cohomology groups}
In this subsection, we give the method for computing cohomology groups of $X$. 

Suppose that $B_X$ is equipped with a structure of a CW complex so that each $p$-cell $e^{(p)}$ is contained in some ${\mathcal S}^{(k)}B_X$ of $B_X$. Let $B_X^{(p)}$ be the $p$-skeleton and $X^{(p)}:={\mu_X}^{-1}(B_X^{(p)})$ its preimage by ${\mu_X}$. We put 
\[
A^{p,q}:=H^{p+q}(X, X^{(p-1)};\Z ),\ E^{p,q}:=H^{p+q}(X^{(p)}, X^{(p-1)};\Z ). 
\]
$i\colon A^{p,q}\to A^{p-1, q+1}$ and $j \colon A^{p,q}\to E^{p,q}$ denote the maps induced from inclusions $(X,X^{(p-2)})\subset (X,X^{(p-1)})$ and $(X^{(p)},X^{(p-1)})\subset (X, X^{(p-1)})$, respectively, and $k:E^{p,q}\to A^{p+1,q}$ also denotes the connecting homomorphism of the exact sequence of the triple $(X,X^{(p)}$, $X^{(p-1)})$. We consider the cohomology Leray spectral sequence of ${\mu_X} :X\to B_X$, namely, the spectral sequence $\{ (E_X)^{p,q}_r, d^X_r \}$ associated with the exact couple 
\[
\displaystyle 
\xymatrix{
\bigoplus_{p,q}A^{p,q}\ar[rr]^{i} & & \bigoplus_{p,q}A^{p,q}\ar[dl]^{j} \\
 & \bigoplus_{p,q}E^{p,q}\ar[ul]^k.  & 
}
\]
For the spectral sequence of cohomology groups, see \cite{H}. For a $p$-cell $e^{(p)}$, we denote by $c^{(p)}$ the center of $e^{(p)}$ and also denote the restriction of $\nu:T_X\to X$ to $\pi^{-1}_{T_X}(c^{(p)})$ by $\nu_{c^{(p)}}:\pi^{-1}_{T_X}(c^{(p)})\to {\mu_X}^{-1}(c^{(p)})$. Let $\left( C^*(B_X ; \SH_T^q), \delta^* \right)$ be the cochain complex of the CW complex $B_X$ with the Serre local system $\SH_T^q$ of the $q$th cohomology with $\Z$-coefficient for the fiber bundle $\pi_{T_X}: T_X\to B_X$. We denote by $C^p(B_X ; \SH_X^q)$ the subset of $C^p(B_X ; \SH_T^q)$ consisting of those cochains which take values in the image $\nu_{c^{(p)}_{\lambda}}^*\left( H^q({\mu_X}^{-1}(c^{(p)}_{\lambda});\Z )\right)$ of $H^q({\mu_X}^{-1}(c^{(p)}_{\lambda});\Z )$ by 
$\nu_{c^{(p)}_{\lambda}}^*$ for each $p$-cell $e^{(p)}_{\lambda}$. 
\begin{thm}\label{H-spectral}
For any $q$, $C^*(B_X ; \SH_X^q)$ is a subcomplex of $\left( C^*(B_X ; \SH_T^q), \delta \right)$. We denote its $p$th cohomology group by $H^p(B_X ;\SH^q_X)$. Then we have the isomorphisms 
\begin{align}
(E_X)_1^{p, q}&\cong C^p(B_X ; \SH^q_X), \ \ \ 
(E_X)_2^{p, q}\cong H^p(B_X ; \SH^q_X), \nonumber \\ 
(E_X)_{\infty}^{p,q} &=F^pH^{p+q}(X; \Z )/F^{p+1}H^{p+q}(X;\Z ), \nonumber 
\end{align}
where $F^lH^k(X; \Z )$ is the image of the natural map $H^k(X, X^{(l-1)};\Z )\to H^k(X;\Z )$. 
\end{thm}
\begin{proof}
The last equality is a consequence of the general theory of a spectral sequence. We prove the first two isomorphisms. For each $p$-cell $e^{(p)}$ of $B_X$, let $\tau \colon D^p\to B_X$ be the characteristic map of $e^{(p)}$ with $\tau (0)=c^{(p)}$, where $D^p (\subset \R^p)$ is the $p$-dimensional closed ball centered at $0$. We put $T_X^{(p)}:=\pi_{T_X}^{-1}(B_X^{(p)})$. Suppose that $e^{(p)}$ is contained in $\CS^{(k)}B_X$. By the same way in the construction of the characteristic bundle, we can see that there exists a rank $n-k$ subtorus bundle $\pi_{S_X}\colon S_X\to \CS^{(k)}B_X$ of $\pi_{T_X}|_{\CS^{(k)}B_X}\colon T_X|_{\CS^{(k)}B_X}\to \CS^{(k)}B_X$ such that the restriction of ${\mu_X} \colon X\to B_X$ to $\CS^{(k)}B_X$ is obtained as the quotient bundle $T_X/S_X\to \CS^{(k)}B_X$. Since the structure group of $\pi_{T_X}\colon T_X\to B_X$ is discrete, there exists a bundle map $\widetilde{\tau}_T\colon (D^p, \partial D^p)\times \pi_{T_X}^{-1}(c^{(p)})\to (T_X^{(p)}, T_X^{(p-1)})$ which covers $\tau$ such that 
\[
\widetilde{\tau}_T|_{\{ 0\}\times \pi_{T_X}^{-1}(c^{(p)})}=\id_{\pi_{T_X}^{-1}(c^{(p)})}
\]
and $\widetilde{\tau}_T$ sends $D^p\times \pi_{S_X}^{-1}(c^{(p)})$ to $S_X$. Hence $\widetilde{\tau}_T$ induces a continuous map $\widetilde{\tau}_X\colon (D^p,$ $\partial D^p)\times {\mu_X}^{-1}(c^{(p)})\to (X^{(p)}, X^{(p-1)})$ such that 
\[
\widetilde{\tau}_X|_{\{ 0\}\times {\mu_X}^{-1}(c^{(p)})}=\id_{{\mu_X}^{-1}(c^{(p)})} 
\]
and the following diagram commutes: 
\begin{equation}\label{widetilde{tau}}
\xymatrix{
D^p\times \pi_{T_X}^{-1}(c^{(p)})\ar[rr]^{\id_{D^p}\times \nu_{c^{(p)}}}\ar[dr]^{\pr_1}\ar[ddr]_{\widetilde{\tau}_T} &   & D^p\times {\mu_X}^{-1}(c^{(p)})\ar[dl]^{\pr_1}\ar[ddr]^{\widetilde{\tau}_X}  & \\
 & D^p\ar[ddr]^{\tau}|\hole &  & \\
 & T_X\ar[rr]^{\nu}\ar[dr]_{\pi_{T_X}} &  & X\ar[dl]^{{\mu_X}} \\
 &  &  B_X. & 
}
\end{equation}
Let $\{ (E_T)^{p,q}_r, d^T_r\}$ be the cohomology Serre spectral sequence of the torus bundle $\pi_{T_X}\colon T_X\to B_X$. By using $\widetilde{\tau}_T$, the excision isomorphism, and the K\"unneth formula, we have the isomorphisms 
\begin{align} 
(E_T)_1^{p,q}&=H^{p+q}((T_X)^{(p)}, (T_X)^{(p-1)}; \Z ) \nonumber \\
&\cong \sum_{\lambda}H^{p+q}(\pi_{T_X}^{-1}(\overline{e^{(p)}_{\lambda}}), 
\pi_{T_X}^{-1}(\overline{e^{(p)}_{\lambda}}-e^{(p)}_{\lambda});\Z ) \nonumber \\
&\cong \sum_{\lambda}H^{p+q}((D^p_{\lambda}, 
\partial D^p_{\lambda})\times \pi_{T_X}^{-1}(c^{(p)}_{\lambda}); \Z )\nonumber \\
&\cong \sum_{\lambda}H^p(D^p_{\lambda}, 
\partial D^p_{\lambda} ;\Z )\otimes H^q(\pi_{T_X}^{-1}(c^{(p)}_{\lambda});\Z ) 
\left( =C^p(B_X ;\SH^q_T)\right) , \nonumber 
\end{align}
\begin{align}
(E_X)_1^{p,q}&=H^{p+q}(X^{(p)}, X^{(p-1)}; \Z ) \nonumber \\
&\cong \sum_{\lambda}H^{p+q}({\mu_X}^{-1}(\overline{e^{(p)}_{\lambda}}), 
{\mu_X}^{-1}(\overline{e^{(p)}_{\lambda}}-e^{(p)}_{\lambda});\Z ) \nonumber \\
&\cong \sum_{\lambda}H^{p+q}((D^p_{\lambda}, 
\partial D^p_{\lambda})\times {\mu_X}^{-1}(c^{(p)}_{\lambda}); \Z )\nonumber \\
&\cong \sum_{\lambda}H^p(D^p_{\lambda}, 
\partial D^p_{\lambda} ;\Z )\otimes H^q({\mu_X}^{-1}(c^{(p)}_{\lambda});\Z ), \nonumber 
\end{align}
where the sum runs over all $p$-dimensional cells $e^{(p)}_{\lambda}$. By the naturality of these isomorphisms and the commutativity of the diagram \eqref{widetilde{tau}}, the following diagram commutes
\begin{equation}\label{induced-nu}
\xymatrix@C=0pt{
(E_X)^{p,q}_1\cong \ar[d]_{\nu^*} & \sum_{\lambda}H^p(D^p_{\lambda}, 
\partial D^p_{\lambda} ;\Z )\otimes H^q({\mu_X}^{-1}(c^{(p)}_{\lambda});\Z )
\ar[d]^{\sum_{\lambda}\id^*_{D^p} \otimes \nu^*_{c^{(p)}_{\lambda}}} \\
(E_T)^{p,q}_1\cong & 
\sum_{\lambda}H^p(D^p_{\lambda}, 
\partial D^p_{\lambda} ;\Z )\otimes H^q(\pi_{T_X}^{-1}(c^{(p)}_{\lambda});\Z ) . 
}
\end{equation}
Moreover, it is easy to see that the homomorphism 
$\nu_{c^{(p)}_{\lambda}}^*\colon H^*({\mu_X}^{-1}(c^{(p)}_{\lambda}) ;\Z )\to H^*(\pi_{T_X}^{-1}(c^{(p)}_{\lambda}) ;\Z )$ 
induced by $\nu_{c^{(p)}}$ is injective. This implies $(E_X)_1^{p, q}\cong C^p(B_X ; \SH^q_X)$. 
The $E_1$-term $\{ (E_T)^{p,q}_1, d^T_1\}$ of the Serre spectral sequence is isomorphic to the CW complex 
$\left( C^p(B_X ; \SH_T^q), \delta \right)$ with the Serre local system $\SH_T^q$ for the torus bundle $\pi_{T_X}: T_X\to B_X$. This fact together with the naturality of the maps in the spectral sequences shows $(E_X)_2^{p, q}\cong H^p(B_X ; \SH^q_X)$. 
\end{proof}
\begin{rem}\label{rem-spectral}
(1) For $q=0$, it is easy to see that $(E_X)_2^{p,0}\cong H^p(B_X;\Z )$. Moreover $(E_X)_1^{p,q}=0$, if $q$ or $p$ is greater than half the dimension of $X$. \\
(2) If $B_X$ is an oriented surface with $\partial B_X \neq \emptyset$, we can take a cell decomposition of $B_X$ so that all zero cells are included in $\partial B_X$. In this case, the Leray spectral sequence $\{ (E_X)^{p,q}_r, d^X_r \}$ degenerates at $E_2$-term. In fact, $\partial B_X\neq \emptyset$ implies $(E_X)_2^{2, 0}\cong H^2(B_X; \Z )=0$, and since $e^{(0)}_{\lambda}\in \partial B_X$, the fiber ${\mu_X}^{-1}(e^{(0)}_{\lambda})$ 
of ${\mu_X}$ on $e^{(0)}_{\lambda}$ is diffeomorphic to the torus whose dimension is equal or less than one. Then $(E_X)_2^{0,2}\cong (E_X)_1^{0,2}\cong C^0(B_X; \SH_X^2)=0$. 
\end{rem}

\begin{cor}\label{fixed point theorem}
Assume that $B_X$ is a finite CW complex. Then the Euler characteristic $\chi (X)$ is equal to the cardinality of $\CS^{(0)}B_X$. 
\end{cor}
\begin{proof}
Let us consider the rational coefficient cohomology Leray spectral sequence $\{ (E_X)^{p,q}_r, d^X_r \}$ of the map ${\mu_X} :X\to B_X$. Define 
\[
\chi ((E_X)_r)=\sum_{p,q}(-1)^{p+q}\dim_{\Q}(E_X)^{p,q}_r. 
\]
Since $(E_X)^{p,q}_1=C^p(B_X ;\SH^q_X)$, 
\begin{align}
\chi ((E_X)_1)&=\sum_{p,q}(-1)^{p+q}\dim_{\Q}(E_X)^{p,q}_1 \nonumber \\
&=\sum_{p,q}(-1)^{p+q}\dim_{\Q}C^p(B_X ;\SH^q_X) \nonumber \\
&=\sum_{p,q}(-1)^{p+q}\sum_{\lambda}\dim_{\Q}H^q({\mu_X}^{-1}(c^{(p)}_{\lambda}) ;\Q ) 
\nonumber \\
&=\sum_{p}\sum_{\lambda}(-1)^p\chi ({\mu_X}^{-1}(c^{(p)}_{\lambda})) , \nonumber
\end{align}
where the summation $\sum_{\lambda}$ runs over all $p$-cells. By the construction of ${\mu_X} \colon X\to B_X$, the fiber ${\mu_X}^{-1}(c^{(p)}_{\lambda})$ is homeomorphic to a compact torus of dimension equal or less than $n$, which is zero-dimensional if and only if $c^{(p)}_{\lambda}\in {\mathcal S}^{(0)}B_X$. Moreover, the assumption of the cell decomposition of $B_X$ implies that $c^{(p)}_{\lambda}\in {\mathcal S}^{(0)}B_X$ if and only if $p=0$. Then $\chi ((E_X)_1)$ is equal to the cardinality of ${\mathcal S}^{(0)}B_X$. On the other hand, it is easy to see 
\[
\chi ((E_X)_r)=\chi ((E_X)_1) 
\]
for any $r$, and by (1) in Remark~\ref{rem-spectral}, for any $r$ greater than $n$, 
\[
(E_X)_{\infty}=(E_X)_r. 
\]
Moreover, from 
$(E_X)_{\infty}^{p,q}=F^pH^{p+q}(X; \Q )/F^{p+1}H^{p+q}(X;\Q )$, 
we can easily check that $\chi (X)=\chi ((E_X)_{\infty})$. This proves the corollary. 
\end{proof}
\begin{rem}
In the case of a locally standard torus action ${\mathcal S}^{(0)}B_X$ corresponds to the fixed point set. In this case Corollary~\ref{fixed point theorem} recovers the well-known fact that $\chi (X)$ is equal to the number of the fixed points. For example, see \cite{Kawakubo} for this fact. 
\end{rem}

\begin{ex}
The orbit map of Example~\ref{ex3.9} is equipped with a section.  Let us compute the cohomology groups of $X$ in Example~\ref{ex3.9}. By cutting $B_X$ along curves $\alpha$ and $\beta$, a cell decomposition of $B_X$ is given. Figure~\ref{fig6} is the development, in which one-cells $e^{(1)}_1$, $e^{(1)}_2$, and $e^{(1)}_3$ correspond to $\alpha$, $\beta$, and the edge arc $\gamma$, respectively. 
\begin{figure}[hbtp]
\begin{center}
\input{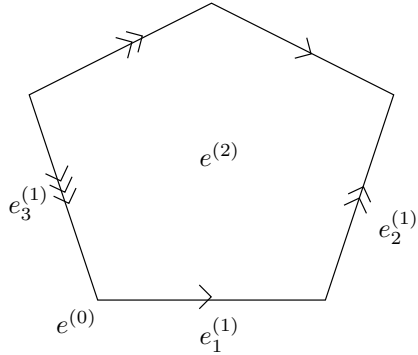}
\caption{A cell decomposition of $B_X$}
\label{fig6}
\end{center}
\end{figure}
Let $\overline{B_X}$ be the pentagon in Figure~\ref{fig6}. $\pi_{T_X}\colon T_X\to B_X$ is obtained from $\overline{B_X}\times T^2$ in the usual manner, namely, by identifying the preimage of $e^{(1)}_1$ (resp. $e^{(1)}_2$) by the first projection of $\overline{B_X}\times T^2$ with the preimage of the corresponding edge of $\overline{B_X}$ fiberwise by using $\rho ([\alpha])$ (resp. $\rho ([\beta])$). Let $\overline{X}$ be the quotient space 
\[
\overline{X}:=\overline{B_X}\times T^2 /\sim ,
\]
where $(b_1, u_1)\sim (b_2, u_2)$ if $b_1=b_2\in e^{(1)}_3$ and $u_2u_1^{-1}\in \{ 1\}\times S^1$, or 
$b_1=b_2$ and it is a vertex of $\overline{B_X}$. $X$ can be obtained from $\overline{X}$ by the same way as $\pi_{T_X}\colon T_X\to B_X$. The natural map $\overline{B_X}\times T^2\to \overline{B_X}\times T^2 /\sim $ descends to 
the map $\nu\colon T_X\to X$. For all cells except for $e^{(0)}$ and $e^{(1)}_3$, $\nu_{c^{(p)}_{\lambda}}\colon \pi_{T_X}^{-1}(c^{(p)}_{\lambda})(=T^2)\to {\mu_X}^{-1}(c^{(p)}_{\lambda})$ is a homeomorphism, whereas for $e^{(0)}$ and $e^{(1)}_3$ it can be identified with the natural projections $T^2\to T^2/T^2\cong \{ 1\}$ and $T^2\to T^2/\{ 1\}\times S^1\cong S^1\times \{ 1\}$, respectively.  

For $q=0$, $H^0(\pi_{T_X}^{-1}(c^{(p)}_{\lambda});\Z)=\Z$ and $\nu_{c^{(p)}_{\lambda}}^*: H^0({\mu_X}^{-1}(c^{(p)}_{\lambda}); \Z )\to H^0(\pi_T^{-1}(c^{(p)}_{\lambda}); \Z )$ is an isomorphism for all cells. Then $H^p(B_X ; \SH_X^0)$ is identified with $H^p(B_X ; \SH_T^0)$, which is naturally isomorphic to $H^p(B_X ; \Z )$. 

For $q=1$, $H^1(\pi_{T_X}^{-1}(c^{(p)}_{\lambda});\Z)=\Z\oplus \Z$ and $\nu_{c^{(p)}_{\lambda}}^*: H^1({\mu_X}^{-1}(c^{(p)}_{\lambda}); \Z )\to H^1(\pi_T^{-1}(c^{(p)}_{\lambda})$ $; \Z )$ is an isomorphism for all cells except for $e^{(0)}$ and $e^{(1)}_3$, whereas $\nu_{c^{(0)}}^*$ is the zero map and $\nu_{c^{(1)}_3}^*$ can be identified with the natural inclusion $\Z\oplus \{ 0\}\to \Z\oplus \Z$. Then $C^0(B_X;\SH_X^1)=0$, $C^1(B_X;\SH_X^1)\cong \Z^{\oplus 5}$, and $C^2(B_X;\SH_X^1)\cong \Z\oplus \Z$ and all coboundary operators vanish except for $\delta^1\colon C^1(B_X;\SH_X^1)\to C^2(B_X;\SH_X^1)$. For $c\in C^1(B_X;\SH_X^1)$, $\delta^1c$ is written as 
\[
\begin{split}
\delta^1c(e^{(2)})=&c(e^{(1)}_1)+ ^t\rho ([\alpha ])^{-1}c(e^{(1)}_2)- ^t\rho ([\alpha ][\beta ][\alpha ]^{-1})^{-1}c(e^{(1)}_1) \\
&- ^t\rho ([\alpha ][\beta ][\alpha ]^{-1}[\beta ]^{-1})^{-1}c(e^{(1)}_2) + ^t\rho ([\alpha ][\beta ][\alpha ]^{-1}[\beta ]^{-1})^{-1}c(e^{(1)}_3). 
\end{split}
\]
With the identification $C^1(B_X;\SH_X^1)\cong \Z^{\oplus 5}$ and $C^2(B_X;\SH_X^1)\cong \Z\oplus \Z$, $\delta^1$ can be identified with the matrix 
\[
\begin{pmatrix}
-1 & 1 & -2 & 2 &3 \\
-1 & 1 & -1 & 1 & 1
\end{pmatrix}. 
\]
Then, the cohomology groups are calculated by 
\[
H^p(B_X ; \SH_X^1)=
\begin{cases}
\Z^{\oplus 3} & p=1 \\
0 & \text{otherwise. }
\end{cases}
\]

For $q=2$, $H^2(\pi_{T_X}^{-1}(c^{(p)}_{\lambda});\Z)=\Z$ and $\nu_{c^{(p)}_{\lambda}}^*: H^2({\mu_X}^{-1}(c^{(p)}_{\lambda}); \Z )\to H^2(\pi_T^{-1}(c^{(p)}_{\lambda}); \Z )$ is an isomorphism for all cells except for $e^{(0)}$ and $e^{(1)}_3$, whereas both of $\nu_{c^{(0)}}^*$ and $\nu_{c^{(1)_3}}^*$ are the zero map. Then $C^0(B_X;\SH_X^2)=0$, $C^1(B_X;\SH_X^2)\cong \Z^{\oplus 2}$, and $C^2(B_X;\SH_X^2)\cong \Z$. All coboundary operators vanish. It is clear except for $p=1$. Since all monodromies along $e^{(1)}_\lambda$ induce the identity map of the second cohomology group of the fiber of $\pi_{T_X}\colon T_X\to B_X$ and $c(e^{(1)}_3)=0$ for $c\in C^1(B_X;\SH_X^2)$, $\delta^1$ also vanishes. Then the cohomology groups are obtained by 
\[
H^p(B_X ; \SH_X^2)=
\begin{cases}
\Z^{\oplus 2} & p=1 \\
\Z & p=2 \\
0 & \text{otherwise. }
\end{cases}
\]
The table for the $E_2$-terms is in Table~\ref{table1}. 
\begin{table}[hbtp]
\begin{tabular}{c|ccccc}
$q$ & & & & & \\
 & $0$ & $\Z^{\oplus 2}$ & $\Z$ & & \\[1mm]
 & $0$ & $\Z^{\oplus 3}$ & $0$ & & \\[1mm]
 & $\Z$ & $\Z^{\oplus 2}$ & $0$ & & \\[1mm]
\hline\tabtopsp{1mm}%
 &  &  &  & & $p$ \\[1mm]
\end{tabular}
\vspace*{3mm}\\
\caption{The table of $(E_X)^{p,q}_2$-terms}
\label{table1}
\end{table}
In particular, the Leray spectral sequence is degenerate at $E^2$-term. Then the cohomology groups of $X$ are given by
\[
H^k(X;\Z )=
\begin{cases}
\Z & k=0, 4 \\
\Z^{\oplus 2} & k=1, 3 \\
\Z^{\oplus 3} & k=2 \\
0  & \text{otherwise. } 
\end{cases}
\]
\end{ex}

\subsection{$K$-groups}\label{K-theory}
By replacing the cohomology functor $H^*(\ )$ by the $K$-functor $K^*(\ )$ in the cohomology Leray spectral sequence of the map ${\mu_X} :X\to B_X$, the similar method is available for computing $K$-groups. Such a spectral sequence is called the Atiyah-Hirzebruch spectral sequence for $K$-groups. We also denote this spectral sequence by the same notation $\{ (E_X)^{p,q}_r, d^X_r \}$. In this subsection, the result is described without proof. But the proof is almost same as in the last subsection. For $K$-theory, see \cite{At2} and for the spectral sequence of $K$-theory, see \cite{AH}. We still assume that $B_X$ is equipped with a structure of a CW complex so that each $p$-cell $e^{(p)}$ is contained in some ${\mathcal S}^{(k)}B_X$ of $B_X$. 
Let $\left( C^p(B_X ; \SK_T^q), \delta \right)$ be the cochain complex of the CW complex $B_X$ with the local system $\SK_T^q$ with respect to the $q$th $K$-group of the fiber of the fiber bundle $\pi_{T_X}: T_X\to B_X$. We denote by $C^p(B_X ; \SK_X^q)$ the subset of $C^p(B_X ; \SK_T^q)$ whose cochain takes a value in the image $\nu_{c^{(p)}}^*\left( K^q ({\mu_X}^{-1}(c^{(p)}))\right)$ of $K^q({\mu_X}^{-1}(c^{(p)}))$ by $\nu_{c^{(p)}}^*$ for each $p$-cell $e^{(p)}$. 
\begin{thm}\label{K-spectral}
For any $q$, $C^*(B_X ; \SK_X^q)$ is a subcomplex of $\left( C^*(B_X ; \SK_T^q), \delta \right)$. 
We denote its $p$th cohomology group by $H^p(B_X ;\SK^q_X)$. Then we have the isomorphisms 
\begin{align}
(E_X)_1^{p, q}&\cong C^p(B_X ; \SK^q_X), \ \ \ 
(E_X)_2^{p, q}\cong H^p(B_X ; \SK^q_X), \nonumber \\ 
(E_X)_{\infty}^{p,q} &=F^pK^{p+q}(X)/F^{p+1}K^{p+q}(X), \nonumber 
\end{align}
where $F^pK^*(X)$ is the image of the map $K^*(X, X^{(p-1)})\to K^*(X)$. 
\end{thm}

\begin{ex}
Let us compute $K$-groups of $X$ in Example~\ref{ex3.9}. We also use the cell decomposition of $B_X$ in Example~\ref{ex3.9}. For even $q$, since a fiber of $\pi_{T_X}\colon T_X\to B_X$ is $T^2$, by Lemma~\ref{K^0}, the $K$-group of its fiber is isomorphic to $\Z \oplus \Z$, and all homomorphisms between them which are induced from monodromies are identity. Also by Lemma~\ref{K^0}, $\nu^*_{c^{(p)}_{\lambda}}:K^q({\mu_X}^{-1}(c^{(p)}_{\lambda}))\to 
K^q(\pi_{T_X}^{-1}(c^{(p)}_{\lambda}))$ is isomorphic for all cells except for $e^{(0)}$ and $e^{(1)}_3$, whereas both of $\nu_{c^{(0)}}^*$ and $\nu_{c^{(1)}_3}^*$ are the natural inclusion $\Z\oplus \{ 0\} \to \Z\oplus \Z$. A similar computation as in Example~\ref{ex3.9} shows that 
\[
H^p(B_X;\SK^q_X)=
\begin{cases}
\Z & p=0, 2 \\
\Z^{\oplus  4} & p=1 \\
0 & \text{otherwise. }
\end{cases}
\]

For odd $q$, by Lemma~\ref{K^1}, we have an isomorphism 
\[
H^p(B_X;\SK^q_X)\cong H^p(B_X;\SH^1_X). 
\] 
$H^p(B_X;\SH^1_X)$ was computed in Example~\ref{ex3.9}. The table of $E^2$-terms is in Table~\ref{table2}.  
\begin{table}[hbtp]
\begin{tabular}{c|ccccc}
$q$ & & & & & \\
 & $\Z$ & $\Z^{\oplus 4}$ & $\Z$ & & \\[1mm]
 & $0$ & $\Z^{\oplus 3}$ & $0$ & & \\[1mm]
 & $\Z$ & $\Z^{\oplus 4}$ & $\Z$ & & \\[1mm]
\hline\tabtopsp{1mm}%
 & $0$ & $\Z^{\oplus 3}$ & $0$ & & $p$ \\[1mm]
 & $\Z$ & $\Z^{\oplus 4}$ & $\Z$ & &
\end{tabular}
\vspace*{3mm}\\
\caption{The table of $(E_X)^{p,q}_2$-terms}
\label{table2}
\end{table}
In particular, the spectral sequence is degenerate at $E^2$-term, and  $K$-groups of $X$ are given by
\[
K^k(X)=
\begin{cases}
\Z^{\oplus 5} & k: \text{even} \\
\Z^{\oplus  4} & k: \text{odd. }
\end{cases}
\]
\end{ex}

\subsection{On signatures in the oriented four-dimensional case}
In this subsection, we shall give the method of computing the signature for the four-dimensional case by using the Novikov additivity. We assume that both of $X$ and the interior of $B_X$ are oriented so that a weakly standard atlas of $X$ used in this subsection and the atlas of $B_X$ induced by it are compatible with the given orientations. In this subsection the assumption $e(X,\CT )=0$ is not necessary. 

For simplicity, suppose that $B_X$ has only one boundary component with $\CS^{(0)}B_X\neq \emptyset$. We divide $B_X$ into two parts $(B_X)_1$ and $(B_X)_2$, where $(B_X)_2$ is the closed neighborhood of the boundary $\partial {B_X}$ such that $\partial {B_X}$ is a deformation retract of $(B_X)_2$ and $(B_X)_1$ is the closure $(B_X)_1=({B_X}\backslash (B_X)_2)^{cl}$ of 
the remainder. We put $X_i={\mu_X}^{-1}((B_X)_i)$ for $i=1,2$, and denote by $\sigma (X_i)$ and $\sigma (X)$ the signature of $X_i$ and $X$, respectively. The Novikov additivity says that
\begin{equation}\label{signX}
\sigma (X)=\sigma (X_1)+\sigma (X_2). 
\end{equation}
First let us compute $\sigma (X_1)$. As we showed in Proposition~\ref{T^n-bundle}, $X_1$ is a $T^2$-bundle. If the genus of ${B_X}$ is equal to zero, then $(B_X)_1$ is contractible. In this case, $\sigma (X_1)$ is zero. 

\begin{figure}[hbtp]
\begin{center}
\input{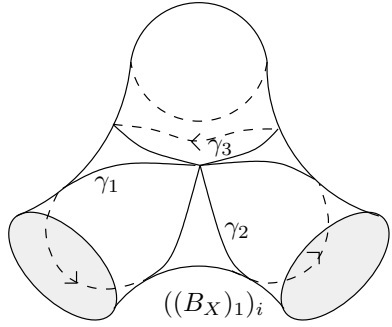}
\caption{$((B_X)_1)_i$ and $\gamma_j$}
\label{fig13}
\end{center}
\end{figure}

Suppose that the genus of ${B_X}$ is greater than zero. We give $(B_X)_1$ a trinion decomposition $\displaystyle (B_X)_1=\cup_{i=1}^k((B_X)_1)_i$, where each $((B_X)_1)_i$ is a surface obtained from $S^2$ by removing three distinct open discs. Let $(X_1)_i={\mu_X}^{-1}(((B_X)_1)_i)$ for $i=1, \ldots ,k$. From the Novikov additivity, we have 
\begin{equation}\label{signX_1}
\sigma (X_1)=\sum_{i=1}^k\sigma ((X_1)_i) ,  
\end{equation}
and each $\sigma ((X_1)_i)$ can be computed as follows. We take oriented loops $\gamma_1$, $\gamma_2$, and $\gamma_3$ of $((B_X)_1)_i$ as in Figure~\ref{fig13} which represent generators of $\pi_1(((B_X)_1)_i)$ with 
$[\gamma_1][\gamma_2][\gamma_3]=1$. Let $\rho :\pi_1(((B_X)_1)_i)\to Sp(2;\Z )$ be the monodromy representation of $T^2$-bundle ${\mu_X} \colon (X_1)_i\to ((B_X)_1)_i$. We put $C_j:=\rho ([\gamma_j])$ for 
$j=1, 2, 3$. 
For $C_1$ and $C_2$, define the vector space $V_{C_1, C_2}$ and the bilinear form $\< \ ,\ \>_{C_1, C_2}$ on $V_{C_1, C_2}$ by 
\[
\begin{split}
&V_{C_1, C_2}:=\{ (x, y)\in \R^2\times \R^2 \colon (C_1^{-1}-I)x+(C_2-I)y=0\} , \\
&\< (x, y), (x', y')\>_{C_1, C_2}:= ^t(x+y)J(I-C_2)y'
\end{split}
\]
for $(x, y)$, $(x', y')\in V_{C_1, C_2}$, where $I=
\begin{pmatrix}
1 & 0 \\
0 & 1
\end{pmatrix}$ and 
$J=
\begin{pmatrix}
0 & 1 \\
-1 & 0
\end{pmatrix}$. 
It is easy to see that $\< \ ,\ \>_{C_1, C_2}$ is symmetric and we denote the signature 
of $\< \ ,\ \>_{C_1, C_2}$ by $\tau_1(C_1, C_2)$. 
\begin{thm}[\cite{E, Me}]\label{sign(X_1)_i}
$\sigma ((X_1)_i)=\tau_1(C_1, C_2)$. 
\end{thm}
We should remark that our orientation of $(X_1)_i$ is different from that in \cite{E, Me}. 
From \eqref{signX_1} and Theorem~\ref{sign(X_1)_i}, we can compute $\sigma (X_1)$.  
\begin{rem}
Meyer shows in \cite{Me} that $\tau_1$ defines the cocycle $\tau_1:Sp(2;\Z )\times Sp(2;\Z )\to \Z$ of $Sp(2;\Z )$. It is called Meyer's signature cocycle. The author was taught Meyer's signature cocycle by Endo~\cite{E}. 
\end{rem}

Next, we shall compute $\sigma (X_2)$. For $X_2$, we can show the following lemma.  
\begin{lem}
${\mu_X}^{-1}(\partial {B_X})$ is a deformation retract of $X_2$.  
\end{lem}
\begin{proof}
Let $h_s\colon (B_X)_2\to (B_X)_2$ $(s\in I)$ be a deformation retraction with $h_0=\id_{(B_X)_2}$, $h_1((B_X)_2)=\partial {B_X}$, and $h_s|_{\partial {B_X}}=\id_{\partial {B_X}}$. By using the homotopy lifting property of the restriction of $\pi_{T_X}\colon T_X\to {B_X}$ to $(B_X)_2$, there exists a lift $\tilde{h}_s^T\colon T_X|_{(B_X)_2}\to T_X|_{(B_X)_2}$ of $h_s$ such that 
$\tilde{h}_0^T=\id_{T_X|_{(B_X)_2}}$ and $\tilde{h}_s^T|_{(T_X|_{\partial {B_X}})}=\id_{T_X|_{\partial {B_X}}}$. $\tilde{h}_s^T$ induces the map $\tilde{h}_s^X\colon X_2\to X_2$ such that $\tilde{h}_0^X=\id_{X_2}$, $\tilde{h}_s^X|_{{\mu_X}^{-1}(\partial {B_X})}=\id_{{\mu_X}^{-1}(\partial {B_X})}$, and the following diagram commutes:
\[
\xymatrix{
 &  &  T_X|_{(B_X)_2}\ar[rr]^{\tilde{h}_s^T}\ar[dll]_{\nu}\ar[ddl]|\hole_{\pi_{T_X}} &  & 
 T_X|_{(B_X)_2}\ar[dll]_{\nu}\ar[ddl]_{\pi_{T_X}} \\
X_2\ar[rr]^{\tilde{h}_s^X\hspace*{3mm}}\ar[dr]_{{\mu_X}} & & X_2\ar[dr]_{{\mu_X}} & & \\
 & (B_X)_2\ar[rr]^{h_s} & & (B_X)_2.  
}
\]
In particular, $\tilde{h}_s^X$ is a required deformation retraction. 
\end{proof}

Suppose that the cardinality of $\CS^{(0)}{B_X}$ is equal to $k$. Then $\CS^{(1)}{B_X}$ has exactly $k$ connected components $(\CS^{(1)}{B_X})_1$, $\ldots$, $(\CS^{(1)}{B_X})_k$, and $\displaystyle {\mu_X}^{-1}(\partial {B_X})
=\cup_{i=1}^k{\mu_X}^{-1}((\CS^{(1)}{B_X})_i^{cl})$, where $(\CS^{(1)}{B_X})_i^{cl}$ is the closure of $(\CS^{(1)}{B_X})_i$. By the construction of the canonical model, it is easy to see that each ${\mu_X}^{-1}((\CS^{(1)}{B_X})_i^{cl})$ is homeomorphic to the two-dimensional sphere $S^2$ if $k\ge 2$, and is homeomorphic to $S^2$ with one self-intersection at north and south points if $k=1$. If $(\CS^{(1)}{B_X})_i^{cl}\cap (\CS^{(1)}{B_X})_j^{cl}\neq \emptyset$ for $i\neq j$, then they have two intersections if $k=2$, and have one intersection if $k>2$. In all cases, every intersection is transversal since a neighborhood of an intersection in $X$ is identified with that of the intersection of $\C\times \{ 0\}$ and $\{ 0\}\times \C$ in $\C^2$. Then ${\mu_X}^{-1}(\partial {B_X})$ looks like a necklace consisting of $k$ spheres and the homology group of $X_2$ is given by 
\[
H_p(X_2;\Z )=H_p({\mu_X}^{-1}(\partial {B_X}); \Z )=
\begin{cases}
\Z & p=0,1 \\
\Z^{\oplus k} & p=2 \\
0 & \text{otherwise.}
\end{cases}
\]
Moreover, by putting $S^2_i:={\mu_X}^{-1}((\CS^{(1)}{B_X})_i^{cl})$ for $i=1, \ldots k$, the homology classes $[S^2_i]\in H_2(X_2;\Z )$ represented by $S^2_i$ are generators of $H_2(X_2;\Z )$. As a summary of the above observation, we obtain the following proposition. 
\begin{prop}\label{int}
For $i\neq j$, the intersection number $[S^2_i]\cdot [S^2_j]$ of $[S^2_i]$ and $[S^2_j]$ is given as follows 
\[
[S^2_i]\cdot [S^2_j]=
\begin{cases}
0 & S^2_i\cap S^2_j=\emptyset \\
1 & S^2_i\cap S^2_j\neq \emptyset \ \text{and}\ k>2 \\
2 & S^2_i\cap S^2_j\neq \emptyset \ \text{and}\ k=2. 
\end{cases}
\]
\end{prop}
Assume that $k>1$. Let us compute the self-intersection number of $[S^2_i]$. We can take a contractible neighborhood $U$ of $(\CS^{(1)}{B_X})_i^{cl}$ in ${B_X}$ so that $U\cap \CS^{(1)}{B_X}$ has exactly two connected components, say $(U\cap \CS^{(1)}{B_X})_1$ and $(U\cap \CS^{(1)}{B_X})_2$, except for $(\CS^{(1)}{B_X})_i$. We may assume that 
$(U\cap \CS^{(1)}{B_X})_1$ and $(U\cap \CS^{(1)}{B_X})_2$ are located as in Figure~\ref{fig12}.    
\begin{figure}[hbtp]
\begin{center}
\input{sigma18.pstex_t}
\caption{$(\CS^{(1)}{B_X})_i$ and $(U\cap \CS^{(1)}{B_X})_a$}
\label{fig12}
\end{center}
\end{figure}
Since $U$ is contractible, there exists a local trivialization $\varphi^\Lambda \colon \pi_{\Lambda_X}^{-1}(U)\cong U\times \Lambda$ of $\pi_{\Lambda}\colon \Lambda_X\to {B_X}$ which sends restrictions of $\pi_{\CL_X}\colon \CL_X \to \CS^{(1)}{B_X}$ to $(\CS^{(1)}{B_X})_i$ and $(U\cap \CS^{(1)}{B_X})_a$ ($a=1$, $2$) fiberwise to the trivial rank one subbundles $(\CS^{(1)}{B_X})_i\times L$ and $(U\cap \CS^{(1)}{B_X})_a\times L_a$, respectively. Since $\pi_{\CL_X}\colon \CL_X \to \CS^{(1)}{B_X}$ is primitive, we can take generators $v$, $v_a$ of $L$ and $L_a$ such that both of the determinants of $(v_1, v)$ and $(v, v_2)$ are equal to one, where $(v_1, v)$ (resp. $(v, v_2)$) denotes the matrix given by arranging the column vectors $v_1$ and $v$ (resp. $v$ and $v_2$) in this order. 
\begin{prop}\label{selfint}
The self-intersection number $[S^2_i]\cdot [S^2_i]$ is equal to the negative determinant $-\det (v_1, v_2)$ of $(v_1, v_2)$. 
\end{prop}
Note that the determinant of $(v_1, v_2)$ does not depend on the choice of the local trivialization $\varphi^\Lambda$ since the structure group of the bundle preserves the orientation of a fiber.  
\begin{proof}
Since the self-intersection number of $S^2_i$ is equal to the Euler number of the normal bundle $\CN_{S^2_i}$ of $S^2_i$ in $X$ (for example, see \cite{BT}), we identify $\CN_{S^2_i}$. For a positive number $\varepsilon$, let $D_{\varepsilon}$ be the two-dimensional closed disc $D_{\varepsilon}:=\{ z\in \C \colon \abs{z}^2\le \varepsilon \}$ 
and $D_{\varepsilon}^{int}$ its interior $D_{\varepsilon}^{int}:=\{ z\in \C \colon \abs{z}^2<\varepsilon \}$. 
We divide $U$ into $U_1$ and $U_2$ as in Figure~\ref{fig26}. 
\begin{figure}[hbtp]
\begin{center}
\input{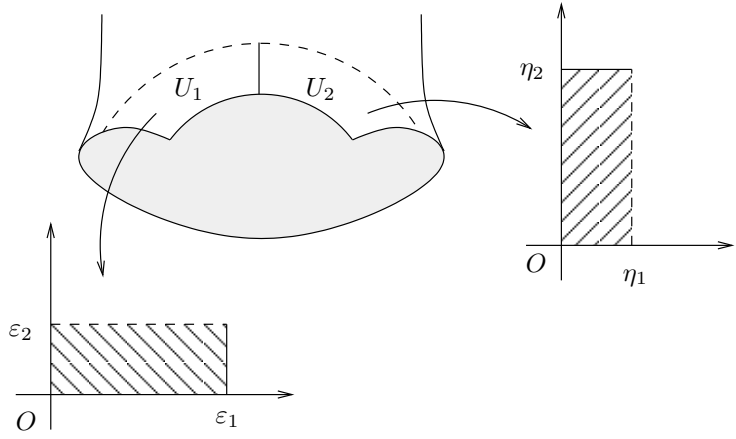}
\caption{$U_1$ and $U_2$}
\label{fig26}
\end{center}
\end{figure}
We may assume that there are homeomorphisms $\varphi^X_1\colon {\mu_X}^{-1}(U_1)\to D_{\varepsilon_1}\times D_{\varepsilon_2}^{int}$, $\varphi^X_2\colon {\mu_X}^{-1}(U_2)\to D_{\eta_1}^{int}\times D_{\eta_2}$, $\varphi^B_1\colon U_1\to [0,\varepsilon_1]\times [0,\varepsilon_2 )$, and $\varphi^B_2\colon U_2\to [0,\eta_1 )\times [0,\eta_2 ]$ such that the overlap map $\varphi^X_{12}:=\varphi^X_1\circ (\varphi^X_2)^{-1}\colon D_{\eta_1}^{int}\times \partial D_{\eta_2}\to \partial D_{\varepsilon_1}\times D_{\varepsilon_2}^{int}$ is $\rho_{12}$-equivariant (with respect to the standard $T^2$-actions) for some positive numbers $\varepsilon_i$, $\eta_i$ ($i=1$, $2$) and $\rho_{12}\in \Aut (T^2)$ and overlap maps satisfy the equation $\varphi^B_{12}\circ \mu_{\C^2}=\mu_{\C^2}\circ\varphi^X_{12}$. Moreover, by applying the same argument used in Proposition~\ref{section} if necessary, we may also assume that overlap maps satisfy the equation $\iota \circ \varphi^B_{12}=\varphi^X_{12}\circ \iota$, where $\iota$ is the section of $\mu_{\C^2}$ defined by \eqref{iota}. Since $\varphi^X_{12}$ sends $\{ 0\} \times \partial D_{\eta_2}$ to $\partial D_{\varepsilon_1}\times \{ 0\}$ and stabilizers of the restriction of the standard $T^2$-action to $\{ 0\} \times \partial D_{\eta_2}$ and $\partial D_{\varepsilon_1}\times \{ 0\}$ are $S^1\times \{ 1\}$ and $\{ 1\}\times S^1$, respectively, $\rho_{12}$ must be of the form 
\begin{equation}\label{rho_12}
\rho_{12}=
\begin{pmatrix}
0 & -1 \\
1 & m
\end{pmatrix}
\end{equation}
for some integer $m$. We show that the Euler number of the normal bundle $\CN_{S^2_i}$ of $S^2_i$ in $X$ is equal to $-m$. In this setting, $S^2_i$ is obtained by gluing $\{ 0\} \times D_{\eta_2}$ and $D_{\varepsilon_1}\times \{ 0\}$ with the overlap map $\varphi^X_{12}\colon \{ 0\}\times \partial D_{\eta_2}\to \partial D_{\varepsilon_1}\times \{ 0\}$. Under the natural identifications of both $\{ 0\}\times D_{\eta_2}$ and $D_{\varepsilon_1}\times \{ 0\}$ with the unit disc $D_1$, the direct computation using \eqref{rho_12} and the fact that $\varphi^X_{12}$ is $\rho_{12}$-equivariant shows that $\varphi^X_{12}\colon \{ 0\}\times \partial D_{\eta_2}\to \partial D_{\varepsilon_1}\times \{ 0\}$ is of the form 
\[
\varphi^X_{12}(z)=z^{-1}
\]
for $z\in \partial D_1$ and the frame bundle of $\CN_{S^2_i}$ can be identified with the $S^1$-bundle which is obtained by gluing two copies of $D_1\times S^1$ with the map sending $(z, u)\in \partial D_1\times S^1$ to $(z^{-1}, z^mu)\in \partial D_1\times S^1$. This implies that the Euler number of $\CN_{S^2_i}$ is equal to $-m$. 
\end{proof}

From Proposition~\ref{int} and \ref{selfint}, we can compute $\sigma (X_2)$ case-by-case for $k>1$. In case of $k=1$, we have a unique point $b$ in $\CS^{(0)}{B_X}$ and the fiber ${\mu_X}^{-1}(b)$ also consists of exactly one point, say $x\in {\mu_X}^{-1}(b)$. By blowing up $X$ at $x$, we can reduce to the case of $k>1$ as in the following example. 

\begin{ex}
Let us compute the signature of Example~\ref{ex3.9}. Recall that ${B_X}$ is a surface of genus one with one corner. As in Figure~\ref{fig14}, we divide ${B_X}$ into two parts $(B_X)_1$ and $(B_X)_2$ and give $(B_X)_1$ the trinion decomposition. The easy computation shows that the value $\tau_1 (\rho ([\alpha^{-1}]), \rho ([\gamma^{-1}]))$ of the Meyer cocycle vanishes. This implies the signature $\sigma (X_1)$ of $X_1$ is zero. 
\begin{figure}[hbtp]
\begin{center}
\input{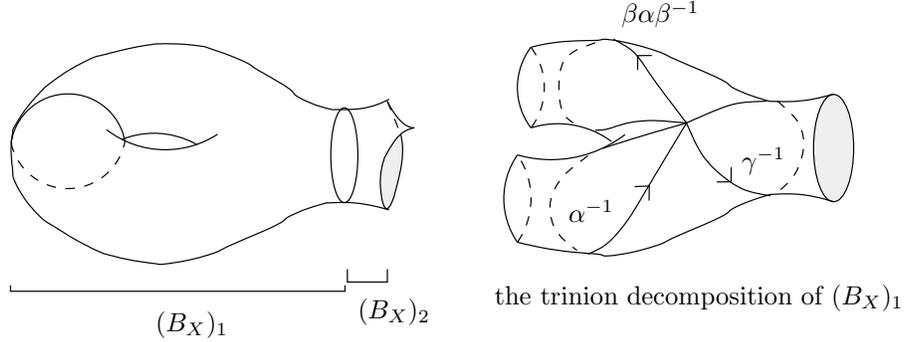}
\caption{${B_X}$, $(B_X)_i$, and the trinion decomposition of $(B_X)_1$}
\label{fig14}
\end{center}
\end{figure}

Next we focus on $X_2$. Since $k=1$, the fiber of ${\mu_X}\colon X\to {B_X}$ at the unique point in $\CS^{(0)}{B_X}$ consists of exactly one point which we denote by $x_0$. By definition, a sufficiently small neighborhood of $x_0$ is identified with that of the origin of $\C^2$. We denote the blowing-up of $X_2$ at $x_0$ by $\widetilde{X_2}$ and also denote the corresponding orbit space by $\widetilde{(B_X)_2}$. (See Example~\ref{blowing-up} for blowing-ups.) Note that $\CS^{(0)}\widetilde{(B_X)_2}$ consists of exactly two points. Let $(\CS^{(1)}\widetilde{(B_X)_2})_1$ and $(\CS^{(1)}\widetilde{(B_X)_2})_2$ denote connected components of $\CS^{(1)}\widetilde{(B_X)_2}$ as in Figure~\ref{fig15}.  
\begin{figure}[hbtp]
\begin{center}
\hspace*{-1.5cm}
\input{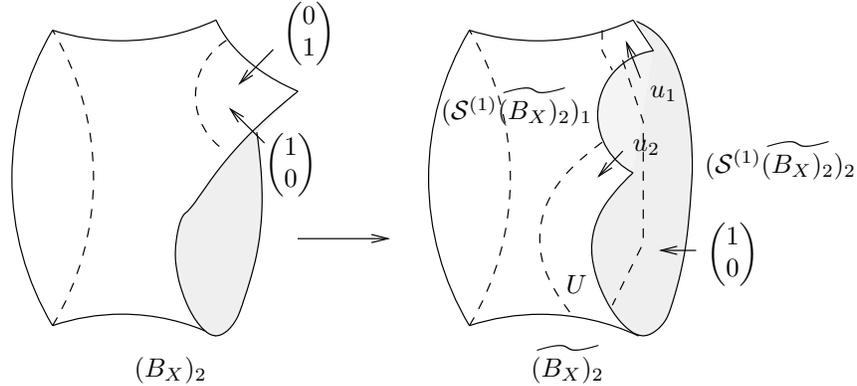}
\caption{The blowing up of $(B_X)_2$}
\label{fig15}
\end{center}
\end{figure}
Since the preimage $S^2_1={\mu_X}^{-1}((\CS^{(1)}\widetilde{(B_X)_2})_1^{cl})$ is an exceptional divisor, its self-intersection number $[S^2_1]\cdot [S^2_1]$ is equal to $-1$. We compute the self-intersection number $[S^2_2]\cdot [S^2_2]$ of $S^2_2={\mu_X}^{-1}((\CS^{(1)}\widetilde{(B_X)_2})_2^{cl})$. By the construction of $X_2$ in Example~\ref{ex3.9}, we can take $u_1$, $u_2$ in Proposition~\ref{selfint} of the forms 
\[
u_1=
\begin{pmatrix}
3 & 1 \\
-1 & 0
\end{pmatrix}
\begin{pmatrix}
1 \\
1
\end{pmatrix}
=
\begin{pmatrix}
4 \\
-1
\end{pmatrix}
,\ \  
u_2=
\begin{pmatrix}
1 \\
1
\end{pmatrix}, 
\]
hence $[S^2_2]\cdot [S^2_2]=-\det (u_1, u_2)=-5$. 
The above computation and Proposition~\ref{int} for $k=2$ show that the intersection 
matrix of $\widetilde{X_2}$ is 
\[
\begin{pmatrix}
-1 & 2 \\
2 & -5
\end{pmatrix}
\]
and $\sigma (\widetilde{X_2})$ is equal to $-2$. Since $\widetilde{X_2}$ is the blowing-up of $X_2$ at $x_0$, $\sigma (X_2)$ is equal to $-1$. Then the computations of $\sigma (X_1)$ and $\sigma (X_2)$ together with $\eqref{signX}$ shows that $\sigma (X)=-1$. 
\end{ex}


\appendix
\section{$K$-theory of low dimensional finite CW complexes}
In this appendix, we shall give some lemmas for $K$-theory of low dimensional finite CW complexes 
which are used in Section~\ref{K-theory}. 
Let $X$ be a finite CW complex. 
\begin{lem}[\cite{Fur}]\label{K^0}
If $\dim X\le 3$, $K(X)$ is isomorphic to $H^0(X;\Z )\oplus H^2(X;\Z )$, where $H^0(X;\Z )\oplus H^2(X;\Z )$ is considered as a subring of the cohomology ring of $X$ with $\Z$-coefficients. 
\end{lem}
\begin{proof}
We denote $H^0(X;\Z )\oplus H^2(X;\Z )$ by $H^0\oplus H^2$. For simplicity we assume that $X$ is connected. Let $\Vect (X)$ be the set of isomorphism classes of vector bundles on $X$. $\Vect (X)$ is a semi-ring with respect to a direct sum and a tensor product. We define the semi-ring homomorphism from $\Vect (X)$ to $H^0\oplus H^2$ by assigning the rank and the first Chern class to an element of $\Vect (X)$. By universality of $K(X)$ it induces a ring homomorphism from $K(X)$ to $H^0\oplus H^2$. We show that this is isomorphism. For any $(r, \alpha )\in H^0\oplus H^2$ there exists a line bundle $L$ on $X$ whose first Chern class is equal to $\alpha$ since $H^2$ is identified with the set of isomorphism classes of line bundles on $X$. Then the above map sends the isomorphism class of $L\oplus \underline{\C}^{r-1}$ to $(r, \alpha )$, where $\underline{\C}^{r-1}$ is the trivial rank $(r-1)$ vector bundle. This proves that it is surjective. Next we show that it is injective. Suppose that an element $[E]-[F]\in K(X)$ lies in the kernel of this map. By construction they have the same rank, say $r$, and the same first Chern class. Then there exists an isomorphism $\varphi \colon \wedge^rE\to \wedge^rF$ between the top exterior bundles because they are line bundles whose first Chern classes are same as that of $E$ and $F$, respectively. For $x\in X$ define $Q_x$ by 
\[
Q_x:=\{ \text{a linear isormophism }f\colon E_x\to F_x\ \colon \wedge^rf=\varphi|_{\wedge^rE_x} \} . 
\]
It is easy to see that $Q_x$ is identified with $\SL_r(\C )$ and the natural projection $\pi \colon \coprod_xQ_x\to X$ is a fiber bundle with fiber $\SL_r(\C )$. To prove the lemma it is sufficient to show that $\pi \colon \coprod_xQ_x\to X$ has a section. We show this by induction on the number of cells contained in $X$. If $X$ consists of only one cell, it is trivial. Assume that it is valid up to $k-1$ and also assume that $X$ consists of exactly $k$ cells. Then $X$ is obtained from some CW complex $X'$ consisting of exactly $k-1$ cells by attaching a cell $D$. By the assumption of the induction, the restriction of $\pi \colon \coprod_xQ_x\to X$ to $X'$ is equipped with a section, and also by the assumption of the lemma, the dimension of $D$ is equal or less than $3$. Since the $i$th homotopy group $\pi_i(\SL_r(\C ))$ vanishes for $i\le 2$ the section can be extended to $X$. This prove the lemma. 
\end{proof}
\begin{cor}
Under the assumption of Lemma~\textup{\ref{K^0}}, by assigning the first Chern class to an element of 
$\widetilde{K}(X)$, we have the ring isomorphism  
\begin{equation}\label{reducedK}
\widetilde{K}(X)\cong H^2(X;\Z ). 
\end{equation}
\end{cor}
\begin{lem}\label{K^1}
If $\dim X\le 2$, there is a ring isomorphism 
\[
K^{-1}(X)\cong H^1(X;\Z ) 
\]
which commutes with maps induced from a continuous map $f:X\to Y$ 
\[
\xymatrix{
K^{-1}(X)\ar[r]^{\cong} & H^1(X;\Z )\ar@{}[dl]|\circlearrowright \\
K^{-1}(Y)\ar[r]^{\cong}\ar[u]^{f^*} & H^1(Y;\Z)\ar[u]_{f^*}. 
}
\]
\end{lem}
\begin{proof}
By definition, 
\[
K^{-1}(X):=\widetilde{K}(\widetilde{S}X^+), 
\]
where $X^+$ is the disjoint union of $X$ with a point $x_0$ and $\widetilde{S}X^+$ is the reduced 
suspension of $X^+$. Let $SX^+$ and $Sx_0$ be unreduced suspensions of $X^+$ and $x_0$, 
respectively. $Sx_0$ is naturally included in $SX^+$ and the quotient space obtained from $SX^+$ 
by collapsing $Sx_0$ is just $\widetilde{S}X^+$. Atiyah shows in \cite{At2} that the natural projection 
$p:SX^+\to \widetilde{S}X^+$ induces the isomorphism 
\begin{equation}\label{projection}
p^*:\widetilde{K}(\widetilde{S}X^+)\cong \widetilde{K}(SX^+). 
\end{equation}
By composing the isomorphisms \eqref{projection}, \eqref{reducedK} for $SX^+$, and 
the suspension isomorphism 
\[
H^2(SX^+;\Z )\cong H^1(X^+;\Z )(\cong H^1(X;\Z )), 
\]
we obtain the isomorphism in the lemma.  
\end{proof}

\begin{rem}
It is also shown in \cite{Fur} that when $\dim X\le 5$, by associating the rank, the first Chern class, 
and the second Chern class with an element of $K(X)$ (by associating the first and second Chern 
classes in the case of $\widetilde{K}(X)$), we have the one-to-one correspondences as a set
\[
\begin{split}
K(X)&\cong H^0(X;\Z )\oplus H^2(X;\Z )\oplus H^4(X;\Z ), \\
\widetilde{K}(X)&\cong H^2(X;\Z )\oplus H^4(X;\Z ). 
\end{split} 
\]
Then the same argument in the proof of Lemma~\ref{K^1} shows that when $\dim X\le 4$, 
there is an one-to-one correspondence 
\[
K^{-1}(X)\cong H^1(X;\Z ) \oplus H^3(X;\Z ) 
\]
as a set.  
\end{rem}

\bibliographystyle{amsplain}
\bibliography{toric-related}
\end{document}